%% file: swenwangnewversion4.tex
\documentclass[11pt,reqno]{amsart}

\usepackage{latexsym}
\usepackage{amsmath,amsthm}
\usepackage{epsfig,psfrag,verbatim,amsfonts,amssymb}

\usepackage{amssymb,amsmath}
\usepackage{hyperref}

\def\bl{\begin{lemma}}
\def\el{\end{lemma}}
\def\bth{\begin{theorem}}
\def\eth{\end{theorem}}
\def\bc{\begin{corollary}}
\def\ec{\end{corollary}}
\def\bcj{\begin{conjecture}}
\def\ecj{\end{conjecture}}
\def\bpr{\begin{proposition}}
\def\epr{\end{proposition}}
\def\bde{\begin{definition}}
\def\ede{\end{definition}}
\def\E{\mathbb{E}}
\def\Pr{\mathbb{P}}
\def\Var{\mbox{\rm Var}}

\newcommand{\be}{\begin{eqnarray}}
\newcommand{\ee}{\end{eqnarray}}
\newcommand{\eps}{{\mbox{$\epsilon$}}}
\newcommand{\R}{{\mathbb R}}

\newcommand{\Z}{{\mathbb Z}}

\newcommand{\F}{{\mathcal F}}

\newcommand{\A}{{\mathcal A}}
\newcommand{\B}{{\mathcal B}}

\newcommand{\EE}{{\mathcal E}}

\newcommand{\sign}{{\rm sign}}
\newcommand{\with}{\hbox{ {\rm with} }}
\renewcommand{\and}{\hbox{ {\rm and} }}

\newcommand{\Tm}{T_{{\rm mix}}}
\newcommand{\C}{{\mathcal{C}}}

\newcommand{\prob}{\mbox{\bf P}}
\renewcommand{\P}{\mbox{\bf P}}
\newcommand{\p}{\mbox{\bf p}}

\newtheorem{theorem}{Theorem}[section]
\newtheorem{definition}{Definition}[section]
\newtheorem{thm}[theorem]{Theorem}
\newtheorem{lemma}[theorem]{Lemma}
\newtheorem{lem}[theorem]{Lemma}
\newtheorem{prop}[theorem]{Proposition}

\newtheorem{corollary}[theorem]{Corollary}
\newtheorem{proposition}[theorem]{Proposition}
\newtheorem{conjecture}[theorem]{Conjecture}

\theoremstyle{definition}
\numberwithin{equation}{section}

\def\m{{m}}
\def\d{{\rm d}}
\def\half{\frac{1}{2}}
\newcommand{\1}{\mbox{\bf 1}}

\begin{document}

\title{A Power Law of Order 1/4 for Critical Mean Field Swendsen-Wang Dynamics}

\author{Yun Long}
\author{Asaf Nachmias}
\author{Weiyang Ning}
\author{Yuval Peres}

\begin{abstract}

The Swendsen-Wang dynamics is a Markov chain widely used by
physicists to sample from the Boltzmann-Gibbs distribution of the
Ising model. Cooper, Dyer, Frieze and Rue proved that on the
complete graph $K_n$ the mixing time of the chain is at most
$O(\sqrt{n})$ for all non-critical temperatures. In this paper we show that the mixing time is $\Theta(1)$ in high temperatures, $\Theta(\log n)$ in low temperatures and $\Theta(n^{1/4})$ at criticality. We also
provide an upper bound of $O(\log n)$ for Swendsen-Wang dynamics for the
$q$-state ferromagnetic Potts model on any tree of $n$ vertices.
\end{abstract}
\maketitle

\tableofcontents

\section{Introduction}

Local Markov chains (e.g. ``Glauber dynamics'') are commonly used to
sample from spin systems on graphs. At low temperatures, however,
their mixing time becomes very large (sometimes exponential in the
size of the graph), making it computationally harder to sample from
the equilibrium measure. In some cases, ``Global'' Markov chains,
which allow moves like cluster flipping, yield much faster mixing
and those are the algorithms of choice when practitioners actually
sample (see for example \cite{S}, \cite{SS1}, \cite{SS2} and
\cite{W}; see \cite{JS} for a different polynomial time algorithm
for sampling from the Ising model). The Swendsen-Wang (SW) algorithm
for the $q$-state ferromagnetic Potts model and its variations are
frequently used in practice.
 Gore and Jerrum \cite{GJ} discovered that for any
$q>2$, on the complete graph $K_n$ there are temperatures where the
SW dynamics has mixing time of order at least
$\exp(\Omega(\sqrt{n}))$.  Borgs, Chayes, Frieze, Kim, Tetali,
Vigoda and Vu \cite{BCFKTVV} proved a similar lower bound on the
mixing time of the SW algorithm at the critical temperatures, on the
$d$-dimensional lattice torus for any $d\geq 2$ and $q$ sufficiently
large.

The natural question remaining is how does the SW algorithm perform
when $q=2$ (i.e., for the Ising model). The first
positive result in this direction is due to Cooper, Dyer, Frieze and
Rue \cite{CDFR}. They proved that the SW algorithm on the complete
graph on $n$ vertices has mixing time at most $O(\sqrt{n})$ for all
non-critical temperatures. In this paper we show that the mixing time is $\Theta(1)$ in high temperatures, $\Theta(\log n)$ in low temperatures and $\Theta(n^{1/4})$ at the critical temperature. The study of the mixing time at criticality is the main effort of this paper. Heuristic arguments for the exponent $1/4$ at criticality
were found  earlier by physicists, see \cite{RTK} and \cite{PBKD}.

It is instructive to compare these results with the mixing time of Glauber dynamics for the critical Ising model on $K_n$. In \cite{LLP}, the authors show that this mixing time is $\Theta(n^{3/2})$. Since in the SW dynamics we update $n$ vertices in each step, the number of vertex updates up to the mixing time is $\Theta(n^{5/4})$, that is, it performs faster by $n^{1/4}$ than Glauber dynamics.

\section{Statement of the results}

The mixing time of a finite Markov chain with transition matrix $\p$
is defined by
$$\Tm = \Tm(1/4) = \min \Big \{ t : \| \p^t(x,\cdot) - \pi(\cdot) \|_{{\rm TV}} \le 1/4 \, , \hbox{{\rm  for all }} x \in V \Big \} \, ,$$
where  $\| \mu - \nu \|_{{\rm TV}} = \max_{A \subset V} |\mu(A) -
\nu(A)|$ is the total variation distance. Before describing the Swendsen-Wang (SW) algorithm on a graph $G=(V,\EE)$, let us first describe its stationary distribution, also known as the Ising model. This is a probability measure (also known as the
Boltzmann-Gibbs distribution) on the set $\Omega = \{1, -1 \}^{V}$ where the probability of each $\sigma \in \Omega$ is
$$ \P(\sigma) = { e^{\beta \sum_{(u,v)\in \EE} \sigma(u) \sigma(v)} \over Z(G) } \, ,$$
where $\beta \in [0,\infty]$ is a parameter usually referred to as the {\em inverse temperature}, and the {\em partition} function $Z(G)$ is defined by
$$ Z(G) = \sum _{\sigma \in \Omega} e^{\beta \sum_{(u,v)\in \EE} \sigma(u) \sigma(v)} \, .$$
For $\sigma \in \Omega$, let $G^+(\sigma)$  be the graph spanned by
the vertices of $G$ which are assigned $1$ by $\sigma$ and similarly
let $G^-(\sigma)$ be the graph spanned by the vertices of $G$ which
are assigned $-1$ by $\sigma$. The SW dynamics on a graph $G$ with
percolation parameter $p\in[0,1]$ is a Markov chain on $\Omega$.
Given the current state of the chain $\sigma_t$, we obtain the next
state $\sigma_{t+1}$ by the following two-step procedure:

\begin{enumerate}
\item[1.] Perform independent $p$-bond percolation on
$G^+(\sigma_t)$ and on $G^-(\sigma_t)$ separately. That is, retain each edge of $G^+(\sigma_t)$ and $G^-(\sigma_t)$ with probability $p$ and erase with probability $1-p$, independently for all edges. Call the obtained graphs $G^+_p$ and $G^-_p$, respectively.
\item[2.] To obtain
$\sigma_{t+1}$, for each connected component $\C$ of $G^+_p$ and of
$G^-_p$, with probability $1/2$ assign all vertices of $\C$ the same
sign $1$ and with probability  $1/2$ assign them all the sign $-1$,
independently for all these components.
\end{enumerate}

It is easy to show using Fortuin and Kasteleyn's Random Cluster
model \cite{FK} (see Edwards and Sokal \cite{ES} for this
derivation) that the Ising model measure is invariant under the SW
dynamics when $p=1-e^{-2\beta}$. Moreover, the SW dynamics is
clearly an aperiodic and irreducible Markov chain. Hence from any
starting configuration $\sigma_0$, the law of $\sigma_t$ obtained
after $t$ updates, converges in distribution to the stationary Ising measure.
Cooper, Dyer, Frieze and Rue \cite{CDFR} investigated the mixing
time of the SW dynamics on the complete graph on $n$ vertices. They
showed that if $p = {c \over n}$ when $c\in (0,\infty)\setminus \{2\}$ is some
constant independent of $n$, then the mixing time of the dynamics is
at most $O(\sqrt{n})$. The following Theorem
improves their result by giving the precise order of the mixing time
at all temperatures.

\begin{thm}\label{mainthm}
Consider the SW dynamics on the complete graph $K_{n}$ on $n$
vertices, with percolation parameter $p=\frac{c}{n}$, where $c$ is a
constant independent of $n$. Then,
\begin{enumerate}
\item[(i)] If $c>2$ then $\Tm = \Theta(\log n)$. \item[(ii)] If
$c=2$ then $\Tm = \Theta(n^{1/4})$. \item[(iii)] If $c<2$ then $\Tm
= \Theta(1)$.
\end{enumerate}
\end{thm}

The $q$-state ferromagnetic Potts model with parameter $\beta>0$ on
a graph $G=(V,\EE)$ is a probability measure on the set
$\Omega=\{1,2,\cdots,q\}^V$ where for each $\sigma\in \Omega$ $$
\P(\sigma) = { e^{\beta \sum_{(u,v)\in
\EE}\1_{\{\sigma(u)=\sigma(v)\}}} \over Z(G) } \, ,$$ and the {\em
partition} function $Z(G)$ is defined by
$$ Z(G) = \sum _{\sigma \in \Omega} e^{\beta \sum_{(u,v)\in \EE} \1_{\{\sigma(u)= \sigma(v)\}}} \,
.$$ We can similarly define the Swendsen-Wang dynamic with parameter $p\in[0,1]$ on $G$ as
the Markov chained defined by first performing $p$-bond percolation on the graphs spanned by the vertices of each state $\{1 ,\ldots ,q\}$, and then color all the components obtained this way uniformly from the $q$ colors and independently for each cluster. In a similar fashion the Potts model is the stationary measure for this chain.
\begin{thm}\label{thm_tree}
The mixing time of the Swendsen-Wang process for the $q$-state
ferromagnetic Potts model at any temperature on any tree with $n$ vertices is $O(\log
n)$, where the constants may depend only on the temperature.
\end{thm}

\subsection{Random graph estimates} Due to the percolative nature of the dynamics we require several estimates about the random graph $G(m,p)$ in non-critical case when $p={\theta\over m}$ where $\theta\neq1$ is a constant and the near-critical case when $mp = 1+o(1)$. We were not able to find such estimates in the vast random graph literature so we provide them in this paper. In the following we highlight some of these estimates which are interesting for the random graph community.

It is a well known result of Pittel \cite{P} that for $G(m,p)$ when $p={\theta\over m}$
where $\theta>1$ is a constant we have that ${|\C_1|-\beta m\over
\sqrt{m}}$ converges in probability to a normal distribution, $\beta$ is the unique positive solution of
the equation
\begin{equation*} 1-e^{-\theta x}=x.
\end{equation*} This does not imply, however, the moderate deviation bound on $|\C_1 - \beta m|$
\begin{equation}
\P(\big||\C_1|-\beta m\big|>A\sqrt{m})\leq
Ce^{-cA^2} \, ,
\end{equation} for any $A>0$, which we prove here, see Lemma \ref{realsupergnp}.
%
%

The study of the random graph ``inside'' the phase transition was
initiated by Bollob\' as \cite{B} where it is shown that if $p={1 +
\eps(m) \over m}$ when $\eps(m)$ is a positive sequence satisfying
$\eps(m) \geq m^{-1/3} \log m$, then with high probability
$|\C_1|=(2+o(1))\eps m$. The logarithmic corrections were removed by
Luczak \cite{L} and this statement holds whenever $\eps(m) \gg
m^{-1/3}$. A stronger result was recently proved by Pittel and
Wormald \cite{PW}. They show that in this regime of $p$, the
distribution of ${|\C_1| - 2\eps m \over \sqrt{m/\eps}}$ converges
to a normal distribution (this is a corollary of Theorem 6 of
\cite{PW}, but in fact the authors prove much more than this
statement).

Surprisingly however, the above results do not give good estimates
on $\E|\C_1|$ and on moderate and large deviations of $|\C_1|-2\eps
m$. These are crucial in our analysis of the Swendsen-Wang chain
since these determine the moments of the increments of the chain. In
Section \ref{secgnplemmas} we prove several such estimates, for
instance
$$ \E |\C_1| = 2\eps m + O(\eps^{-2} + \eps^2 m) \, ,$$
see the more accurate inequality in Theorem \ref{supergnp1}. Another
interesting estimate is a bound on the deviation of $|\C_1|$,
$$ \prob \Big ( \Big | |\C_1| - 2\eps m \Big | > A \sqrt{ {m \over
\eps} } \Big ) \leq Ce^{-{cA^2}} \, ,$$ for any $A$ satisfying $0 \leq
A \leq \sqrt{\eps^3 m}$. See Theorem \ref{gnpout2}.

\section{Mixing time preliminaries}

\begin{lem}
Suppose $(X_t, Y_t)$ is a coupling of two copies of the same Markov
chain with $X_0=x$, $Y_0=y$. We have
\begin{equation}\label{eqn_coupleprob}
\|\P^t(x,\cdot)-\P^t(y,\cdot)\|_{TV}\leq \P(X_t\not =Y_t).
\end{equation}
\end{lem}
One can apply triangle inequality to (\ref{eqn_coupleprob}) to get
$$ \|\P^t(x,\cdot)-\pi(\cdot)\|_{TV}\leq \max_{y}\P(X_t\not =Y_t),$$
and by taking maximum over all $x$, we have
\be\label{mixingtimethm}\max_{x}\|\P^t(x,\cdot)-\pi(\cdot)\|_{TV}\leq
\max_{x,y}\P(X_t\not =Y_t).\ee This gives the following theorem,
which we will use as a main technique to get the mixing time.
\begin{lemma}\label{thm_coupling}
If for every two state $x,y\in\Omega$, we could couple $\{X_t,
Y_t\}$ with $X_0=x$, $Y_0=y$ with a constant probability
$\epsilon>0$ after $L$ steps. Then, $\Tm(X,1-\epsilon)\leq L$.
\end{lemma}
As long as one can obtain the order of an upper bound of
$\Tm(1-\epsilon)$, the same order of upper bound holds for
$\Tm(1/4)$. We refer Section 4.4 and 4.5 of \cite{LPW} of detailed
discussion on this.

\section{Outline of the proof of Theorem \ref{mainthm}}
Due to the length of the proof of Theorem \ref{mainthm}, we provide
here a ``road-map'' of whole argument. The reader is advised to
follow this outline to get the general idea of the proof and go to
the main contents for further details whenever needed. Let
$\{\sigma_t\}_{t=0}^\infty$ be the SW Markov chain. Consider the
chain $X_t$ defined by
\begin{equation}\label{one_step}
 X_t = { \Big | \sum _{v} \sigma_t(v) \Big | } \, .
\end{equation}
Since the underlying graph is complete, $\{X_t\}$ is a Markov
chain with state space $\{0,\ldots,n\}$. Given $X_0$,
the random variable $X_1$ is determined by two independently drawn
random graphs $G({n+X_0\over 2},{c\over n})$ and $G({n-X_0\over
2},{c\over n})$. If we denote by $\{ \C^+ _j\}_{j \geq 1}$
and $\{ \C^- _j \}_{j \geq 1}$ the connected components of the
corresponding two random graphs, then $X_1$ is distributed as \be
\label{nextstep} \Big | \sum_{j\geq 1} \eps_j |\C^+_j| + \sum_{j
\geq 1} \eps '_j |\C^-_j| \Big | \, ,\ee where $\{\eps_j\}$ and
$\{\eps '_j\}$ are i.i.d. random variables taking $1$ with
probability $1/2$ and $-1$ otherwise. This is the reason that the
moments and large deviation estimates of random graph component sizes are
useful in our approach.

Frequently, to obtain upper bounds on the mixing time of the SW
chain we will obtain a bound on the mixing time of the chain $X_t$
and then use the following lemma, which appears in \cite{CDFR} and
is based on the path coupling idea of Bubley and Dyer \cite{BD}.
\begin{lemma}\label{pathcoupling}
Suppose $\{\sigma_t\}$ and $\{\sigma '_t\}$ are two SW chains such
that $X_0=X'_0$. There exists a coupling of the two chains such that
with probability at least ${1 \over 2}$ the two chains meet after
$O(\log n)$ steps.
\end{lemma}

\subsection{Outline of the Proof of Theorem \ref{mainthm} (i)}
By Lemmas \ref{thm_coupling} and \ref{pathcoupling} it suffices to
show that we can couple two copies of the magnetization chain $X_t$
and $X_t'$ such that they meet in $O(\log n)$ with probability $\Omega(1)$ which is uniform over all initial values $x_0$
and $x_0'$. It turns out that the stationary distribution is
concentrated in a window of length $\sqrt{n}$ around $\gamma_0 n$
for some $\gamma_0=\gamma_0(c)\in [0,1]$. In fact, the one step
evolution of $X_t$ essentially contracts the second moment of
$|X_t-\gamma_0n|$. That is, we have \be \E(X_1-\gamma_0 n)^2\leq
\delta (x_0 n-\gamma_0 n)^2+Bn \, ,\ee for some constants
$\delta\in(0,1)$ and large $B$. See Theorem \ref{lem_contrsup}. It follows
quickly that there exists an interval of values
$I=[\gamma_0n-A\sqrt{n},\gamma_0n+A\sqrt{n}]$ for some large constant
$A$ such that for any initial value $x_0$ we have that $X_t\in I$
with probability $\Omega(1)$ whenever $t=\Theta(\log n)$.

Once the two chains are both in the interval $I$ one can show that
they can be coupled to meet in the next step with
probability $\Omega(1)$. This is the content of Theorem \ref{lem_coupleinside}.
The main idea of that argument is that the random graph $G(n,
{c\over n})$ has $\Theta(n)$ isolated vertices with high probability
and that the difference of the sums of the spins of the two chains
{\em before} we assign spins to the isolated vertices is
$O(\sqrt{n})$. Thus, we one can couple the two chains to correct the
$O(\sqrt{n})$ error by assignment of those isolated vertices. This
follows from the classical local central limit theorem for the simple
random walk.


\subsection{Outline of the Proof of Theorem \ref{mainthm} (iii)} Since we need to prove an $O(1)$ upper bound we cannot use Lemma \ref{pathcoupling} here. However, the study of $X_t$'s evolution will still be useful. As in the supercritical case, the stationary measure is concentrated in a window of width $\Theta(\sqrt{n})$, but this time around $0$ and mixing occurs much faster. We will show, as before, that we have contraction, that is,
\begin{equation}
\E\big(X_1^2 \, |\, X_0\big) \leq \delta X_0^2 + Bn
\end{equation} for some constants $\delta\in(0,1)$ and large $B$ and for all $x_0 \in [0,1]$.  Moreover, if $0\leq X_0 \leq \frac{1}{c}-\frac{1}{2}$, we have
\begin{equation}
\E X_1^2 \leq B n \, ,
\end{equation} see Theorem \ref{lem_contrsub}. The first inequality
implies that $X_t$ will be in the window $[0,({1\over c}-{1\over
2})n]$ in $O(1)$ steps with probability $\Omega(1)$. The
second inequality implies that from this window $X_t$ jumps into
$[0,A\sqrt{n}]$ with high probability in just one more step. This
gives the mixing time upper bound on the chain $X_t$.

To go further and obtain a mixing time of the SW chain one needs to
consider the following two-dimensional chain.  For a starting
configuration $\sigma_0$, let $G_1$ denote the vertices with
positive spin and $G_2$ be its complement. Let $(Y_t, Z_t)$ be a
two-dimensional Markov chain, where $Y_t$ records the number of
vertices with positive spin in $G_1$ and $Z_t$ records the number
vertices with positive spin in $G_2$. By symmetry, the probability of the SW chain of being at $\sigma$ at time $t$ is the same for all $\sigma$ which have the same two-dimensional chain value. Consequently, the total variation distance of $\sigma_t$
from stationarity is the same as the total variation distance of the
two-dimensional chain from its stationary distribution. By Lemma
\ref{thm_coupling}, it suffices to provide a coupling of such
two-dimensional chains so that they meet in $O(1)$ steps
with probability $\Omega(1)$. By our previous argument,
$Y_t+Z_t$ will be in the window
$[{n\over2}-A\sqrt{n},{n\over2}+A\sqrt{n}]$ within $O(1)$
steps. One can show that once inside this window, such a coupling does
exist. See Proposition \ref{twodimensioncouple}. The idea is similar
to the proof of part (i) by considering the isolated vertices in
the two random graphs.

\subsection{Outline of the Proof of Theorem \ref{mainthm} (ii)} \label{roadmapcrit}
We again use Lemma \ref{pathcoupling} and Lemma
\ref{thm_coupling} to bound the mixing time of the magnetization
chain. However, to simplify our calculations, we will consider a
slight modification to the magnetization chain $X_t$. Instead of
choosing a random spin for each component after the percolation
step, we assign a positive spin to the {\em largest} component and random
spins for all other components. Let $X'_t$ be the sum of spins at
time $t$ (notice that we do {\em not} take absolute values here),
that is,
\begin{equation}
X_{t+1}' \stackrel{d}{=} \max\{|\C_1^+(t)
|,|\C_1^-(t)|\}+\eps\min\{|\C_1^+(t) |,|\C_1^-(t)|\}+\sum_{j\geq
2}\epsilon_j |\C_j^+(t) | + \sum_{j\geq 2}\epsilon_j' |\C_j^-(t)|\,
,
\end{equation}
where as usual $\eps, \{\eps_j\}$ and $\{\epsilon_j'\}$ are
independent mean zero $\pm$ signs. This chain has state space
$\{-n,\ldots,n\}$ and its absolute value has the same distribution
as our original chain $X_t$. As a consequence, any upper bound on
the mixing time of $X_t'$ implies the same upper bound on the mixing
time of $X_t$. For convenience, we will now denote this modified
chain by $X_t$. Let $X_t$ and $Y_t$ be two such chains such that
$X_t$ start from an arbitrary location $X_0$ and $Y_t$ starts from
the stationary distribution. We will show that we can couple $X_t$
and $Y_t$ so that they meet in $O(n^{1/4})$ steps with probability $\Omega(1)$. It will become evident that it suffices to restrict the
attention to $X_0\in[0,n]$. We will divide this into two subcases:
\begin{enumerate}
\item[(i)]$X_0\in[n^{3/4},n],$\item[(ii)]$X_0\in[0,n^{3/4}],$\end{enumerate}
and consider them separately, let us begin with case (i). In this
case the coupling strategy is as follows.

Consider the first {\em crossing time} of $X_t$ and $Y_t$, that is,
the first time $t$ such that
$\sign(X_t-Y_t)\neq\sign(X_{t-1}-Y_{t-1})$. We will show that this
is likely to occur only when the two chains take values
$\Theta(n^{3/4})$ and, more importantly, the distance between the
chains one step before the crossing time is of order $n^{5/8}$. This
is the content of Theorem \ref{thm_shoot}. The fact that one time
step before the crossing time is not a stopping time is problematic
and requires an {\em overshoot} estimate stating that the two chains
are not likely to cross each other from distance larger than $O(n^{5/8})$.
For random walks, these kind of estimates are classical (see for instance \cite{LL}).
The key estimate here is Theorem \ref{prop_nearhit}.

Next we show that when the chains take values $\Theta(n^{3/4})$ they
satisfying a local central limit theorem in scale $n^{5/8}$. In
particular, the chain has probability $\Omega(n^{-5/8})$ to move
to any point $x$ in an interval of size $\Theta(n^{5/8})$ around the
starting point.  We use the standard characteristic function technique to show this, see Lemma \ref{onelocalclt}. Now we are ready to conclude the proof in this case since we know that a step before the crossing time the chains have already been at distance $O(n^{5/8})$ from each other, so the local CLT provides a way to couple them in a few additional steps after the crossing time. See the proof of Theorem \ref{allplus}. \\

Let us consider now case (ii) in which $X_0\in[0,n^{3/4}]$. To
handle this case define $I=[-An^{2/3}, An^{2/3}]$ and proceed in two steps.
\begin{enumerate}\item With high probability $X_t$ will visit the interval $I$ by time $O(n^{1/4})$. This is proved in Theorem \ref{pushdownuse} and is based on the fact that the drift of the chain $|X_t|$ in this regime is approximately $-n^{1/2}$ (this is a small negative drift).
\item Once the chain is inside $I$, it will be pushed above $\Omega(n^{3/4})$ within $O(n^{1/4})$ steps. See Theorem \ref{pushup2}.
\end{enumerate}
Thus, with these two claims we see that in at most $O(n^{1/4})$
steps the chain is pushed into the $n^{3/4}$ regime and we may use
the theorems of case (i) to conclude. Let us briefly expand on the
proofs of (1) and (2).

The proof of (1) relies on the fact that the chain has a negative
drift of magnitude $\Omega(n^{1/2})$ as long as $X_t \not \in I$.
This follows rather easily from the random graph estimate Theorem
\ref{critgnp}. Note, however, that Theorem \ref{critgnp} estimates
the expected size of the
cluster discovered in time $\delta \eps m$ in the exploration
process for some small $\delta>0$ and {\em not} of the largest cluster $\C_1$.
We denote the former cluster by
$\C_{\delta \eps m}$ and remark that it has high probability of
being the largest. However, we were unable to prove the estimate of
Theorem \ref{critgnp} for $\C_1$ but only for $\C_{\delta \eps m}$.
This is the reason we need to consider yet another slight
modification of the magnetization chain: instead of giving a plus
sign to $\C_1$ and drawing random signs for the rest of the
clusters, we give the plus sign to $\C_{\delta \eps m}$ and the rest
receive random signs. From this point on the proof of (1) is rather
straightforward.

For the proof of (2) one has to show that when $X_t$ is in $I$, even
though the drift is negative there is still enough noise to
eventually push $X_t$ to the $n^{3/4}$ regime. We were unable to
pursue this strategy since it involves very delicate random graph
estimates we were unable to obtain. Instead we use the following
coupling idea. Since the stationary distribution normalized by $n^{3/4}$ has a weak limit with positive density at $0$, the expected number of visits to $I$ by the stationary chain before time $T$ is $\Theta\big ({n^{2/3} \over n^{3/4}} T \big )$. In Lemma \ref{pushup1} we show that when $T=\Theta(n^{1/4})$ the actual number of visits to $I$ is positive with high probability.
Next, to show that
$X_t$ is pushed upwards we start a stationary chain $Z_t$ and wait
until it enters $I$. We then couple $X_t$ and $Z_t$ such that they
meet inside $I$ and from that point they stay together. The only
technical issue with this strategy is how to perform the coupling of
$Z_t$ and $X_t$ inside $I$. This will follow, as before, from a
uniform lower bound stating that for any $x,x_0 \in I$ we have
$$\P(X_1=x \, \mid \, X_0=x_0)\geq cn^{-2/3} \, .$$
This estimate is done inside the scaling window of the random
graph phase transition and so the proofs are different from the previous
ones and require some combinatorial estimates. See Lemmas
\ref{criticalgnp} and Lemma \ref{uniformbd}.


%

\input{gnplemmas.tex}

\section{Supercritical case}
In this section we show that the mixing time of the Swendsen-Wang
chain is $\Theta(\log n)$ in the supercritical case $c>2$. This is
part (i) of Theorem \ref{mainthm}. Let $\{X_t=x_0n\}_{t\geq 0}$ be
the one dimensional chain defined in (\ref{one_step}). For
$x>\frac{2}{c}-1$ (so that $\frac{c(1+x)}{2}>1$), define
\begin{equation}\label{eqn_iteration}
\phi(x)=\beta\Big(\frac{c(1+x)}{2}\Big)\frac{1+x}{2},
\end{equation}
where $\beta(\cdot)$ is defined in (\ref{eqn_giantsize}). Since
$\beta:\R^+\to \R$ we have that $\phi:[-1,\infty]\to \R$. We begin with some preparations for the proof.


\begin{lem}\label{cor_superwalk}
For $c>2$, there exists an unique fixed point
$\gamma_0\in(1-\frac{2}{c},1)$ of $\phi(x)$. Furthermore, we have
\begin{equation} \frac{1}{2}<\phi'(x)<1 \qquad \hbox{ {\rm for} }
x>2c^{-1}-1 \, ,
\end{equation} and there exists a constant $\delta\in(0,1)$ such that for every
$x\in [0,1]\setminus\{\gamma_0\}$ we have
\begin{equation}\label{eqn_supwalkbound}
\half<\frac{\phi(x)-\gamma_0}{x-\gamma_0}\leq \delta.
\end{equation}
\end{lem}

\begin{thm}\label{lem_contrsup}
There exist constants $\delta\in(0,1)$ and $B>0$ such that
\be\label{eqn_contrsup} \E(X_{1}-\gamma_0 n)^2\leq \delta (X_0 - \gamma_0 n)^2 + Bn \, .\ee
\end{thm}

\begin{prop}\label{cor_concensup} We have
$$\E \Big ( X_1-\phi(x_0) n \, \mid X_0 \in [\gamma_0 n, n] \Big)^2 = O(n) \, .$$
\end{prop}

\begin{prop}\label{prop_supising}
If $X$ is distributed as the stationary distribution of the
magnetization Swendsen-Wang chain, then
$$\E(X-\gamma_0 n)^2= O(n).$$
\end{prop}

\begin{thm}\label{lem_coupleinside}
Suppose $X_0,Y_0$ are two magnetization Swendsen-Wang chains such
that $X_0,Y_0\in[\gamma_0 n-A\sqrt{n}, \gamma_0 n+A\sqrt{n}]$ where
$A$ is a constant, we can couple $X_1$ and $Y_1$ such that $X_1=Y_1$
with probability $\Omega(1)$ (which may depend on $A$).
\end{thm}

\noindent{\bf Proof of part (i) of Theorem \ref{mainthm}:} Rearranging Theorem \ref{lem_contrsup} and taking expectations gives $$\E(X_{t+1}-\gamma_0 n)^2-\frac{B}{1-\delta}n \leq\delta\Big[\E(X_t-\gamma_0 n)^2 -\frac{B}{1-\delta}n\big]$$ for all $t$.
We apply this inductively and get
$$\E(X_{C\log n}-\gamma_0 n)^2-\frac{B}{1-\delta}n\leq \delta^{C\log n}\Big[\E(X_0-\gamma_0 n)^2 -\frac{B}{1-\delta}n\Big] \, .$$ Hence, when $C=C(\delta)$ is large enough we get that
$$ \E(X_{C\log n}-\gamma_0 n)^2 = O(n) \, ,$$
and so Markov's inequality gives
\begin{equation}\label{superwindow1}
\P(|X_{C\log n}-\gamma_0n|\leq A \sqrt{n})\geq \frac{3}{4}
\end{equation} for some large constant $A$. Let $X_t'$ be a magnetization SW chain starting at stationarity. By Theorem \ref{prop_supising} and Markov's inequality we have
\begin{equation}\label{superwindow2}
\P(|X'_{C\log n}-\gamma_0n|\leq A\sqrt{n})\geq \frac{3}{4}
\end{equation} for some large constant $A$. Now, to couple $X_t$ and $X_t'$ we first run them independently until
time $C\log n$. By (\ref{superwindow1}) and (\ref{superwindow2}), we
have that $X_{C\log n}, X_{C\log n}'\in
[\gamma_0n-A\sqrt{n},\gamma_0n+A\sqrt{n}]$ with probability at least
$1/2$. By Theorem \ref {lem_coupleinside}, we can couple $X_{C\log
n+1}$ and $X_{C\log n+1}'$ such that $X_{C\log n+1}=X_{C\log n+1}'$
with probability $\Omega(1)$. Then by Lemma
\ref{pathcoupling}, we have that $\{\sigma_t\}$ and $\{\sigma_t'\}$
can be coupled such that $\sigma_t=\sigma_t'$ in $O(\log n)$ steps
with probability $\Omega(1)$. The upper bound of mixing time
follows from Lemma \ref{thm_coupling}.

For the lower bound, we will show that if $X_0=n$, then
$$\|X_{\alpha\log n}-X_{\pi}\|_{TV}\geq 1/4$$ for some small constant $\alpha>0$, where $X_{\pi}$ is the stationary distribution of magnetization Swendsen-Wang chain. By (\ref{eqn_supwalkbound}), we have that
\begin{eqnarray*}
\P\Big(X_{t+1}-\gamma_0 n\leq \frac{1}{4}(X_t-\gamma_0 n)\Big)
&\leq& \P\Big(X_{t+1}-\gamma_0 n\leq \frac{1}{2}(\phi(X_t/n) n - \gamma_0 n )\Big)\\
&=&\P\Big(X_{t+1}-\phi(X_t/n) n \leq
\frac{1}{2}(\gamma_0 n -\phi(X_t/n)n) \Big ) \, .
\end{eqnarray*}
When $X_t \geq \gamma_0 n$ we have that $\phi(X_t/n) n \geq \gamma_0 n$ by (\ref{eqn_supwalkbound}), hence
Proposition \ref{cor_concensup} and Markov's
inequality imply that
\begin{equation}\label{markov1}\P\Big(X_{t+1}-\gamma_0 n\leq
\frac{1}{4}(X_t-\gamma_0 n) \mid X_t \geq \gamma_0 n \Big) \leq
\frac{O(n)}{(\phi(X_t/n)n-\gamma_0 n)^2}.\end{equation} Furthermore, if
$X_t-\gamma_0 n \geq n^{\frac{3}{4}}$, then
$\phi(X_t/n)n-\gamma_0 n \geq n^{\frac{3}{4}}/2$ by
(\ref{eqn_supwalkbound}). Plugging this into (\ref{markov1}) gives
\begin{equation}\label{inductive}\P\Big(X_{t+1}-\gamma_0 n \geq
\frac{1}{4}(X_t-\gamma_0 n)\big|X_t-\gamma_0 n\geq
n^{\frac{3}{4}}\big) \geq 1- O(n^{-\frac{1}{2}}).\end{equation}
Starting from $X_0=n$, by applying (\ref{inductive}) iteratively we
have \begin{equation}\P\big(X_{\alpha\log n}-\gamma_0 n\geq
n^{\frac{3}{4}}\big) \geq (1- O(n^{-\frac{1}{2}}))^{\alpha\log
n} = 1-o(1) \, ,\end{equation} when $\alpha>0$ is small enough constant. On the
other hand, by Proposition \ref{prop_supising} and the Markov's
inequality, we have $\P\Big(|X_{\pi}-\gamma_0 n|\geq
A\sqrt{n}\Big)\leq \frac{1}{4}$ for some constant $A$. Putting the
two inequalities together, we get
\begin{equation}\label{eqn_tvdistsup}
\|X_{\alpha\log _4 n}-X_{\pi}\|_{TV}\geq \frac{3}{4}-o(1)\geq
\frac{1}{4},
\end{equation} which gives a lower bound on the mixing time of magnetization SW chain $X_t$. This concludes the proof since any lower bound of the mixing time of $X_t$ implies the same lower bound of mixing time of $\sigma_t$.\qed\\

\noindent{\bf Proof of Lemma \ref{cor_superwalk}:} By the definition
of $\beta(\cdot)$ in equation (\ref{eqn_giantsize}), we know
$\phi(x)$ is the positive solution of
\begin{equation}\label{eqn_phi}
1-e^{-c \phi(x)}=\frac{2\phi(x)}{x+1}
\end{equation} for all $x>\frac{2}{c}-1$.
Taking derivative of both sides yields
$$c e^{-c \phi}\phi'=
\frac{2(x+1)\phi'-2\phi}{(x+1)^2}.$$ By plugging in
$x+1=\frac{2\phi}{1-e^{-c \phi}}$ we get
\begin{equation}\label{eqn_dydx}
\phi'=\frac{1-2e^{-c \phi}+e^{-2c \phi}}{2(1-e^{-c \phi}-c\phi e^{-c
\phi})}.
\end{equation}
By (\ref{eqn_dydx}), we have that $\half<\phi'$ if and only if $e^{-c\phi}> 1-c \phi$ which is true for all $c \phi>0$. We also have
that \be\label{derivativeref}\phi'<1 \Longleftrightarrow c \phi<\sinh (c
\phi)\ee which holds for all $c \phi>0$.

Since $c>2$ (which implies $\frac{2}{c}-1<0$), we have that $\phi
'<1$ for all $x\in[0,1]$. Since $\phi'$ is continuous, we have a
constant $\delta_1\in(0,1)$ such that $\phi'(x)<\delta_1$ for all
$x\in[0,1]$. Note that $\phi (0)=\frac{1}{2} \beta(\frac{c}{2})>0$,
$\phi (1)=\beta (c)<1$ and $\phi$ is strictly increasing in
$[0,1]$, we have by Rolle's theorem that there exists an unique
point $\gamma_0\in(0,1)$ such that $\phi(\gamma_0)=\gamma_0$. By
plugging in $x=1-2/c$ into (\ref{eqn_iteration}) and the definition
of $\beta(\cdot)$, we get
\begin{eqnarray*}
\phi(1-2/c)>1-2/c&\Leftrightarrow&\beta(c-1)>1-\frac{1}{c-1}\\&\Leftrightarrow&e^{-(c-2)}<\frac{1}{c-1}
\end{eqnarray*}
which is always true for $c>2$. It follows immediately that $\gamma_0>1-\frac{2}{c}$. \qed \\

Recall that given $X_0=x_0n$, we have that $X_1$ is distributed
as in (\ref{nextstep}). To prove Theorem \ref{lem_contrsup} we first
state a useful lemma.

\begin{lem}\label{lem_supsmall} Let $c>2$.
\begin{itemize}
\item[(i)] There exists a non-negative function $h(\cdot)$ with $h(\eps)\rightarrow0$ as $\eps\rightarrow 0$ such that if
$|x_0-(1-\frac{2}{c})|\leq \eps$, then $\E\Big(\sum_{j\geq 1}|\C^-_j|^2\Big) \leq h(\eps)n^2$.
\item[(ii)] For any fixed $\eps>0$, if $x_0\in [0, (1-\frac{2}{c}-\eps)]$, then $\E\Big(\sum_{j\geq 1}|\C^-_j|^2\Big)\leq(\phi^2(-x_0)+o(1))n^2$.
\item[(iii)] For any fixed $\eps>0$, if $x_0\in [1-\frac{2}{c}+\eps,1]$, then $\E\Big(\sum_{j\geq 1}|\C^-_j|^2\Big)\leq O(n)$.
\end{itemize}
\end{lem}

\noindent{\bf Proof of Theorem \ref{lem_contrsup}:} By
(\ref{nextstep}), we have
\begin{equation}\label{eqn_secmoment}
\E X_1^2=\E\Big(\sum_{j\geq 1} |\C^+_j|^2\Big) +\E\Big(\sum_{j\geq
1}|\C^-_j|^2\Big)
\end{equation}
and
\begin{eqnarray}\label{eqn_absolute}
\E X_1&=&\E\Big | \sum_{j\geq 1} \eps_j |\C^+_j| + \sum_{j \geq 1}
\eps '_j |\C^-_j| \Big |=\E\Big | \eps_1\big(\sum_{j\geq 1} \eps_j
|\C^+_j| + \sum_{j \geq 1} \eps '_j |\C^-_j| \big)\Big
|\nonumber\\&\geq&\E\Big( \eps_1\big(\sum_{j\geq 1} \eps_j
|\C^+_j| + \sum_{j \geq 1} \eps '_j |\C^-_j| \big)\Big)\nonumber\\
&=& \E\Big[|\C_1^+|+\eps_1\Big(\sum_{j\geq 2} \eps_j |\C^+_j| +
\sum_{j \geq 1} \eps '_j |\C^-_j|\Big)\Big]=\E |\C_1^+|.
\end{eqnarray}

Combining (\ref{eqn_secmoment}) and (\ref{eqn_absolute}), we get
\be\label{eqn_supersplit} \E\Big(X_1-\gamma_0 n\Big)^2\leq
\E\Big(\sum_{j\geq 1} |\C^+_j|^2\Big) +\E\Big(\sum_{j\geq
1}|\C^-_j|^2\Big)-2\gamma_0n\cdot \E |\C_1^+|+\gamma_0^2n^2.\ee

The random graph $G(\frac{1+x_0}{2}n,\frac{c}{n})$ is supercritical
with $\theta=\frac{1+x_0}{2}n\frac{c}{n}\geq\frac{c}{2}>1$. By Corollary
\ref{lem_superdeviation} we have
\begin{equation}\label{eqn_supersquare}
\E\Big(\sum_{j\geq 1}|\C^+_j|^2\Big)\leq
\Big(\E|\C^+_1|\Big)^2+O(n).
\end{equation}
Plugging (\ref{eqn_supersquare}) into (\ref{eqn_supersplit}), we get
\be\label{eqn_supsplit2} \E\Big(X_1-\gamma_0 n\Big)^2\leq
\Big(\E|\C^+_1|-\gamma_0 n\Big)^2 +\E\Big(\sum_{j\geq
1}|\C^-_j|^2\Big)+O(n).\ee
By Corollary \ref{lem_giant}, we have $\Big|\E|\C^+_1|-\phi (x_0)
n\Big|\leq O(\sqrt{n})$. Thus,
\begin{eqnarray}\label{eqn_supsplit2.5}
\Big(\E|\C^+_1|-\gamma_0 n\Big)^2 &\leq& \Big|\E|\C^+_1|-\phi(x_0)
n\Big|^2+ \big|\phi(x_0)n-\gamma_0
n\Big|^2\nonumber\\&+&2\Big|\E|\C^+_1|-\phi(x_0) n\Big|\Big|\phi(x_0)n-\gamma_0 n\Big|\nonumber\\
&\leq& \Big|\phi(x_0)n-\gamma_0
n\Big|^2+O(\sqrt{n})\Big|\phi(x_0)n-\gamma_0 n\Big|+O(n).
\end{eqnarray}
Applying Lemma \ref{cor_superwalk} gives that
\begin{eqnarray}\label{eqn_supsplit3}
\Big(\E|\C^+_1|-\gamma_0 n\Big)^2 &\leq& \delta^2_1|x_0-\gamma_0|^2
n^2+ |x_0-\gamma_0|O(n^{3/2})+O(n).
\end{eqnarray}
If $|x_0-\gamma_0| = O(n^{-\frac{1}{2}})$, then
$|x_0-\gamma_0|n^{3/2}=O(n)$. If
$|x_0-\gamma_0|n^{\frac{1}{2}}\rightarrow\infty$, we have
$|x_0-\gamma_0|O(n^{3/2})=o(|x_0-\gamma_0|^2 n^2)$. Plugging these
back into (\ref{eqn_supsplit3}), we get \be\label{eqn_supsplit4}
\E(X_1-\gamma_0
n)^2\leq(\delta_1^2+o(1))|x_0-\gamma_0|^2n^2+O(n)+\E\Big(\sum_{j\geq
1}|\C^-_j|^2\Big). \ee

To estimate $\E\Big(\sum_{j\geq 1}|\C^-_j|^2\Big)$, choose a small
constant $\eps$ such that $\delta_1^2+h(\eps)<1$ where $h(\cdot)$
is defined in part (i) of Lemma \ref{lem_supsmall}. If $\Big| x_0-(1-\frac{2}{c})\Big|< \eps$, we have that
\be\label{overallcase1}\E(X_1-\gamma_0 n)^2\leq
(\delta_1^2+h(\eps))|x_0-\gamma_0|^2 n^2+O(n)\ee by plugging part (i) of Lemma
\ref{lem_supsmall} into (\ref{eqn_supsplit4}).

If $x_0\in [0, (1-\frac{2}{c}-\eps)]$, we have that
$|\phi(x_0)-\gamma_0|$ is uniformly bounded from below by Lemma
\ref{cor_superwalk}. As a result, we have that
$\Big(\E|\C^+_1|-\gamma_0 n\Big)^2\leq
(\phi(x_0)-\gamma_0)^2n^2+O(n)$ in (\ref{eqn_supsplit2.5}). Plugging
this and part (ii) of Lemma \ref{lem_supsmall} into
(\ref{eqn_supsplit2.5}) gives
\be\label{overallcase2.5}\E(X_1-\gamma_0 n)^2\leq
((\phi(x_0)-\gamma_0)^2+\phi^2(-x_0)+o(1))n^2+O(n).\ee By Lemma
\ref{cor_superwalk} and Rolle's Theorem, we have
$$\frac{\phi(x_0)-\phi(-x_0)}{2x_0}\geq \frac{1}{2},$$
which leads to $\phi^2(-x_0)\leq \Big(\phi(x_0)-x_0\Big)^2$. This
gives \be\label{eqn_closed} \frac{\Big(\phi(x_0)-\gamma_0\Big)^2
+\phi^2(-x_0)}{(x_0-\gamma_0)^2}\leq
\frac{\Big(\phi(x_0)-\gamma_0\Big)^2
+\Big(\phi(x_0)-x_0\Big)^2}{(x_0-\gamma_0)^2}<1, \ee since
$x_0<\phi(x_0)<\gamma_0$. The left hand side of (\ref{eqn_closed})
is smaller than $1$ for all $x_0 \in[0, (1-\frac{2}{c}-\beta)]$, so
it is smaller than some constant $\delta_2<1$ uniformly. Plugging
this into (\ref{overallcase2.5}), we get \be\label{overallcase2}
\E(X_1-\gamma_0 n)^2\leq \delta_2(x_0-\gamma_0)^2 n^2+O(n).\ee If
$x_0\in [1-\frac{2}{c}+\eps,1]$, we plug (iii) of Lemma
\ref{lem_supsmall} into (\ref{eqn_supsplit4}) and obtain
\be\label{overallcase3} \E(X_1-\gamma_0
n)^2\leq(\delta_1^2+o(1))|x_0-\gamma_0|^2n^2+O(n).\ee
Combining (\ref{overallcase1}),(\ref{overallcase2}) and (\ref{overallcase3}) concludes our proof.\qed\\

\noindent{\bf{Proof of Lemma \ref{lem_supsmall}.}} We begin with case (ii). In this regime, the random
graph $G(\frac{1-x_0}{2}{n},\frac{c}{n})$ is supercritical with
$\theta>1+\frac{c\eps}{2}$. In the same way we obtained (\ref{eqn_supersquare}) we also have
\begin{equation}\label{eqn_supersmall}
\E\Big(\sum_{j\geq 1}|\C^-_j|^2\Big)\leq \Big(\E|\C^-_1|\Big)^2
+O(n).
\end{equation}
By Corollary \ref{lem_giant} we have that $|\E|\C^-_1|-\phi (-x_0)
n|\leq O(\sqrt{n})$ showing that
\begin{eqnarray}\label{eqn_case2}
\nonumber \E\Big(\sum_{j\geq 1}|\C^-_j|^2\Big)&\leq&
\Big(\phi(-x_0)n+O(\sqrt{n})\Big)^2
+O(n)\nonumber\\
&\leq&\phi^2(-x_0)n^2+O(n)+\phi(-x_0)O(n^{3/2})\nonumber\\&=&\big(\phi^2(-x_0)+o(1)\big)n^2,
\end{eqnarray}
since $|\phi(x_0)|$ is uniformly bounded from below, as required.

We now prove case (i). Note that we have
$$\E\Big(\sum_{j\geq
1}|\C^-_j|^2\Big)=\Big(\frac{1-x_0}{2}n\Big)\E|C_v|$$ as in
(\ref{eqn_useful}). Since $\E|C_v|$ is decreasing in $x_0$, we have
that $\E\Big(\sum_{j\geq 1}|\C^-_j|^2\Big)$ reaches its maximum at
$x_0=1-\frac{2}{c}-\eps$. Plugging in this value into
(\ref{eqn_case2}) gives \be\label{eqn_supsup2} \E\Big(\sum_{j\geq
1}|\C^-_j|^2\Big) \leq \big(\beta^2(1+\frac{\eps
c}{2})+o(1)\big)n^2 \ee Note that $\beta(x)\rightarrow 0$ as
$x\rightarrow1$, so we can take $h(x)=\beta^2(1+\frac{cx}{2})+o(1)$.

To prove case (iii) note that
$\frac{c(1-x_0)}{2} \leq 1-\frac{\eps c}{2}$, so the
random graph $G(\frac{1-x_0}{2}{n},\frac{c}{n})$ is subcritical in
this regime with $\theta$ bounded from above away from $1$. Applying
Lemma \ref{lem_subdeviation}, we get
\begin{equation}\label{eqn_high}
\E\Big[\sum_{j\geq 1}|\C_j^- |^2\Big] = O(n).
\end{equation} \qed

\noindent{\bf{Proof of Proposition \ref{cor_concensup}:}} Note that
(\ref{eqn_supsplit2}) is valid for all $\gamma_0\in[0,1]$ and in
particular for $\phi(x_0)$. Thus, \be\label{eqn_supsplit6}
\E\Big(X_1-\phi(x_0) n\Big)^2\leq \Big(\E|\C^+_1|-\phi(x_0) n\Big)^2
+\E\Big(\sum_{j\geq 1}|\C^-_j|^2\Big)+O(n).\ee

Recall that $x_0\geq\gamma_0>1-\frac{2}{c}$. By Corollary
\ref{lem_giant} we have that $\big(\E|\C^+_1|-\phi(x_0) n\big)^2 =
O(n)$. The random graph of $G(\frac{1-x_0}{2}n,\frac{c}{n})$ is in
regime (iii) of Lemma \ref{lem_supsmall}. Plugging (\ref{eqn_high})
into (\ref{eqn_supsplit6}), we get
\begin{equation}\label{eqn_supconcen}
\E\Big(X_1-\phi(x_0) n\Big)^2 = O(n) \, ,
\end{equation}
as required. \qed

\noindent{\bf{Proof of Proposition \ref{prop_supising}:}} If $X_0$
follows the stationary distribution of the magnetization SW chain,
so does $X_1$. Taking expectation of both sides of
(\ref{eqn_contrsup}) gives
\begin{eqnarray*}
\E(X_1-\gamma_0n)^2\leq\delta\E(X_0-\gamma_0n)^2+Bn \, ,
\end{eqnarray*}
as required. \qed \\

To prove Theorem \ref{lem_coupleinside} we need the following lemma.
\begin{lem}\label{lem_CLT} Let $Y$ and $Z$ be two random variables distributed as the sum of $n$ independent random $\pm$ signs. Then for any fixed constant $a$, there exists a constant $\kappa(a)\in(0,1]$ such that
for any $-a\sqrt{n}\leq y\leq a\sqrt{n}$, we can couple $Y$ and
$Z$ such that $Y-y=Z$ with probability at least $\kappa$.
\end{lem}
\noindent{\bf{Proof.}} Direct corollary of the local central
limit theorem of simple random walk. \qed \\

\noindent{\bf{Proof of Theorem \ref{lem_coupleinside}.}}  To couple
$X_1$ and $Y_1$, we first apply the percolation step of the Swendsen-Wang
dynamics in both chains independently. By Lemma \ref{lem_isolated},
with probability $1-O(\frac{1}{n})$, the number of isolated points
after percolation is bigger than $\frac{1}{3 e^{c}}n$ in both
chains. Conditioned on this, we assign each component a $\pm$ spin using the following procedure.

First assign the spins of components independently in descending
order of their size until there are $\frac{1}{3 e^{c}}n$ components
left. Note the remaining components are all isolated vertices. Denote
by $\bar{X}_1$ and $\bar{Y}_1$ as the absolute value of the sum of
spins at this time respectively.

Note that (\ref{eqn_absolute}) and (\ref{eqn_supersplit}) are still
valid if we replace $X_1$ by $\bar{X}_1$. Consequently, Theorem
\ref{lem_contrsup} is also valid if replacing $X_1$ by $\bar{X}_1$.
Hence, since $|X_0-\gamma_0n| \leq A\sqrt{n}$ we have
$$\E(\bar{X}_1-\gamma_0n)^2 = O(n) \, .$$ By Markov's inequality, there
exists a constant $A_1$ such that
\begin{equation}\label{eqn_xbar}
\P\Big(|\bar{X}_1-\gamma_0 n|\geq A_1\sqrt{n}\Big)\leq \frac{1}{4} \, ,
\end{equation}
and similarly
\begin{equation}\label{eqn_ybar}
\P\Big(|\bar{Y}_1-\gamma_0 n|\geq A_1\sqrt{n}\Big)\leq \frac{1}{4} \, .
\end{equation}

Consider the event
$$\mathcal{A}:=\{|\bar{X}_1-\gamma_0 n|<A_1\sqrt{n}\}\cap\{ |\bar{Y}_1-\gamma_0 n|<
A_1\sqrt{n}\}\cap\{\hbox{{\rm There are at least ${n \over 3e^c}$ isolated vertices}}\} \, .$$
By (\ref{eqn_xbar}) and (\ref{eqn_ybar}) we
have that $\P(\mathcal{A})\geq \frac{1}{4}$.

Conditioned on $\mathcal{A}$, we have $|\bar{X}_1-\bar{Y}_1|\leq
2A_1\sqrt{n}$. Denote by $\hat{X}_1$ and $\hat{Y}_1$ the sum of
spins of the rest of the components (all of them being isolated vertices) of the two
chains respectively. Note $\hat{X}_1$ and $\hat{Y}_1$ are i.i.d. sums
of $\pm$ spins. By Lemma
\ref{lem_CLT} we can couple $\hat{X}_1$ and $\hat{Y}_1$ so that
$\hat{X}_1+\bar{X}_1=\hat{Y}_1+\bar{Y}_1$ with
probability $\Omega(1)$. Finally, notice that $X_1\stackrel{(d)}{=}|
\bar{X}_1+\hat{X}_1 |$  and $Y_1\stackrel{(d)}{=}|
\bar{Y}_1+\hat{Y}_1 |$, concluding the proof.\qed


\section{Subcritical case}
In this section, we prove that in the subcritical case $c<2$, the
mixing time of Swendsen-Wang chain is $\Theta(1)$. This is part
(iii) of Theorem \ref{mainthm}.
\begin{lem}\label{cor_subwalk}
For $c\in(1,2)$ there exists a constant $\delta\in(0,1)$ such that
for all $x\in[\frac{2}{c}-1,1]$, we have
\begin{equation}\label{eqn_subwalkbound}
\frac{\phi(x)}{x}\leq \delta
\end{equation} where $\phi(\cdot)$ is defined in (\ref{eqn_iteration}).
\end{lem}

\begin{theorem}\label{lem_contrsub}
There exist two constants $\delta\in(0,1)$ and $B>0$ such that
\begin{equation}\label{eqn_contrsub}
\E\big(X_1^2\mid X_0 \big)\leq \delta X_0^2 + Bn \, .
\end{equation}
Moreover, if $0\leq x_0\leq \frac{1}{c}-\frac{1}{2}$, we have
\begin{equation}\label{subwindow1}
\E X_1^2\leq B n.
\end{equation}
\end{theorem}

To get the constant upper bound of mixing time we need to consider
the following two-dimensional chain. Let $G_1$ be a fixed subset of
the vertices and $G_2$ its complement. Let
$(Y_t, Z_t)$ be a two-dimensional Markov chain, where $Y_t$ record
the number of vertices with positive spin in $G_1$ and $Z_t$ record
the number vertices with positive spin in $G_2$.

\begin{prop}\label{twodimensioncouple}
Let $(Y_t, Z_t)$ and $(\widetilde{Y}_t, \widetilde{Z}_t)$ be two
two-dimensional chains as defined above. Suppose $Y_0+Z_0$ and
$\widetilde{Y}_0+\widetilde{Z}_0$ lie in the window
$I=[\frac{n}{2}-A\sqrt{n}, \frac{n}{2}+A\sqrt{n}]$ where $A$ is a
constant. Then we can couple $(Y_1,Z_1)$ and $(\widetilde{Y}_1,
\widetilde{Z}_1)$ such that $(Y_1,Z_1)=(\widetilde{Y}_1,
\widetilde{Z}_1)$ with probability $\Omega(1)$ (which may depend on $A$).
\end{prop}

\noindent{\bf{Proof of part (iii) of Theorem \ref{mainthm}:}} For
any starting configuration $\sigma$, let $G_1$ be the vertices with
positive spin and $G_2$ be its complement. Let $X_t$ be the
magnetization chain and $(Y_t, Z_t)$ be the two-dimensional chain as
described above. Let $\pi$ be the stationary distribution of
Swendsen-Wang chain and $\widetilde{\pi}$ be the stationary
distribution of $(Y_t, Z_t)$. By symmetry, configurations with same
two-dimensional chain value have same distributions for any $t$.
Consequently
\be\label{twodimensiontv}
||\sigma\P^t-\pi||_{TV}=||(|G_1|,0)\P^t,\widetilde{\pi}||_{TV} \, . \ee

Thus, by Lemma \ref{thm_coupling} it suffices to couple the chains $(Y_t, Z_t)$ and $(\widetilde{Y}_t,
\widetilde{Z}_t)$ such that they meet with
probability $\Omega(1)$ in time $t=\Theta(1)$. By Lemma \ref{lem_contrsub}, we have
\begin{equation*}
\E (X_{t+1}^2) -\frac{B}{1-\delta}n\leq \delta\Big[\E
(X_t^2)-\frac{B}{1-\delta}n\Big].
\end{equation*}
Applying this inductively we get
\begin{equation*}
\E (X_{t}^2) -\frac{B}{1-\delta}n\leq \delta^t\E (X_0^2)\leq
\delta^t n^2.
\end{equation*}
For $t\geq 2\log_{\delta}\frac{1}{8}(\frac{1}{c}-\frac{1}{2})$ and
large $n$, we have $$\E (X_{t}^2)\leq
\frac{1}{4}\big(\frac{1}{c}-\frac{1}{2}\big)^2 n^2.$$
For such $t$ Markov's
inequality gives \be\label{subwindow}\P \Big(X_t \geq
\big(\frac{1}{c}-\frac{1}{2}\big)n\Big)\leq \frac{1}{4}.\ee

By Theorem \ref{lem_contrsub} and Markov's
inequality, if $X_t\in[0,(\frac{1}{c}-\frac{1}{2})n]$, then $X_{t+1} \in [0,A\sqrt{n}]$ with
probability at least $1/2$ for some large constant $A$. Combining this and
(\ref{subwindow}), we have that after constant number of steps, the
chain $X_t$ will jump into the window $I=[0,A\sqrt{n}]$ with
probability $\Omega(1)$.

For any two Swendsen-Wang chains $\sigma$ and $\tilde{\sigma}$, Let
$X_t$ and $\widetilde{X}_t$ be the corresponding magnetization
chains. Running the two Swenden-Wang dynamics independently first,
by the argument above, we have that $X_t$ and $\widetilde{X_t}$ both
jump into $[0,A\sqrt{n}]$ after constant steps with probability $\Omega(1)$. By Proposition \ref{twodimensioncouple}, we
can couple the two two-dimensional chains so that $(Y_t,
Z_t)=(\widetilde{Y}_t, \widetilde{Z}_t)$ with probability $\Omega(1)$, which concludes the whole proof. \qed \\

\noindent{\bf{Proof of Lemma \ref{cor_subwalk}:}} Note $\phi$ is
differentiable on $[\frac{2}{c}-1,1]$. Recalling
(\ref{derivativeref}), we have $\phi'<1$ for all $x>\frac{2}{c}-1$.
By Rolle's Theorem, we have $\phi (x)-0\leq x-(\frac{2}{c}-1)$ for
all $x>\frac{2}{c}-1$. So
$$\frac{\phi (x)}{x}\leq 1-\frac{\frac{2}{c}-1}{x}\leq
1-(\frac{2}{c}-1)$$ for all $x\in[\frac{c}{2}-1,1]$. \qed \\

\noindent{\bf{Proof of Theorem \ref{lem_contrsub}:}} We use the fact that
(\ref{eqn_secmoment}) is still valid. The random graph
$G(\frac{1-x_0}{2}n,\frac{c}{n})$ is subcritical with
$\theta=(\frac{1-x_0}{2}n)\frac{c}{n} = \frac{c}{2}$. By Lemma
\ref{lem_subdeviation}, we have
\begin{equation}\label{eqn_qq}
\E\Big(\sum_{j\geq 1}|\C^-_j|^2\Big)= O(n).
\end{equation} If $c<1$, the random graph $G(\frac{1+x_0}{2}n,\frac{c}{n})$ is
subcritical with $\theta=(\frac{1+x_0}{2}n)\frac{c}{n}\leq c<1$. By
Lemma \ref{lem_subdeviation}, we have
\begin{equation}\label{eqn_qq1}
\E\Big(\sum_{j\geq 1}|\C^+_j|^2\Big) = O(n).
\end{equation} If $c \geq 1$, then let $\eps>0$ be a
small constant that we will determine later and consider the
following three cases.

(i) $0\leq x_0 \leq \frac{2}{c}-1-\eps$. In this case, the random
graph $G(\frac{1+x_0}{2}n,\frac{c}{n})$ is subcritical with
$\theta\leq 1-\frac{\eps c}{2}$. By Lemma \ref{lem_subdeviation},

\begin{equation}\label{subcase1}
\E\Big(\sum_{j\geq 1}|\C^+_j|^2\Big)= O(n).
\end{equation}

(ii) $\frac{2}{c}-1+\eps\leq x_0\leq 1$(in case $c>1$). In this
case, the random graph $G(\frac{1+x_0}{2}n,\frac{c}{n})$ is
supercritical with $\theta\geq1+\frac{\eps c}{2}$. By Corollary
\ref{lem_superdeviation}, we have
\begin{eqnarray*}
\E\Big(\sum_{j\geq 1}|\C^+_j|^2\Big)&\leq& \big(\E |\C^+_1|\big)^2+
O(n)\\
&=& \Big(\phi(x_0)n\Big)^2+\Big(\E|\C^+_1|-\phi(x_0)n\Big)
\Big(\E|\C^+_1|+\phi(x_0)n\Big)+O(n).
\end{eqnarray*}
By Corollary \ref{lem_giant}, we have $\Big|
\E|\C^+_1|-\phi(x_0)n\Big|= O(\sqrt{n})$. By Lemma
\ref{cor_subwalk}, we have $\phi(x_0)n\leq \delta x_0 n$. So we have
\begin{eqnarray*}
\E\Big(\sum_{j\geq 1}|\C^+_j|^2\Big)&\leq& \delta^2 x_0^2
n^2+O(n^{3/2}) \leq (\delta^2+o(1))x_0^2n^2.
\end{eqnarray*}

(iii) $\frac{2}{c}-1-\eps\leq x_0\leq \frac{2}{c}-1+\eps$ (or
$1-\eps\leq x_0\leq 1$ in case $c=1$). Recall that
$\E\Big(\sum_{j\geq 1}|\C^+_j|^2\Big)=\frac{1+x_0}{2}n\E|C_v|$. So
$\E\Big(\sum_{j\geq 1}|\C^+_j|^2\Big)$ reaches its maximum at
$x_0=\frac{2}{c}-1+\eps$ for $1<c<2$ or $x_0=1$ for $c=1$. In the
former case, by the estimate in case (ii), we get
\begin{equation*}
\E\Big(\sum_{j\geq 1}|\C^+_j|^2\Big)\leq
(\delta^2+o(1))(\frac{2}{c}-1+\eps)^2 n^2.
\end{equation*}
Now we choose $\eps$ to be small enough such that
$\delta^2\Big(\frac{2/c-1+\eps}{2/c-1-\eps}\Big)^2< 1$, then we
choose a constant $\delta_1$ such that $\delta^2
\Big(\frac{2/c-1+\eps}{2/c-1-\eps}\Big)^2<\delta_1<1$. Then We have
\begin{equation*}
\E\Big(\sum_{j\geq 1}|\C^+_j|^2\Big)\leq \delta x_0^2 n^2.
\end{equation*}

In the latter case, by Theorem 1 of \cite{NP1}, we have that
\begin{equation*} \E\Big(\sum_{j\geq 1}|\C^+_j|^2\Big)=o(n^2).
\end{equation*} The
Lemma follows from combining case (i), (ii) and (iii).\qed \\

\noindent{\bf{Proof of Proposition \ref{twodimensioncouple}:}}
Suppose without lost of generality that $|G_2|\leq |G_1|$. Since
$Y_0+Z_0\in I$, the random graphs $G(Y_0+Z_0,\frac{c}{n})$ and
$G(n-(Y_0+Z_0),\frac{c}{n})$ are both subcritical for large $n$. The
same is true for the chain
$(\widetilde{Y_t},\widetilde{Z_t})$. In the first chain after the
percolation step, denote by $\{\A_j\}_{j\geq1}$ and
$\{\B_j\}_{j\geq1}$ the components with vertices completely in $G_1$
and $G_2$ respectively. Note that there are also components that
have vertices in both $G_1$ and $G_2$. Denote such components by
$\{\C_j\}_{j\geq1}$. In the second chain, we denote by
$\widetilde{\A_j}$,$\widetilde{\B_j}$ and $\widetilde{\C_j}$ to be
these components. Lemma \ref{lem_isolated} implies that for some $c>0$ with probability $\Omega(1)$ we have that the number of isolated vertices in $\{A_j\}$ is at least $c|G_1$ and at least $c|G_2$ for $\{\widetilde{A_j}\}$. Denote this event by $\A$.

Furthermore, by Lemma
\ref{lem_subdeviation} we have
\begin{eqnarray}
&&\E(\sum_{j\geq1}|\A_j|^2+\sum_{j\geq1}|\C_j\cap G_1|^2)=O(|G_1|)\label{subcomponent1}\\
&&\E(\sum_{j\geq1}|\B_j|^2+\sum_{j\geq1}|\C_j\cap G_2|^2)=O(|G_2|)\label{subcomponent2}\\
&&\E(\sum_{j\geq1}|\tilde{A_j}|^2+\sum_{j\geq1}|\tilde{\C_j}\cap
G_1|^2)=O(|G_1|)\label{subcomponent3}\\
&&\E(\sum_{j\geq1}|\tilde{B_j}|^2+\sum_{j\geq1}|\tilde{\C_j}\cap
G_2|^2)=O(|G_2|)\label{subcomponent4}.
\end{eqnarray}

Now, we first assign spins to all components
except the isolated vertices in $\{A_j\}$ and $\{\widetilde{A_j}\}$ independently in
both chains. Let $M_1$, $N_1$ be the sum of spins in $G_1$ and $G_2$
respectively in first chain before assigning the rest of the spins, and similarly $\widetilde{M_1}$, $\widetilde{N_1}$
be the same for the second chain at this time. By
(\ref{subcomponent1}),(\ref{subcomponent2}),(\ref{subcomponent3}),(\ref{subcomponent4})
and Markov's inequality that we have
$$\big\{\mathcal{A},|M_1-\widetilde{M_1}|=O(\sqrt{|G_1|}),|N_1-\widetilde{N_1}|=O(\sqrt{|G_2|})\big\} $$ occurs with probability $\Omega(1)$. Then by Lemma
\ref{lem_CLT}, we can couple the sum of spins in both $G_1$ and
$G_2$ so that they are the same in both chains with probability $\Omega(1)$. This gives the required coupling of $(Y_1,Z_1)$ and
$(\widetilde{Y_1},\widetilde{Z_1})$. \qed


\section{Critical Case}
In this section, we prove that the mixing time for the Swendsen-Wang
dynamics in the critical case $c=2$ is of order $n^{1/4}$. This is
part (ii) of Theorem \ref{mainthm}.

Let $X_t$ and $Y_t$ be two magnetization chains such that $X_t$
starts from an arbitrary location and $Y_t$ starts from the
stationary distribution. To prove an upper bound of order $n^{1/4}$
to the mixing time we show that we can couple $X_t$ and $Y_t$ so
that they meet in time $O(n^{1/4})$ with
probability $\Omega(1)$. For a high level view of this coupling strategy we
refer the reader to Section \ref{roadmapcrit}.

Consider the following slight modification to the magnetization
chain $X_t$. Instead of choosing a random spin for each component
after the percolation step, we assign a positive spin to the largest
component and random spins for all other components. Let $X'_t$ be
the sum of spins at time $t$ (notice that we do {\em not} take
absolute values here), that is,
\begin{equation}\label{eqn_criticalstep}
X_{t+1}' \stackrel{d}{=} \max\{|\C_1^+(t)
|,|\C_1^-(t)|\}+\eps\min\{|\C_1^+(t) |,|\C_1^-(t)|\}+\sum_{j\geq
2}\epsilon_j |\C_j^+(t) | + \sum_{j\geq 2}\epsilon_j' |\C_j^-(t)|\,
,
\end{equation}
where as usual $\eps, \{\eps_j\}$ and $\{\epsilon_j'\}$ are
independent mean zero $\pm$ signs. This chain has state space
$[-n,n]$ and its absolute value is distributed as our original
chain. As a consequence, any upper bound on the mixing time of the
modified chain implies the same upper bound on the original chain.

The bulk of this section is devoted to the proof of the upper bound
on the mixing time (the corresponding lower bound is much easier to
prove and this is done in subsection \ref{critlowerbound}). To ease the notation, in
this section we will refer to this modified chain by $X_t$ and
$Y_t$. The only exception to this in this section is Theorem
\ref{pushdownuse} where another modification to the chain was
required for the proof.
%

The upper bound asserted in part (ii) of Theorem \ref{mainthm} will
follow immediately by the following two theorems. Though their
statement is almost identical, the difference in the starting point
$X_0$ give rise to completely different proof methods so we chose to
specify them as two separate theorems for convenience.

\begin{theorem}\label{allplus} Let $X_t$ and $Y_t$ be two SW
magnetization chains such that $X_0\geq n^{3/4}$ and
$Y_0\stackrel{d}{=}\pi$. Then we can couple $X_t$ and $Y_t$ so that
they meet each other within $O(n^{1/4})$ steps probability $\Omega(1)$.
\end{theorem}

\begin{theorem}\label{startzero}
Let $X_t$ and $Y_t$ be two SW magnetization chains such that $0\leq
X_0\leq n^{3/4}$ and $Y_0\stackrel{d}{=}\pi$. Then we can couple
$X_t$ and $Y_t$ so that they meet each other within $O(n^{1/4})$
steps with probability $\Omega(1)$.
\end{theorem}

\noindent{\bf{Proof of the upper bound of part (ii) of Theorem
\ref{mainthm}}:} Theorem \ref{allplus} and Theorem \ref{startzero}
give that for any $X_0\geq0$ we can couple $X_t$ and $Y_t$ so that
they meet within $O(n^{1/4})$ steps. If $X_0<0$, then by
(\ref{eqn_criticalstep}) and symmetry we have that
$$P(X_1\geq0)\geq\frac{1}{2},$$ so we may apply Theorem \ref{allplus} and Theorem
\ref{startzero} again. This shows that the mixing time of $X_t$ is
bounded above by $O(n^{1/4})$. Note that $|X_t|$ and the original
magnetization chain has the same distribution. Now Lemma
\ref{pathcoupling} gives the required upper bound and concludes the
proof..\qed

\subsection{Starting at the $[n^{3/4},n]$ regime: Proof of Theorem \ref{allplus}}

\begin{theorem} \label{thm_critical1}{\rm [Crossing and overshoot]} Let $X_t$ and $Y_t$ be two SW magnetization chains with $X_0\geq n^{3/4}$ and
$Y_0\stackrel{d}{=}\pi$. Put
$$ T = \min \big \{ t : X_t, Y_t\in
[A^{-1} n^{3/4}, A n^{3/4}] \and |X_t-Y_t|\leq h n^{5/8} \big \} \,
,$$ for some constant $h>0$ and large constant $A$. Then we can choose positive
constants $h,q,K$ depending only on $A$ such that
$$\P(T\leq Kn^{1/4})\geq q \,.$$
\end{theorem}

\begin{theorem} \label{localclt} {\rm [Local CLT]} For any constants $A>1$ and $h>0$, there exist constants $\delta=\delta(A,h)>0$ and $k=k(A,h)\in \mathbb{N}$
such that for any $x_0 \in [A^{-1} n^{3/4}, A n^{3/4}]$ and any $x
\in n+2\mathbb{Z}$ with $|x-x_0| \leq hn^{5/8}$, we have
$$ \prob ( X_k = x|X_0=x_0 ) \geq \delta n^{-5/8}.$$
\end{theorem}

\noindent{\bf{Proof of Theorem \ref{allplus}}:} By Theorem
\ref{thm_critical1}, the event $T\leq Kn^{1/4}$ occurs with
probability at least $q$. By Theorem \ref{localclt} and the strong
Markov Property we learn that there exist $\delta >0$ and $k\in
\mathbb{N}$ such that for any $x \in n+2\mathbb{Z}$ with
$|x-X_{T}|\leq hn^{5/8}$ and $|x-Y_{T}|\leq hn^{5/8}$, we have
$$\prob ( X_{T +k} = x \, \mid \, T \leq Kn^{1/4}) \geq \delta
n^{-5/8} \, ,$$ and $$\prob ( Y_{T +k} = x \,| \tau\leq Kn^{1/4})
\geq \delta n^{-5/8}.$$ Thus, for any such $x$ we can couple $X_t$
and $Y_t$ so that $X_{T +k}=Y_{T +k}=x$ with probability at least
$\delta n^{-5/8}$. We have at least $\frac{hn^{5/8}}{2}$ such $x$'s
so in this coupling we have that $X_{T +k}=Y_{T +k}$ with
probability at least $h\delta /2$. Lemma \ref{thm_coupling} concludes the proof. \qed \\


\subsubsection{\bf{Crossing and overshoot: Proof of Theorem
\ref{thm_critical1}}}

For any two magnetization chains $X_t$ and $Y_t$, define
$J_t=X_t-Y_t.$ Let $\tau$ be the first time the two chains cross
each other, i.e. \be\label{crossingtime}\tau:=\min\{t: \sign J_t\neq
\sign J_0\}.\ee

The following theorem implies Theorem \ref{thm_critical1}
immediately.

\begin{theorem}\label{thm_shoot}
Let $X_t$ and $Y_t$ be two independent magnetization SW chain with
$X_0\geq n^{3/4}$ and $Y_0\stackrel{d}{=}\pi$. There exists
positive constants $\delta,K,A$ and $h$ such that
$$\P\Big(\tau\leq Kn^{1/4}\,;\ X_{\tau-1},
Y_{\tau-1}\in [A^{-1} n^{3/4},\,A n^{3/4}]\,; J_{\tau-1}\leq
hn^{5/8}\Big)\geq \delta\, .$$
\end{theorem}

To prove Theorem \ref{thm_shoot} we will use the following results.

\begin{thm}\label{thm_statregime}
The stationary distribution $\pi$ of the modified magnetization
chain satisfies $$
\lim_{n\rightarrow\infty}\pi[a_1n^{3/4},\,a_2n^{3/4}]=\frac{1}{Z}\int_{a_1}^{a_2}
\exp(-\frac{1}{12}x^4)dx \, , $$ for any constants $a_2\geq a_1\geq
0$ where $Z=\int_0^{\infty}exp(-\frac{1}{12}x^4)dx$ is the
normalizing constant.
\end{thm}

\begin{lem}\label{gnptosw} For any constant $A>0$ there exists $N$ such that for all $n\geq N$ we have that if $X_0\in [A^{-1} n^{3/4},An^{3/4}]$, then the following hold:
\begin{itemize}
    \item[(i).] $-C n^{1/2} \leq \E X_1 -X_0\leq 0.$
    \item[(ii).] $\E |X_1-x_0|^k\leq Cn^{5k/8}$ for $k=2,3,4.$
    \item[(iii).]$\E \sum _{j \geq 1} |\C_j^-|^2 \geq
    cn^{5/4}.$
\end{itemize}
where $C=C(A)$ and $c=c(A)$ are constants.
\end{lem}

\begin{thm}\label{prop_nearhit}
Let $X_t$ and $Y_t$ be two independent magnetization chains with
$X_0, Y_0\in [b_1 n^{3/4}, b_2 n^{3/4}]$ for constants $b_2 > b_1 >
0$. Put $h=\frac{x_0-y_0}{n^{5/8}}$ and suppose that $h>0$ and that
$h=o(n^{1/8})$. Let $\tau$ be the crossing time of $X_t$ and $Y_t$
defined in (\ref{crossingtime}). Then there exist positive constants
$M$ and $\delta$ which only depend on $b_1$ and $b_2$ such that
$$\P(\tau\leq M h^2)\geq \delta \, .$$
\end{thm}

\begin{lemma}\label{lem_nojump} Let $X_t$ be a magnetization SW chain and $I=[a_1 n^{3/4},\,a_2
n^{3/4}]$ where $a_2>a_1>0$ are two constants. Let $h\in(0, a_1)$
and $\xi\in[0, a_1/4]$ be two constants. Then for any $b\in I$, we
have
$$\P\big(\sign (X_1-b)\neq \sign (X_0-b)\,\big| \, \, X_0>-\xi n^{3/4}\, , |X_0-b|\geq hn^{3/4}\big)\leq D n^{-1/3},$$ where $D=D(a_1, a_2, h,\xi)$ is a
constant.
\end{lemma}

\begin{thm}\label{lem_windowrest}
For any fixed constants $b_2>b_1>0$, $q<1$ and $K>0$, there exists a
constant $B=B(b_1,b_2,q,K)$ such that for every $X_0\in [b_1
n^{3/4},\,b_2 n^{3/4}]$, we have
$$\P\Big(X_t\leq Bn^{3/4}\ {\rm for\ all}\ t\in [0, Kn^{1/4}]\Big)\geq q.$$
\end{thm}

\begin{thm}\label{lem_hitting} Let $X_t$ be a magnetization SW chain with $X_0 > a n^{3/4}$ where $a>0$ is a constant.
Define $\tau_a = \min\{t\,:\, X_t \leq an^{3/4}\}$. Then for any
positive constant $b>0$ we have \be \label{eqn_crossmarkov} \P (
\tau_a
> bn^{1/4} )\leq \sqrt{{6 \over ab}} \, .\ee
\end{thm}

We begin by showing how these results imply the main theorem of this subsection. \\

\noindent{\bf{Proof of Theorem \ref{thm_shoot}:}} Let $a_1, K$ and
$C$ be three positive constants to be selected later. Define
$$\tau_1:=\min\{t: X_t< a_1 n^{3/4}\} \, ,$$
and define $\mathcal{A}$ to be the event that
\begin{enumerate}
\item $Y_0\in[\frac{a}{2}n^{3/4},an^{3/4}]$ and
\item $\tau_1\leq Kn^{1/4}$ and
\item $Y_{\tau_1}\geq a_1n^{3/4}$ and
\item $Y_t\leq Cn^{3/4}$ for all $t\leq Kn^{1/4}$.
\end{enumerate}
First we determine constants $a_1, \delta, K$ and $C$ so that
$\P(\mathcal{A})\geq \delta>0$. By Theorem \ref{thm_statregime},
there exists a constant $q>0$ such that $$\P\Big(\pi_n\in
[\frac{n^{3/4}}{2},\,n^{3/4}]\Big)\geq q.$$ By Theorem
\ref{thm_statregime} again, we can choose $a_1>0$ such that
$$\P\Big(\pi_n\in [-n,\,a_1 n^{3/4}]\Big)\leq \frac{q}{2}.$$ Since
$X_t$ and $Y_t$ are independent we have that
$Y_{\tau_1}\stackrel{d}{=}\pi_n$. Thus
$$\P\big(Y_0\in[\frac{n^{3/4}}{2}, n^{3/4}],
Y_{\tau_1}>a_1n^{3/4}\big)\geq\frac{q}{2}.$$  By Lemma
\ref{lem_hitting} there exists a constant $K=K(a_1,q)$ such that
\be\label{probA1}\P\big(Y_0\in[\frac{n^{3/4}}{2}, n^{3/4}],
Y_{\tau_1}>a_1n^{3/4},\tau_1\leq Kn^{1/4}\big)\geq\frac{q}{4}.\ee By
Lemma \ref{lem_windowrest}, there is a constant $C=C(K,q)$ such that
\be\label{probA2}\P\big(Y_t\leq C n^{3/4} \hbox{ for all } t\leq
Kn^{1/4}\Big)\geq1-\frac{q}{8}.\ee Combining (\ref{probA1}) and
(\ref{probA2}) shows that $\P(\mathcal{A})\geq\frac{q}{8}$. Note
that if $\mathcal{A}$ occurs, then
\be\label{crosstime}\tau\leq\tau_1\leq Kn^{1/4}.\ee

Next we show $\big\{J_{\tau-1}\leq \frac{a_1}{2}n^{3/4}\big\}\cap\A$
has positive probability. We do this by proving $J_{\tau-1}\leq
\frac{a_1}{2}n^{3/4}$ occurs with high probability on $\mathcal{A}$.
Note that $\big\{J_{\tau-1}\leq \frac{a_1}{2}n^{3/4}\big\}\cap\A$
implies $X_{\tau-1},Y_{\tau-1}\in [A^{-1}n^{3/4},An^{3/4}]$ for some
large constant $A$.

If $\{J_{\tau-1}>\frac{a_1}{2}n^{3/4}\}\cap\A$ occurs, then there
exists some $t\leq Kn^{1/4}$ such that $J_t>\frac{a_1}{2}n^{3/4}$
and $J_{t+1}<0$. This implies that there is a point $y\in
[\frac{a_1}{2}n^{3/4},(C+\frac{a_1}{2})n^{3/4}]$ with
$|X_t-y|\geq \frac{a_1}{4}n^{3/4}$ and $|Y_t-y|\geq \frac{a_1}{4}n^{3/4}$ and at least one of
$X_{t+1}$ and $Y_{t+1}$ crosses $y$. Suppose first that
$Y_{\tau-1}\geq -\xi n^{3/4}$ where $\xi$ is a small positive
constant. Then Lemma \ref{lem_nojump} and the union bound give that
\begin{equation}\label{eqn_shoottoobig}
\P\Big(\mathcal{A},J_{\tau -1}> \frac{a_1}{2}n^{3/4},
Y_{\tau-1}\geq-\xi n^{3/4} \Big)\leq Dn^{-1/3} Kn^{1/4}=o(1).
\end{equation}

Next suppose that $Y_{\tau-1}<-\xi n^{3/4}$. Then there is a
$t\in[0,Kn^{1/4}]$ such that $Y_t\leq-\xi n^{3/4}$. By
(\ref{eqn_criticalstep}), for any starting location, we have
$$\P(X_1<-\xi n^{3/4})\leq\P\Big(\eps\min\{|\C_1^+(t)
|,|\C_1^-(t)|\}+\sum_{j\geq 2}\epsilon_j |\C_j^+(t) | + \sum_{j\geq
2}\epsilon_j' |\C_j^-(t)|<-\xi n^{3/4}\Big).$$ By Theorem
\ref{gnp2moments} we have that
\be\label{momenta}\E\Big(\eps\min\{|\C_1^+(t)
|,|\C_1^-(t)|\}+\sum_{j\geq 2}\epsilon_j |\C_j^+(t) | + \sum_{j\geq
2}\epsilon_j' |\C_j^-(t)|\Big)^4 = O(n^{8/3}) \, ,\ee
so Markov's inequality gives that \be\label{nextstep1}\P(X_1<-\xi n^{3/4}) =
O(n^{-1/3}).\ee The union bound implies now that
\begin{equation}\label{eqn_shoottoobig1} \P\Big(\mathcal{A},J_{\tau
-1}> \frac{a_1}{2}n^{3/4}, Y_{\tau-1}< -\xi n^{3/4}\Big)\leq
O(n^{-1/3}) Kn^{1/4}=o(1) \, , \nonumber
\end{equation}
and so together with (\ref{eqn_shoottoobig}) we conclude that
$\{\mathcal{A},X_{\tau-1},Y_{\tau-1}\in[A^{-1}n^{3/4},An^{3/4}]\}$
occurs with probability $\Omega(1)$ for some constant $A$. We
denote this event by $\B$.

\vspace{.3cm} It remains to prove that $\{J_{\tau -1}\leq hn^{5/8}\}
\cap \B$ occurs with probability $\Omega(1)$ for some constant $h>0$.
Suppose first $J_{\tau-1}>n^{23/32}$. Notice that
$\{\B,J_{\tau-1}>n^{23/32}\}$ implies there is a $t<Kn^{1/4}$ such
that $X_t,Y_t\in[A^{-1}n^{3/4},An^{3/4}]$, $J_t>n^{23/32}$ and
$J_{t+1}<0$. This implies at least one of $X_t$ and $Y_t$ has to
make a huge jump of order at least $n^{23/32}$. By part (ii) of
Lemma \ref{gnptosw} with $k=4$, Markov's inequality and the union
bound we have
\begin{equation}\label{eqn_shootlarge}
\P(J_{\tau-1}>n^{23/32}, \mathcal{B}) \leq
O(n^{4(5/8-23/32)} n^{1/4})=o(1).
\end{equation}
To handle the case $J_{\tau-1}<n^{23/32}$, let
$$W_k=[2^{k}n^{5/8}, 2^{k+1}n^{5/8}],$$and consider the probability $\P(J_{\tau-1}\in W_k, \mathcal{B})$. Let $T_m$ be the first time
that $J_t\in W_k$ for $m$-th time and $$\A_m=\bigcap_{1\leq m'\leq
m}\{X_{T_{m'}},Y_{T_{m'}}\in[A^{-1}n^{3/4},A^{3/4}]\}.$$ Note that
$\A_m\in \F_{T_m}$. We have
\begin{equation}\label{eqn_shootsplit}
\P(J_{\tau-1}\in W_k, \mathcal{B}) \leq \sum_{m=1}^{\infty}
\P(T_m\leq\tau-1\, , J_{T_m +1}<0 \, , \mathcal{B}).
\end{equation} Notice that $\{T_m\leq\tau-1,J_{T_m +1}<0, \mathcal{B}\}$ implies that for all $m'\leq m$, we have $X_{T_{m'}}>a_1n^{3/4}$, $Y_{T_{m'}}<Cn^{3/4}$ and $|X_{T_{m'}}-Y_{T_{m'}}|\leq 2^{k+1}n^{5/8}$. This in particular implies that $X_{T_{m'}},Y_{T_{m'}}\in[A^{-1}n^{3/4},A^{3/4}]$. Hence $\{T_m\leq\tau-1\, , J_{T_m +1}<0\, , \mathcal{B}\}$ implies $\{T_m\leq\tau-1 \, , J_{T_m +1}<0\, , \A_m\}$. Also, by part (ii) of Lemma \ref{gnptosw} and Markov's inequality, we
have
$$\P(|X_{t+1}-X_t|\geq 2^{k-1} n^{5/8}\,|\,X_t=\Theta(n^{3/4}))\leq
\frac{C}{2^{4k}}.$$ The same inequality holds for $Y_t$ by the same
reason. Thus
\begin{equation*}
\P\Big(|J_{t+1}-J_t|\geq 2^{k}
n^{5/8}\Big|X_t,Y_t=\Theta(n^{3/4})\Big)\leq \frac{C}{2^{4k}}.
\end{equation*}
We now use the strong Markov property on the stopping time $T_m$ and
plug the above estimate in (\ref{eqn_shootsplit}) to get that
\begin{eqnarray}
\P(J_{\tau-1}\in W_k,\,\mathcal{B})&\leq&\sum_{m=1}^{\infty}
\P(T_m\leq\tau-1,J_{T_m +1}<0,
\A_m)\nonumber\\&\leq&\sum_{m=1}^{\infty}
\frac{C}{2^{4k}}\P(T_m\leq\tau-1,\,\A_m)\, .\label{eqn_shoot2}
\end{eqnarray}

If $\{T_m\leq\tau-1 \, , \A_m\}$ occurs, then for any $l\leq m$, we
have $T_{m-l}\leq\tau-1$ and
$X_{T_{m-l}},Y_{T_{m-l}}\in[A^{-1}n^{3/4},An^{3/4}]$ and most
importantly, the chains do not cross between time $T_{m-l}$ and
$T_m$, which is at least $l$ steps.  Now, let $M$ and $r$ be the
constants from Lemma \ref{prop_nearhit} and put $l=M2^{2k}$. The
strong Markov property on the stopping time $T_{m-M2^{2k}}$ and
Lemma \ref{prop_nearhit} gives that
\begin{eqnarray*}\P(T_m\leq \tau-1,\,\A_m)&\leq&
(1-r) \P(T_{m-M2^{2k}}\leq \tau-1,\,\A_{m-M2^{2k}}) \, .
\end{eqnarray*}
Applying this recursively gives that
$$\P(T_m\leq \tau-1,\,\A_m)\leq (1-r)^{[\frac{m}{M2^{2k}}]}.$$
Plugging this into (\ref{eqn_shoot2}), we get
$$\P(J_{\tau-1}\in
W_k,\,\mathcal{B})\leq \frac{C}{2^{4k}} \sum_{m=1}^{\infty}
(1-r)^{[\frac{m}{M2^{2k}}]}=\frac{C}{2^{2k}}.$$ Combining this and
(\ref{eqn_shootlarge}) we have for large enough $k_0$, $\{J_{\tau
-1}\leq 2^{k_0}n^{5/8},\mathcal{B}\}$ occurs with
probability $\Omega(1)$, which concludes the proof of the theorem.\qed \\

We now proceed with proving the statements we have used so far in
the proof of Theorem \ref{thm_shoot}. To prove Theorem
\ref{thm_statregime} we will use the following small lemmas.

\begin{lemma}[Simon and Griffiths (1973)]\label{isingregime}
Denote by $S_n$ the sum of spins for Ising model on the complete
graph. If the inverse temperature $\beta ={1\over n}$, then there
exists a random variable $X$ with density proportional to $\exp
(-\frac{1}{12}x^4)$ such that
$$\frac{S_n}{n^{3/4}}\stackrel{d}{\rightarrow} X,$$  as $n\rightarrow
\infty$.
\end{lemma}

\begin{corollary}\label{aymptotic}
Consider Ising model on the complete graph with inverse temperature
$\beta={1\over n}+O(\frac{1}{n^2})$. For any fixed constants
$a_2\geq a_1\geq 0$, we have
$$\lim_{n\rightarrow\infty}\P(|S_n|\in[a_1n^{3/4},a_2n^{3/4}])={1\over
A}\int_{a_1}^{a_2}exp(-{1\over 12}x^4)dx,$$ where
$A=\int_0^{\infty}exp(-{1\over 12}x^4)dx$ is the normalizing
constant.
\end{corollary}

\noindent{\bf{Proof of Corollary \ref{aymptotic}:}} By Lemma
\ref{isingregime} we have that the conclusion of the corollary holds for $\beta_1={1\over
n}$. Thus it suffices to prove that for any configuration $\sigma$ in which $|S_n(\sigma)| \in [a_1 n^{3/4}, a_2 n^{3/4}]$ we have
$$ \prob_\beta (\sigma) = (1+o(1)) \prob_{\beta_1}(\sigma) \, ,$$

Observe that on complete graph we have that
$$\sum_{u,v, u\neq v}\sigma(u)\sigma(v)=\frac{S_n^2-n}{2}.$$ Thus, for any $\sigma$ with $S_n(\sigma)\in[a_1n^{3/4}, a_2n^{3/4}]$, we have
\be\label{ratio}
\frac{\P_{\beta}(\sigma)}{\P_{\beta_1}(\sigma)}=\frac{e^{\beta(\frac{S_n^2-n}{2})}/Z(\beta)}{e^{\beta_1(\frac{S_n^2-n}{2})}/Z(\beta_1)}
= (1+o(1))\frac{Z(\beta_1)}{Z(\beta)} \, ,\ee
so it is enough to show $Z(\beta)=(1+o(1))Z(\beta_1)$. Indeed, Lemma
\ref{isingregime} implies that
$$\P_{\beta_1}(|S_n| \geq n^{7/8}) = o(1) \, ,$$ but for any configuration $\sigma$ with $|S_n(\sigma)| \leq n^{7/8}$ we have that
$$ e^{\beta  ( {S_n^2 - n \over 2} )} = (1+o(1)) e^{\beta_1 ({S_n^2 - n \over 2} )} \, ,$$
and the assertion follows.\qed

\begin{lemma}\label{nonnegative}
Let $\pi_n$ be the stationary distribution of the modified
magnetization SW chain $X_t$, then we have
$$\lim_{n\rightarrow\infty} \pi_n[-\infty,0]=0.$$
\end{lemma}

\noindent{\bf{Proof of Lemma \ref{nonnegative}:}} Recall that SW
dynamics with parameter $p$ has stationary distribution of Ising
model with $p=1-e^{-2\beta}$. Plugging in $p={2\over n}$ we get
$\beta={1\over n}+O({1\over n^2})$. By Corollary \ref{aymptotic}, for
any $\eps>0$, there exists constant $b_1$ and $b_2$ such that
$0<b_1<b_2$ and
$$\pi_n \big ([b_1n^{3/4},b_2n^{3/4}] \cup
[-b_2n^{3/4},-b_1n^{3/4}]\big )>1-\eps.$$ By definition of
stationarity, for any set $S$ we have

\begin{equation}\label{stationary}
\sum_{y\in [-n,n]}\pi_n(y)\P(y,S)=\pi_n(S).
\end{equation} Put $S=[-n,0]$ and
denote $\pi_n[-n,0]$ by $\delta_n$. For any $X_0$ we have
$$\P(X_1\leq0)\leq 1/2 \, ,$$
by symmetry. For $X_0\in[b_1n^{3/4},b_2n^{3/4}]$ Lemma \ref{lem_nojump} gives that
$$P(X_1<0) \leq Dn^{-1/3} \, .$$ Plugging these into (\ref{stationary}),
we have $$\delta_n=\sum_{y\in
[-n,n]}\pi_n(y)\P(y,S)\leq\frac{1}{2}\delta_n+\eps+Dn^{-1/3},$$ which
gives $$\delta_n\leq2(\eps+Dn^{-1/3}) \, ,$$ concluding the proof. \qed \\

\noindent{\bf{Proof of Theorem \ref{thm_statregime}:}} Directly follows from Corollary \ref{aymptotic}
and Lemma \ref{nonnegative}.\qed \\

The following is an easy estimate which use frequently to show that
the main contribution from the first term of
(\ref{eqn_criticalstep}) comes from the $|\C_1^+|$ element rather
than the $|\C_1^-|$ element.
\begin{prop} \label{c1plus} If $X_0 \geq C n^{2/3} \log^2 n$ for some large constant $C$, then
$$ \prob(|\C_1^-| \geq |\C_1^+|) \leq  O(e^{-c\log^2 n}) \, .$$
\end{prop}
\noindent {\bf Proof.} By our condition on $X_0$ we have that
$|\C_1^+|$ is distributed as the size of the largest component in a
supercritical random graph $G(m,p)$ with $m={n+X_0 \over 2}$ and
$p={1+\eps \over m}$ with $\eps = X_0/n = \Omega(n^{-1/3} \log^2
n)$. Theorem \ref{gnpout2} gives that
$$ \prob(|\C_1^+| \geq c n^{2/3} \log^2 n) \geq 1 - C e^{-c\log^2 n} \, ,$$ for some small $c>0$. On the other hand $|\C_1^-|$ is distributed as a subcritical random graph. Theorem 1 of \cite{NP1} gives that
$$ \prob(|\C_1^-| \geq c n^{2/3} \log^2 n ) \leq C e^{-c\log^2 n} \, ,$$
which finishes the proof.\qed

\begin{lemma}\label{driftestimate} If $X_t \geq C n^{2/3} \log
n$ for some large constant $C$, then \be\label{eqn_firstdrift}
\E[X_{t+1} \mid X_t ] \leq X_t\Big (1-{X_t \over 6n} \Big ). \ee
\end{lemma}
\noindent{\bf Proof.} By (\ref{eqn_criticalstep}) we have $\E[
X_{t+1} \mid X_t ] =\E \big[|\max\{|\C_1^+|,|\C_1^-|\} |\big|
X_t\big]$, hence Proposition \ref{c1plus} gives that
$$ \E[ X_{t+1} \mid X_t ] = \E |\C_1^+| + O(e^{-c\log^2 n}) \, .$$
Thus, Theorem \ref{supergnp1} yields that
$$\E[ X_{t+1} \mid X_t ] \leq {2 \frac{X_t}{n}\frac{n+X_t}{2}} - {\frac{7}{3} \frac{X_t^2}{n^2} \frac{n+X_t}{2}} + O(e^{-c\log^2 n}) \leq X_t \Big(1- {X_t \over 6n}\Big) \, ,$$
when $n$ is large enough. \qed \\


\noindent{\bf Proof of Lemma \ref{gnptosw}:} As in the previous
proof we have
$$ \E X_1 = \E |\C_1^+| + O(e^{-c\log^2 n}) \, .$$ Since
$\eps=\frac{x_0}{n}=\Theta(n^{-1/4})$ Theorem \ref{supergnp1} gives
that
\begin{eqnarray*}
\E |\C^+_1 | &=& 2\frac{x_0}{n}\frac{n + x_0}{2} - {8 \over 3}\Big(\frac{x_0}{n}\Big)^2\frac{n + x_0}{2} + O\Big(\Big(\frac{x_0}{n}\Big)^3\frac{n + x_0}{2}\Big)\\
&=& x_0 -\frac{x_0^2}{3n} +O\Big(\frac{x_0^3}{n^2}\Big) = x_0 -\frac{x_0^2}{3n} + O(n^{1/4}) \, ,
\end{eqnarray*} which gives part (i) of the lemma since $x_0 \in [A^{-1} n^{3/4}, An^{3/4}]$.
We now prove part (ii). For $k=2,3,4$, by (\ref{eqn_criticalstep})
and Jensen's inequality we have that
\begin{eqnarray}\label{eqn_stepsplit}
\E |X_1-x_0|^k &=&\E\big||\C_1^+|-x_0+ \sum_{j\geq 2}\epsilon_j
|\C_j^+ | + \sum_{j\geq 1}\epsilon_j'
|\C_j^-|\big|^k+\Theta(e^{-cn^{1/8}})\nonumber\\&\leq&
2^{k-1}\Big(\E\Big||\C_1^+ |-x_0\Big|^k+ \E\Big|\sum_{j\geq
2}\epsilon_j |\C_j^+ | + \sum_{j\geq 1}\epsilon_j'
|\C_j^-|\Big|^k\Big).
\end{eqnarray}
Theorem \ref{supergnpmoments} now gives that
\begin{eqnarray*}
\E\Big||\C_1^+ |-2\frac{x_0}{n}\frac{n+x_0}{2}\Big|^k &\leq&
C\Big(\frac{n+x_0}{2}\Big/\frac{x_0}{n}\Big)^{k/2}\\
&\leq& C\Big(\frac{n^2}{x_0}\Big)^{k/2}\leq O(n^{5k/8}).
\end{eqnarray*}
Another application of Jensen's inequality gives that
\begin{eqnarray}\label{eqn_giantk}
\E\big||\C_1^+|-x_0\big|^k&=&\E\big|\big(|\C_1^+|-2\frac{x_0}{n}\frac{n+x_0}{2}\big)+
\big (2\frac{x_0}{n}\frac{n+x_0}{2}-x_0\big)\big|^k\nonumber
\\&\leq&
2^{k-1}\Big(O(n^{5k/8})+\Big(\frac{x_0^2}{n}\Big)^k\Big)\leq
O(n^{5k/8}).
\end{eqnarray}
To bound the rest of (\ref{eqn_stepsplit}), notice that by Holder's
inequality, we only need to consider the case $k=4$. We have
\begin{eqnarray*}\label{eqn_kk4}
\E\Big|\sum_{j\geq 2}\epsilon_j |\C_j^+ | + \sum_{j\geq
1}\epsilon_j' |\C_j^-|\Big|^4 &\leq& \sum_{j\geq 2}\E |\C_j^+ |^4 +
\sum_{j\geq 1}\E |\C_j^-|^4 +  \Big(\sum_{j\geq 2}\E
|\C_j^+|^2\big)\big(\sum_{j\geq 1}\E |\C_j^-|^2\Big)  \\  &+& \E
\sum_{i,j\geq2,i\neq j}|\mathcal{C}_i^+|^2|\mathcal{C}_j^+|^2+\E
\sum_{i,j\geq1,i\neq j}|\mathcal{C}_i^-|^2|\mathcal{C}_j^-|^2 \, .
\end{eqnarray*}
By Theorem \ref{supergnpmoments} we have
\begin{equation*}
\sum_{j\geq 2}\E |\C^+_j|^2 \leq C_2 n \Big(\frac{x_0}{n}\Big)^{-1}=
O(n^{5/4})
\end{equation*} and
\begin{equation*}
\sum_{j\geq 2}\E |\C^+_j|^4 \leq C_4 n \Big(\frac{x_0}{n}\Big)^{-5}=
O(n^{9/4}).
\end{equation*}
By Theorem \ref{subgnpmoments}, we have
\begin{equation*}
\sum_{j\geq 1}\E |\C^-_j|^2 \leq C_2 n \Big(\frac{x_0}{n}\Big)^{-1}
= O(n^{5/4})
\end{equation*} and
\begin{equation*}
\sum_{j\geq 1}\E |\C^-_j|^4 \leq C_4 n \Big(\frac{x_0}{n}\Big)^{-5}
= O(n^{9/4}).
\end{equation*}
These together with Theorem \ref{gnp2moments} to handle the cross
terms finishes the proof of part (ii) of the lemma. Part (iii)
follows immediately by Theorem \ref{subgnpmoments},
$$\E\sum_{j\geq 1}|\C^-_j |^2 \geq c_2\frac{n-x_0}{2}\frac{n}{x_0}\geq c_2\frac{n}{4}A^{-1}n^{1/4}\geq cn^{5/4} \, .$$
\qed

\begin{lemma}\label{lem_moment}
Let $X$ be a real valued random variable with $\E X=0$ and $\E
X^2\geq h^2$ and $\E X^4\leq b h^4$ where $b\geq 1$. Then for any
$\rho\in[0,1]$ we have
$$ \P(X\leq -\rho h)\geq \frac{(1-\rho ^2)^2}{2b}\, .$$
\end{lemma}
\noindent{\bf Proof of Lemma \ref{lem_moment}:} By Cauchy-Schwartz
$$\E [X^2 \1_{\{X^2\geq \rho^2h^2\}}] \leq \sqrt{\E X^4 \E
\1_{\{X^2\geq\rho^2 h^2\}}} \leq \sqrt{b h^4\P(X^2\geq\rho^2 h^2)}
\, .$$ Hence,
$$ h^2\leq\E X^2 \leq \rho^2 h^2 +\E [X^2 \1_{\{X^2\geq \rho^2h^2\}}] \leq \rho^2 h^2 + \sqrt{b h^4\P(X^2\geq\rho^2 h^2)} \, .$$
We conclude that
$$ \P(|X|\geq \rho h)\geq \frac{(1-\rho ^2)^2}{b} \, ,$$
and the assertion follows by symmetry since $\P(X\leq-\rho
h)=\P(-X\leq-\rho h)$. \qed\\

The following will be used in the proof of Theorem
\ref{prop_nearhit}.

\begin{theorem}\label{cor_1over2}
Let $X_t$ be a magnetization chain with $X_0\in [b_1n^{3/4},b_2
n^{3/4}]$ where $b_2>b_1>0$ are two constants. Let $\tau_1$ be the
first time that
$X_t\not\in[\frac{b_1}{2}n^{3/4},\,(b_2+\frac{b_1}{2})n^{3/4}]$.
Then there exists a constant $C=C(b_1,b_2)>0$ such that for all
constant $\delta>0$ we have
$$\P (\tau_1 \leq \delta n^{1/4})\leq C\delta^2.$$
\end{theorem}

\noindent{\bf{Proof of Theorem \ref{cor_1over2}}} Denote by $I$ the
interval $[\frac{b_1}{2}n^{3/4},\,(b_2+\frac{b_1}{2})n^{3/4}]$. Part
(ii) of Lemma \ref{gnptosw} gives
\begin{equation}\label{eqn_only234}
\E\Big[(X_{(t+1)\wedge\tau_1}-X_{t\wedge\tau_1})^k\Big|
\F_t\Big]\leq Cn^{5k/8}
\end{equation}
for $k=2,3,4$. Define
$$Z:=X_{(t+1)\wedge\tau_1}-X_{t\wedge\tau_1}-(\E
X_{(t+1)\wedge\tau_1}-\E X_{t\wedge\tau_1}).$$ Note that $\big|\E
X_{(t+1)\wedge\tau_1}-\E X_{t\wedge\tau_1}\big|\leq Cn^{1/2}$ by
part (i) of Lemma \ref{gnptosw}, hence
\begin{equation}\label{eqn_234}
\E\Big[Z^k\Big| \F_t\Big] \leq C n^{\frac{5k}{8}}
\end{equation}
for $k=2,3,4$. Also, for $k=1$, part (i) of Lemma \ref{gnptosw}
gives that
\begin{equation}\label{eqn_only1}
\E[Z |\F_t]\leq C\sqrt{n}.
\end{equation}
Denote $$f(t)=\Big(\E [X_{t\wedge\tau_1}-\E
X_{t\wedge\tau_1}]^4\Big)^{1/2}.$$ Note that
\begin{equation}\label{eqn_Mexpansion1}
f(t+1)^2 =\E [X_{(t+1)\wedge\tau}-\E
X_{(t+1)\wedge\tau_1}]^4=\E\Big[(X_{t\wedge\tau_1}-\E
X_{t\wedge\tau_1})+Z\Big]^4.
\end{equation}
For $k=1,2,3,4$, we have
\begin{eqnarray*}
\E\Big[\Big(X_{t\wedge\tau_1}-\E
X_{t\wedge\tau_1}\Big)^{4-k}Z^k\Big]&=&\E\Big(\E\Big[\Big(X_{t\wedge\tau_1}-\E
X_{t\wedge\tau_1}\Big)^{4-k}Z^k\Big|\F_t\Big]\Big)\\
&=&\E\Big[\Big(X_{t\wedge\tau_1}-\E
X_{t\wedge\tau_1}\Big)^{4-k}\E[Z^k|\F_t]\Big].
\end{eqnarray*}
H\"{o}lder's inequality implies that
\begin{equation}\label{eqn_a2341}
\E\Big[\Big(X_{t\wedge\tau_1}-\E
X_{t\wedge\tau_1}\Big)^{4-k}Z^k\Big]\leq C
n^{\frac{5k}{8}}f(t)^{\frac{4-k}{2}},
\end{equation}
for $k=2,3,4$ and by (\ref{eqn_only1})
\begin{equation}\label{eqn_aonly11}
\E\Big[\Big(X_{t\wedge\tau_1}-\E X_{t\wedge\tau_1}\Big)^3 Z\Big]\leq
C \sqrt{n}f(t)^{\frac{3}{2}}.
\end{equation}
Expanding the right hand side of (\ref{eqn_Mexpansion1}) and
plugging (\ref{eqn_a2341}) and (\ref{eqn_aonly11}) into it, we get
\be\label{eqn_ft1} f(t+1)^2\leq
f(t)^2+C\sqrt{n}f(t)^{3/2}+Cn^{5/4}f(t)+Cn^{15/8}f(t)^{1/2}+Cn^{5/2}.
\end{eqnarray} Comparing the right hand side of (\ref{eqn_ft1}) with \be\label{eqn_ft22}\Big(f(t)+Cn^{1/2}f(t)^{1/2}+Cn^{5/4} \Big)^2,\ee we find that the first, second, third and fifth term of
(\ref{eqn_ft1}) is dominated by expanding (\ref{eqn_ft22}). For the
forth term, if $f(t)=O(n^{5/4})$, then it is dominated by
$(Cn^{5/4})^2$. Otherwise it is dominated by $Cn^{5/4}f(t)$. The
conclusion is that
\be\label{eqn_recursive}f(t+1)^2\leq\Big(f(t)+Cn^{1/2}f(t)^{1/2}+Cn^{5/4}
\Big)^2.\ee Thus, if $f(t)=O(n^{3/2})$, then we have
\begin{equation}\label{eqn_shoot4xstep}
f(t+1)\leq f(t)+Cn^{5/4}.
\end{equation} Since $f(0)=0$, by iterating (\ref{eqn_shoot4xstep}) we get that $f(t)\leq
Ctn^{5/4}$ for all $t\leq \delta n^{1/4}$ where $\delta>0$ is a
constant. Put $t=\delta n^{1/4}$. Markov's inequality gives that
\begin{equation}\label{eqn_forthconstra}
\P(| X_{\delta n^{1/4}\wedge\tau_1}-\E X_{\delta
n^{1/4}\wedge\tau_1}|\geq\frac{b_1}{4}n^{3/4})\leq
\frac{(C\delta)^2}{\big(\frac{b_1}{4}\big)^4}.
\end{equation}
By part (i) of Lemma \ref{gnptosw} we have that $$|\E X_{\delta
n^{1/4}\wedge\tau_1}-X_0|\leq C\delta n^{3/4}.$$ Thus, for small
enough $\delta$ we have
\begin{equation}\label{eqn_forthconstra2}
\P\Big(\Big| X_{\delta
n^{1/4}\wedge\tau_1}-x_0\Big|\leq\frac{b_1}{2}n^{3/4}\Big)\geq
1-C\delta^2,
\end{equation} which means that $X_t$ has not jumped out of the window $I$ within
$\delta n^{1/4}$ steps with probability at least
$1-C\delta^2$. \qed \\

\noindent{\bf Proof of Theorem \ref{prop_nearhit}:} Recall that
$J_t=X_t-Y_t$. Let $M$ be a large constant that will be chosen
later. Assume without loss of generality that $J_0 \geq 0$, we will prove that $J_{Mh^2}$ is negative with probability $\Omega(1)$, which implies the theorem. Denote by $I$ the
interval $[\frac{b_1}{2}n^{3/4}, (b_2+\frac{b_1}{2})n^{3/4}]$ and
define
$$\tau_1=\min\{t: X_t \not \in I \hbox{ or } Y_t \not \in I \} \,
.$$ We will prove out claim by precisely estimating the first,
second and forth moment of $J_{Mh^2\wedge\tau_1}$ and then apply
Lemma \ref{lem_moment} to $J_{Mh^2\wedge\tau_1}-\E
J_{Mh^2\wedge\tau_1}$. We start with first moment estimate. By part
(i) of Lemma \ref{gnptosw} and the optional stopping theorem we get
\begin{equation}\label{eqn_shoot1xstep}
\E X_{t\wedge\tau_1}-C\sqrt{n}\leq\E X_{(t+1)\wedge\tau_1}\leq \E
X_{t\wedge\tau_1}
\end{equation} for some constant $C=C(b_1,b_2)>0$.
Applying (\ref{eqn_shoot1xstep}) recursively gives that
\begin{equation}\label{eqn_shoot1x}
X_0-CMh^2 \sqrt{n}\leq\E X_{Mh^2\wedge\tau_1}\leq X_0.
\end{equation}
The same formula holds for $Y_t$, hence
\begin{equation}\label{eqn_shoot1j}
\E J_{Mh^2\wedge\tau_1}\leq hn^{5/8}+CMh^2n^{1/2}.
\end{equation}

We proceed with the second moment estimate. Notice that if $X_0\in
I$, we have that
\begin{eqnarray*}
\E(X_1-\E(X_1|\mathcal{F}_0))^2&=&\E\Big[\max\{|\C_1^+|,|\C_1^-|\}-\E\max\{|\C_1^+|,|\C_1^-|\}\\&+&\eps\min\{|\C_1^+|,
|\C_1^-|\}+\sum_{j\geq 2}\epsilon_j |\C_j^+| +\sum_{j\geq 2}\epsilon_j' |\C_j^-|\Big]^2
\end{eqnarray*} by (\ref{eqn_criticalstep}). In the $n^{3/4}$ regime, we have $\P(|C^-_1|\geq|C^+_1|)=O(e^{-c\log^2 n})$ by Proposition \ref{c1plus}, hence $$\E(X_1-\E(X_1|\mathcal{F}_0))^2\geq (1-Ce^{-c\log^2 n}))\E\sum_{j\geq1}|C^-_j|^2\geq c_1n^{5/4}\, ,$$ by part (iii) of Lemma (\ref{gnptosw}). Also, by part (ii) of Lemma \ref{gnptosw} we have that $$\E(X_1-\E(X_1|\mathcal{F}_0))^2\leq Cn^{5/4}.$$
Now let
$$A_t=\sum_{i=0}^{t-1}(X_{(i+1)\wedge\tau_1}-X_{i\wedge\tau_1})-\E(X_{(i+1)\wedge\tau_1}-X_{i\wedge\tau_1}|\mathcal{F}_i)$$
and $$B_t=X_0-\E
X_{t\wedge\tau_1}+\sum_{i=0}^{t-1}\E(X_{(i+1)\wedge\tau_1}-X_{i\wedge\tau_1}|\mathcal{F}_i).$$
Then it is easy to verify that $A_t+B_t=X_{t\wedge\tau_1}-\E
X_{t\wedge\tau_1}$. Moreover, since the martingale increments are
orthogonal we have that $$\E
A_t^2=\sum_{i=0}^{t-1}\E(X_{(i+1)\wedge\tau_1}-\E(X_{(i+1)\wedge\tau_1}|\mathcal{F}_i))^2.$$
Since $h=o(n^{1/8})$ Theorem \ref{cor_1over2} gives that
$$\P(\tau_1 \leq M h^2) = o(1) \, .$$ This implies that
$$cMh^2n^{5/4}\leq\E A_{Mh^2}^2\leq CMh^2n^{5/4}.$$ By (\ref{eqn_shoot1x}) and
part (i) of Lemma \ref{gnptosw}, we get that $|B_t|\leq Ctn^{1/2}$.
This gives $$\E B_t^2\leq Ct^2n.$$ Cauchy-Schwarz inequality gives
$$\E|A_{Mh^2}B_{Mh^2}|\leq CMh^2n^{9/8}.$$
Thus, we have $$\Var X_{Mh^2\wedge\tau_1}=\E A_{Mh^2}^2+\E
B_{Mh^2}^2+2\E A_{Mh^2}B_{Mh^2}\geq (c-o(1))Mh^2n^{5/4}.$$ The
same estimates hold for $Y_t$. Since $X_t$ and $Y_t$ are independent we have \be\label{eqn_shoot2j}\Var
J_{Mh^2\wedge\tau_1}\geq c_1Mh^2n^{5/4}.\ee

For the fourth moment estimate, by (\ref{eqn_shoot4xstep}) we have
$$\E [X_{Mh^2\wedge\tau_1}-\E X_{Mh^2\wedge\tau_1}]^4\leq
(Mh^2Cn^{5/4})^2$$ and $$\E [Y_{Mh^2\wedge\tau_1}-\E
Y_{Mh^2\wedge\tau_1}]^4\leq (Mh^2Cn^{5/4})^2.$$ By the Jensen's
inequality, we get
\begin{equation}\label{eqn_shoot4j}
\E [J_{Mh^2\wedge\tau_1}-\E
J_{Mh^2\wedge\tau_1}]^4\leq16(Mh^2Cn^{5/4})^2.
\end{equation}
Putting (\ref{eqn_shoot2j}) and (\ref{eqn_shoot4j}) together, taking
$\rho=\frac{1}{\sqrt{Mc_1}}$ and using Lemma \ref{lem_moment}, we
get $$\P(J_{Mh^2\wedge\tau_1}-\E J_{Mh^2\wedge\tau_1}\leq
-hn^{5/8})\geq \delta,$$ where $\delta>0$ is a constant. Combining
this with (\ref{eqn_shoot1j}), we get
\begin{equation}\label{eqn_willhit}
\P \Big(J_{Mh^2\wedge\tau_1}\leq 0\Big)\geq \delta.
\end{equation} Here we choose $M$
so that $Mc_1\geq 2$, concluding the proof.\qed \\

\noindent{\bf{Proof of Lemma \ref{lem_nojump}:}} If $X_0\leq
b-hn^{3/4}$, then assume first $X_0\in[\xi n^{3/4}, b- h n^{3/4}]$.
In this regime, by part (ii) of Lemma \ref{gnptosw} with $k=4$ and
Markov's inequality we have
$$\P\big(\sign (X_1-b)\neq \sign (X_0-b)\big)\leq
\frac{Cn^{5/2}}{h^4n^3}=O(n^{-1/2}).$$ Assume now $X_0\in[-\xi
n^{3/4}, \xi n^{3/4}]$. Recall the distribution of $X_1$ in
(\ref{eqn_criticalstep}) and that $b\geq a_1n^{3/4}$ and $\xi\leq
a_1/4$. If $X_1\leq a_1n^{3/4}$, then either
$$\max\{|\C_1^+|,|\C_1^-|\}>2\xi n^{3/4},$$ or $$\eps\min\{|\C_1^+(t) |,|\C_1^-(t)|\}+\sum_{j\geq 2}\epsilon_j
|\C_j^+(t) | + \sum_{j\geq 2}\epsilon_j' |\C_j^-(t)|\geq
\frac{a_1}{2}n^{3/4}.$$ By Theorem \ref{gnpout2} and monotonicity of
$|\C_1|$, we have
$$\P(\max\{|\C_1^+|,|\C_1^-|\}>2\xi n^{3/4})\leq Ce^{-cn^{1/8}}.$$ By Theorem
\ref{gnp2moments} and Markov's inequality, we have
\begin{eqnarray*}
\P(\eps\min\{|\C_1^+(t) |,|\C_1^-(t)|\}+\sum_{j\geq 2}\epsilon_j
|\C_j^+(t) | + \sum_{j\geq 2}\epsilon_j' |\C_j^-(t)|\geq
\frac{a_1}{2}n^{3/4})= O(n^{-1/3}).\end{eqnarray*} Thus, we have
$$\P(X_1\geq b)= O(n^{-1/3}).$$

If $X_0\geq b+hn^{3/4}$, then assume first $X_0\in[b+h^{3/4},
Bn^{3/4}]$ for some large constant $B$. By part (ii) of Lemma
\ref{gnptosw} with $k=4$ and Markov's inequality, we have
$$\P(X_1\leq b)\leq \frac{Cn^{5/2}}{h^4n^3}=O(n^{-1/2}).$$ Assume
$X_0\geq Bn^{3/4}$. If $X_1\leq bn^{3/4}$, then either
$$\max\{|\C_1^+|,|\C_1^-|\}\leq\frac{B}{2} n^{3/4},$$ or $$\eps\min\{|\C_1^+(t)
|,|\C_1^-(t)|\}+\sum_{j\geq 2}\epsilon_j |\C_j^+(t) | + \sum_{j\geq
2}\epsilon_j' |\C_j^-(t)|\leq -(\frac{B}{2}-a_2)n^{3/4}.$$ By
Theorem \ref{gnpout2}, we have
$$\P(\max\{|\C_1^+|,|\C_1^-|\}\leq\frac{B}{2} n^{3/4})\leq
Ce^{-cn^{1/8}}.$$ By Theorem \ref{gnp2moments} and Markov's
inequality we have
\begin{eqnarray*}\P(\eps\min\{|\C_1^+(t) |,|\C_1^-(t)|\}+\sum_{j\geq 2}\epsilon_j
|\C_j^+(t) | + \sum_{j\geq 2}\epsilon_j' |\C_j^-(t)|\leq
O(n^{-1/3}).\end{eqnarray*} Thus, we have
$$\P(X_1\leq b)=O(n^{-1/3})\, .$$\qed

\noindent{\bf{Proof of Theorem \ref{lem_windowrest}:}} Denote by $I$
the interval $[b_1 n^{3/4},b_2 n^{3/4}]$ and let $B$ be a large
constant to be chosen later. For any $X_0\in I$, define
$$A_{t,X_0} =\P(X_t\ \rm{exceeds\ }Bn^{3/4}\ \rm{within\ }t\ \rm{steps} \,\, |\,\, X_0).$$
Let
$$A_t=\max_{X_0\in I} A_{t,X_0}.$$
Then $A_t$ is increasing in $t$. Let $$\tau=\min\{t:
X_t\not\in[\frac{b_1}{2}n^{3/4},Bn^{3/4}]\}\, .$$ Then
$X_{t\wedge\tau}$ is a supermartingale by part (i) of Lemma
\ref{gnptosw}. Thus we have \be\label{upperboundre}\E
X_{Kn^{1/4}\wedge\tau}\leq b_2n^{3/4}.\ee
For simplicity denote $g(B)=\max_{X_0\in I} \P(X_{Kn^{1/4}\wedge\tau}\geq B
n^{3/4}\mid X_0)$. We get from the above estimate and (\ref{upperboundre})
that $g(B)\rightarrow 0$ as $B\rightarrow \infty$. For all $X_0\in
I$ and $t\leq Kn^{1/4}$, we have
\begin{equation}\label{eqn_split}
A_{t,X_0}\leq g(B,X_0)+ \P\Big(X_{Kn^{1/4}\wedge\tau}\leq \frac{b_1}{2}
n^{3/4}, X_t\ \rm{exceeds}\ Bn^{3/4}\ \rm{before}\ t\Big).
\end{equation} Denote $$\A=\{X_{Kn^{1/4}\wedge\tau}\leq \frac{b_1}{2} n^{3/4}, X_t\ \rm{exceeds}\
Bn^{3/4}\ \rm{before}\ t\}.$$ Let $\tau_1$ be the exit time of
$[\frac{b_1}{2}n^{3/4},(b_2+\frac{b_1}{2})n^{3/4}]$. By Theorem
\ref{cor_1over2}, we have
$$\P(\tau_1 > \delta n^{1/4})\geq 1-C\delta^2,$$ for any sufficiently  small constant $\delta>0$. On the event $\{\tau_1 > \delta n^{1/4}\}$, there
are three cases:
\begin{itemize}
\item[(i)]$X_{\delta n^{1/4}}\in[\frac{b_1}{2}n^{3/4},b_1n^{3/4}]$.
\item[(ii)]$X_{\delta n^{1/4}}\in[b_1n^{3/4},b_2n^{3/4}]$.
\item[(iii)]$X_{\delta n^{1/4}}\in[b_2n^{3/4},(b_2+\frac{b_1}{2})n^{3/4}]$.
\end{itemize}
For case (ii), by the Markov property at time $\delta n^{1/4}$, we
have that
$$\P(\A\, \big |\, \tau_1 > \delta n^{1/4}, X_{\delta
n^{1/4}}\in[b_1n^{3/4},b_2n^{3/4}])\leq A_{t-\delta n^{1/4}}.$$ For
case (i), define
$$T=\min\{t > \delta n^{1/4} : X_t\in[b_1n^{3/4},b_2n^{3/4}]\}\, .$$ By monotonicity of $A_t$ and the strong Markov property on $T$ we
have
$$\P\big(\A\,\big |\,\tau_1
> \delta n^{1/4}, X_{\delta n^{1/4}}\in[\frac{b_1}{2}n^{3/4},b_1n^{3/4}], T<t\big)\leq A_{t-\delta n^{1/4}}.$$ The event $\{\A, \tau_1
> \delta n^{1/4}, X_{\delta n^{1/4}}\in[\frac{b_1}{2}n^{3/4},b_1n^{3/4}], T\geq t\}$ implies that there exists $t\leq Kn^{1/4}$ such that $X_t<\frac{b_1}{2}n^{3/4}$ and
$X_{t+1}>b_2n^{3/4}$. By Lemma \ref{lem_nojump} and the union bound,
this happens with probability at most
$$Dn^{-\frac{1}{3}}Kn^{1/4}=O(n^{-1/12}).$$ For case (iii), the event
$\big\{\A, \tau_1 > \delta n^{1/4}, X_{\delta
n^{1/4}}\in[b_2n^{3/4},(b_2+\frac{b_1}{2})n^{3/4}]\big\}$ implies
that $X_t$ first goes below $\frac{b_1}{2}n^{3/4}$ and then goes
above $Bn^{3/4}$. Let $$T'=\min\{t: t>\tau, X_t\in I \}.$$ By
monotonicity of $A_t$ and the strong Markov property on $T'$, we
obtain
$$\P\big(\mathcal{A}\,\big| \, \tau_1
> \delta n^{1/4}, X_{\delta n^{1/4}}\in[b_2n^{3/4},(b_2+\frac{b_1}{2})n^{3/4}], T'<t\big)\leq A_{t-\delta
n^{1/4}}.$$ By similar argument in case (ii), we have
$$\P\big(\mathcal{A}, \tau_1
> \delta n^{1/4},X_{\delta
n^{1/4}}\in[b_2n^{3/4},(b_2+\frac{b_1}{2})n^{3/4}], T'\geq t\big)=
O(n^{-1/12}).$$ Summing up the above estimates, we obtain
\be\label{estimate1}\P(\A, \tau_1 > \delta n^{1/4})\leq A_{t-\delta
n^{1/4}}+O(n^{-1/12}).\ee

On the event $\{\A, \tau_1 \leq \delta n^{1/4}\}$, which happens
with probability at most $C\delta^2$, there are two cases to
consider:
\begin{itemize}
\item[(i)]$X_{\tau_1}<\frac{b_1}{2}n^{3/4}$,
\item[(ii)]$X_{\tau_1}> (b_2+\frac{b_1}{2})n^{3/4}$.
\end{itemize}
In case (i), let $$T_1=\min\{t:t>\tau_1, X_t\in I\}.$$ By
monotonicity of $A_t$ and the strong Markov property on $T_1$, we
have $$\P(\A \, | \, \tau_1 \leq \delta
n^{1/4},X_{\tau_1}<\frac{b_1}{2}n^{3/4}, T_1<t)\leq A_t.$$ A similar
argument as before gives us $$\P(\A,\tau_1 \leq \delta
n^{1/4},X_{\tau_1}<\frac{b_1}{2}n^{3/4}, T_1\geq t)= O(n^{-1/12}).$$
In case (ii), let $$T_2=\min\{t:t>\tau, X_t\in I\}.$$ Similar
arguments gives $$\P(\A|\tau_1 \leq \delta
n^{1/4},X_{\tau_1}<\frac{b_1}{2}n^{3/4}, T_2<t)\leq A_t,$$ and
$$\P(\A,\tau_1 \leq \delta n^{1/4},X_{\tau_1}<\frac{b_1}{2}n^{3/4},
T_2\geq t)= O(n^{-1/12}).$$ Summing over these estimate, we obtain
\be\label{estimate2} \P(\A, \tau_1 \leq \delta n^{1/4})\leq
C\delta^2(A_t+O(n^{-1/12})).\ee Plugging (\ref{estimate1}) and
(\ref{estimate2}) into (\ref{eqn_split}), we get $$A_{t,X_0}\leq
g(B)+ (A_{t-\delta
n^{1/4}}+O(n^{-1/12}))+C\delta^2(A_t+O(n^{-1/12})).$$ Maximizing
over $X_0$ and rearranging gives
$$A_t\leq \frac{1}{1-C\delta^2}(A_{t-\delta n^{1/4}}+g(B)+O(n^{-1/12})).$$
Telescoping gives
$$A_{Kn^{1/4}}\leq \frac{1}{1-C\delta^2}^{\Big\lceil\frac{K}{\delta}\Big\rceil}\Big(C\delta^2+g(B)+O(n^{-1/12})\Big).$$ Since $\frac{1}{1-C\delta^2}^{\Big\lceil\frac{K}{\delta}\Big\rceil}$ converges as $\delta$ goes to $0$, we conclude that we can choose $\delta>0$ small enough and $B$ so large to make $A_{Kn^{1/4}}$ arbitrarily small, as required.
\qed \\

\noindent{\bf{Proof of Theorem \ref{lem_hitting}:}} Notice that
\begin{eqnarray*}
  \E \Big [X_{(t+1)\wedge \tau_a}\Big |\F_t\Big ] &=& \E\Big[ X_{t+1}\1_{\{\tau_a\geq t+1\}} +X_{\tau_a}\1_{\{\tau_a\leq t\}}\Big|\F_t\Big] \\
   &=&  \E[X_{t+1}|\F_t]\1_{\{\tau_a\geq t+1\}}+X_{\tau_a}\1_{\{\tau_a\leq t\}}.
\end{eqnarray*}
By Lemma \ref{driftestimate}, we have
\begin{equation}\label{eqn_lastcondi}
  \E \Big [X_{(t+1)\wedge \tau_a}\Big |\F_t\Big ] \leq X_t\Big(1-\frac{X_t}{6n}\Big)\1_{\{\tau_a\geq t+1\}}+X_{\tau_a}\1_{\{\tau_a\leq t\}}
   = X_{t\wedge \tau_a}-\frac{X_t^2}{6n}\1_{\{\tau_a\geq t+1\}}.
\end{equation} Taking expectations on both sides of (\ref{eqn_lastcondi}), we get \be\label{expectation}\E X_{(t+1)\wedge \tau_a}\leq\E
X_{t\wedge \tau_a}-\frac{1}{6n}\E X_t^2\1_{\{\tau_a\geq t+1\}}.\ee
Note that $$\E\Big(X_t^2\1_{\{\tau_a\geq t+1\}}\Big)\geq
a^2n^{3/2}\P(\tau_a\geq t+1),$$ and $$a^2n^{3/2}\geq\frac{\E\Big(
X_{\tau_a}^2\1_{\{\tau_a\leq t\}}\Big)}{\P(\tau_a\leq t)}.$$ Hence
we have $$\E\Big(X_t^2\1_{\{\tau_a\geq t+1\}}\Big)\geq\frac{\E\Big(
X_{\tau_a}^2\1_{\{\tau_a\leq t\}}\Big)}{\P(\tau_a\leq
t)}\P(\tau_a\geq t+1),$$ which implies $$\E\Big(
X_{\tau_a}^2\1_{\{\tau_a\leq t\}}\Big)\leq \frac{\P(\tau_a\leq
t)}{\P(\tau_a\geq t+1)}\E\Big(X_t^2\1_{\{\tau_a\geq t+1\}}\Big).$$
Adding $\E\Big(X_t^2\1_{\{\tau_a\geq t+1\}}\Big)$ to both sides, we
obtain
\begin{equation}\label{eqn_fracratio}
\frac{\E\Big[ X_t^2\1_{\{\tau_a\geq t+1\}}\Big]}{\P(\tau_a\geq
t+1)}\geq\E X_{t\wedge \tau_a}^2\geq (\E X_{t\wedge \tau_a})^2.
\end{equation}
Plugging into (\ref{expectation}), we get \be\label{recursive1}\E
X_{(t+1)\wedge\tau_a}\leq\E
X_{t\wedge\tau_a}-\frac{1}{6n}\P(\tau_a\geq t+1)(\E
X_{t\wedge\tau_a})^2.\ee Note that $\E X_{(t+1)\wedge \tau_a}>0$.
Taking the inverse of (\ref{recursive1}) leads to
$$\frac{1}{\E X_{(t+1)\wedge \tau_a}}\geq\frac{1}{\E X_{t\wedge \tau_a}}+\frac{1}{6n}\P(\tau_a\geq t+1).$$
Summing $t$ from $0$ to $\lceil b n^{1/4}\rceil-1$, we get
\be\label{eqn_1}\frac{1}{\E X_{\lceil b n^{1/4}\rceil\wedge
\tau_a}}\geq\frac{1}{6n}\sum_{t=0}^{\lceil b n^{1/4}\rceil
-1}\P(\tau_a\geq t+1)\geq \frac{1}{6n}\P(\tau_a\geq b n^{1/4}) b
n^{1/4}.\ee On the other hand, for any $x\in[0,n]$, observe that
$X_{bn^{1/4}\wedge\tau_a}\leq -x$ implies there exists $t\leq
bn^{1/4}$ such that $X_t>an^{3/4}$ and $X_{t+1}<-x$. This implies
either $$\max\{|\C^+_1|,|\C_1^-|\}\leq \frac{a}{2}n^{3/4},$$ or
$$\eps\min\{|\C_1^+(t)
|,|\C_1^-(t)|\}+\sum_{j\geq 2}\epsilon_j |\C_j^+(t) | + \sum_{j\geq
2}\epsilon_j' |\C_j^-(t)|\leq -x-\frac{a}{2}n^{3/4}.$$ By Theorem
\ref{gnpout2}, we have $\P(\max\{|\C^+_1|,|\C_1^-|\}\leq
\frac{a}{2}n^{3/4})= O(e^{-cn^{1/8}}).$ By Theorem \ref{gnp2moments}
and Markov's inequality, we have
\begin{eqnarray*}\P(\eps\min\{|\C_1^+(t)
|,|\C_1^-(t)|\}+\sum_{j\geq 2}\epsilon_j |\C_j^+(t) | + \sum_{j\geq
2}\epsilon_j' |\C_j^-(t)|\leq -x-\frac{a}{2}n^{3/4}))\leq
\frac{Cn^{8/3}}{(x+\frac{a}{2}n^{3/4})^4}.\end{eqnarray*} Hence by
union bound we obtain $$\P(X_{bn^{1/4}\wedge\tau_a}\leq -x)\leq
\frac{Cn^{8/3}}{(x+\frac{a}{2}n^{3/4})^4}bn^{1/4}.$$ By a direct
computation we obtain
$$\E(|X_{bn^{1/4}\wedge\tau_a}|\1_{\{X_{bn^{1/4}\wedge\tau_a}\leq0\}})\leq \sum_{x=0}^nbn^{1/4}\frac{Cn^{8/3}}{(x+\frac{a}{2}n^{3/4})^4}=O(n^{2/3}).$$ Thus we get $$\E X_{\lceil b n^{1/4}\rceil\wedge \tau_a}\geq\P(\tau_a
>\lceil b n^{1/4}\rceil) a n^{3/4}-O(n^{2/3}).$$ Multiplying this and
(\ref{eqn_1}) we get
$$1\geq\frac{1}{6n}b n^{1/4}a n^{3/4}\Big[\P(\tau_a >\lceil b n^{1/4}\rceil)\Big]^2,$$
which gives (\ref{eqn_crossmarkov}).\qed


\subsubsection{Coupling inside the scaling window: Proof of Theorem \ref{localclt}}

\begin{lemma} \label{onelocalclt} For any fixed constant $A>1$ there
exist positive constants $q=q(A), \beta=\beta(A)$, such that if
$X_0\in [A^{-1} n^{3/4},A n^{3/4}]$, then
$$ \prob ( X_1 = x|X_0 ) \geq q n^{-5/8}$$ for any $x\in n+2\mathbb{Z}$
and $|x-X_0|\leq \beta n^{5/8}$.
\end{lemma}

\noindent{\bf{Proof of Theorem \ref{localclt}}} We will use
induction to prove that for any $\ell>0$ and any $x\in
n+2\mathbb{Z}$ such that $|X_0-x|\leq \beta(1/2+\ell /2)n^{5/8}$, we
have
\begin{equation}\label{eqn_inductclt}
\prob ( X_{\ell} = x|X_0 ) \geq
q^{\ell}\Big(\frac{\beta}{2}\Big)^{\ell -1} n^{-5/8}.
\end{equation} This implies Theorem \ref{localclt} immediately.

We prove this assertion by induction on $\ell$. Lemma
\ref{onelocalclt} implies (\ref{eqn_inductclt}) is true for $\ell =
1$. Suppose now (\ref{eqn_inductclt}) holds for $\ell$ and we prove for $\ell+1$. If $x\in
n+2\mathbb{Z}$ and $|x-X_0|\leq \beta\Big(1/2+(\ell+1) /2\Big)n^{5/8}$,
then the number of $y$ such that $y\in n+2\mathbb{Z}$ and $|y-x|\leq
\beta n^{5/8}$ and $|y-X_0|\leq \beta(1/2+\ell/2)n^{5/8}$ is at least
$\frac{\beta}{2}n^{5/8}$. Thus, we have
\begin{eqnarray}
\nonumber \P\Big(|X_\ell-x|\leq \beta n^{5/8}\Big)&=&\sum_{|y-x|\leq
\beta n^{5/8}}\P(X_{\ell}=y)\\
&\geq& q^\ell\Big(\frac{\beta}{2}\Big)^{\ell-1}
n^{-5/8}\frac{\beta}{2}n^{5/8}=
q^\ell\Big(\frac{\beta}{2}\Big)^\ell,\label{eqn_xll}
\end{eqnarray} where we used the induction hypothesis.
Since $|x-X_0|\leq \beta\Big(1/2+(\ell+1) /2\Big)n^{5/8}$, we get
$$\P\Big(X_{\ell+1}=x\,\Big|\,|X_\ell-x|\leq
\beta n^{5/8}\Big)\geq q n^{-5/8}$$ by Lemma \ref{onelocalclt}.
Together with (\ref{eqn_xll}) we get (\ref{eqn_inductclt}) for
$\ell+1$, concluding the proof. \qed\\

Recall that conditioned on the cluster sizes, $X_1$ is a summation
of independent but not identically distributed random variables. The
following is a local central limit theorem for such sums, tailored to
our particular needs, and is used to prove Lemma \ref{onelocalclt}.
We have not found in the literature a statement general enough to be
valid in our setting. The proof is the standard proof of the local
CLT using characteristic function.

\begin{lemma}\label{lem_localclt}
Suppose $K_n$ are positive integers such that $K_n\geq qn$ for some
constant $q>0$ and $a_1, a_2, \cdots, a_{K_n}$ are positive integers
such that $a_j=1$ for $1\leq j\leq qn$ and $a_j\leq
\sqrt{\frac{qn}{2}}$ for all $j$. Let $b(n)=\sum_{j=1}^{K_n}a_j$ and
$c(n)= \sqrt{\sum_{j=1}^{K_n} a_j^2 /n^{5/4}}$. Assume that there
are two positive constants $\delta$ and $C$ such that $\delta<c(n)<C$
for all $n$. Let $X_n=\sum^{K_n}_{j=1}\eps_j a_j$ where $\{\eps_j\}$
is independent random $\pm$ signs. Then for any $x\in b(n)+2\Z$ and large
enough $n$, we have
\begin{equation}\label{eqn_localclt}
\P\Big(X_n=x\Big)\geq \frac{\sqrt{2}}{\sqrt{\pi}c(n)
n^{5/8}}\Big(e^{-\frac{x^2}{2}}-1/2\sqrt{2}\Big).
\end{equation}
\end{lemma}

\noindent{\bf{Proof of Lemma \ref{onelocalclt}:}} We need to show
that with probability $\Omega(1)$ the percolation configuration fits
the setting of Lemma \ref{lem_localclt}. Define $\A_1$ and $\A_2$ as
the following events:
$$\A_1 = \big\{|\C_1^+|\in \big[X_0-\frac{c}{4}
n^{5/8},X_0+\frac{c}{4} n^{5/8}\big],\ \sum_{j\geq 2}|\C_j^+|^2\leq
D n^{5/4},\ \sum_{j\geq 2}|\C_j^+|^3\leq D n^{7/4}\big\},$$
$$\A_2 = \big\{\sum_{j\geq 1}|\C_j^-|^2\leq D n^{5/4},\ \sum_{j\geq
1}|\C_j^-|^3\leq D n^{7/4},\sum_{|\C^-_j|\leq
\frac{\sqrt{n}}{6}}|\C^-_j|^2\geq c^2 n^{5/4},\
\big|\{j:|\C^-_j|=1\}\big|\geq \frac{n}{18}\big\}$$ where $D$ and
$c$ are constants to be selected later. First we prove that $\A_1$
and $\A_2$ both happen with probability $\Omega(1)$. To bound
from below the probability of $\A_2$, take
$\delta=\frac{1}{3\sqrt{2}}$ in Theorem \ref{gnpsmalls}. We get
\begin{equation}\label{eqn_boundc1}
\P\Big(\sum_{|\C^-_j|\leq \frac{1}{6}\sqrt{n}}|\C^-_j|^2\geq
cn^{5/4}\Big)\geq q=q(A)>0,
\end{equation} for some $c=c(A)>0$. By Theorem \ref{subgnpmoments}, for $k=2,3$ we have
$$\E\sum_{j\geq 1}|\C^-_j|^k\leq C n(A^{-1} n^{-1/4})^{-2k+3} \, .$$
Thus, for
\begin{equation}\label{eqn_D2}
D\geq \frac{4CA^3}{q} \, ,
\end{equation}
we have by Markov's inequality that
\begin{equation}\label{eqn_F22}
\P\Big(\sum_{j\geq 1}|\C^-_j|^2\geq D n^{5/4}\Big)\leq \frac{q}{4}
\end{equation} and
\begin{equation}\label{eqn_F23}
\P\Big(\sum_{j\geq 1}|\C^-_j|^3\geq D n^{7/4}\Big)\leq \frac{q}{4}.
\end{equation} By Lemma \ref{lem_isolated}, we have
\begin{equation}\label{eqn_F24}
\P\Big(\big|\{j: |\C^-_j|=1\}\big|\geq \frac{n}{18}\Big)\geq
1-C/n\geq 1-\frac{q}{4}.
\end{equation} Putting (\ref{eqn_boundc1}), (\ref{eqn_F22}),
(\ref{eqn_F23}) and (\ref{eqn_F24}) together, we get
$$\P(\A_2)\geq \frac{q}{4}.$$

To bound from below the probability of $\A_1$, we apply Theorem
\ref{giantclt} to get that
$$\P\Big(\Big| |\C_1^+|-x_0(1+\frac{x_0}{n})\Big|\leq \frac{c}{4}(\frac{n+x_0}{2})^{5/8}
\Big)\geq q=q(A)>0.$$ Since $\frac{x_0^2}{n}=o(n^{5/8})$ and
$\frac{n+x_0}{2}\leq \frac{3}{4}n$, we get
\begin{equation}\label{eqn_F11}
\P\Big(|\C_1^+|\in \Big[x_0-\frac{c}{4}n^{5/8},\,
x_0+\frac{c}{4}n^{5/8}\Big]\Big)\geq q.
\end{equation}
By Theorem \ref{gnp2moments}, for $k=2,3$ we have
$$\E\sum_{j\geq 2}|\C^+_j|^k\leq C_kn(A^{-1} n^{-1/4})^{-2k+3}.$$ Again, when $D$ satisfies \ref{eqn_D2} we get
by Markov's inequality that
\begin{equation}\label{eqn_F12}
\P\Big(\sum_{j\geq 2}|\C^+_j|^2\geq D n^{5/4}\Big)\leq \frac{q}{4}
\end{equation} and
\begin{equation}\label{eqn_F13}
\P\Big(\sum_{j\geq 2}|\C^+_j|^3\geq D n^{7/4}\Big)\leq \frac{q}{4}.
\end{equation}
By (\ref{eqn_F11}), (\ref{eqn_F12}) and (\ref{eqn_F13}), we have
\be\label{couplepart1}\P(\A_1)\geq \frac{q}{2}.\ee

Since $\A_1$ and $\A_2$ are independent, we get
$$\P(\A_1\cap\A_2)\geq \frac{q^2}{8},$$ providing $D$ satisfies (\ref{eqn_D2}). By Proposition \ref{c1plus} we have $\P(|C_1^-|\geq|C_1^+|)= O(e^{-c\log^2 n})$. Hence the event $$\A=\{\A_1,\A_2, |C_1^-|<|C_1^+|)\}\, ,$$ occurs with probability $\Omega(1)$.

Next we prove that for every $x\in n+2\mathbb{Z}$ and $|x-X_0|\leq
\frac{c}{2}n^{5/8}$, there exist a constant $\delta>0$ such that
\begin{equation}\label{eqn_condclt}
\P(X_1 =x\,|\A)\geq \delta n^{-5/8} \, ,
\end{equation} which will conclude the proof. Denote
$$M_1=|\C^+_1|+\sum_{j\geq 2}\eps_j |\C^+_j |
+\sum_{|\C_j|>\sqrt{n}/6} \eps'_j |\C^-_j |$$ and
$$M_2=\sum_{|\C_j|\leq\sqrt{n}/6} \eps'_j |\C^-_j |.$$ Note that $M_1$ and $M_2$ are independent conditioned on $\A$. We will first prove that there exist a constant $\alpha>0$ such that \begin{equation}\label{eqn_bothsmall}
\P\Big(|M_1-X_0|\leq \frac{c}{2} n^{5/8}\big|\A\Big)\geq \alpha.
\end{equation} On $\A$ we have that \be\label{eqn_secondmo}\sum_{j\geq 2} |\C^+_j |^2
+\sum_{|\C_j|>\sqrt{n}/6} |\C^-_j |^2\leq 2D n^{5/4}\ee and
\be\label{eqn_thirdmo}\sum_{j\geq 2} |\C^+_j |^3
+\sum_{|\C_j|>\sqrt{n}/6} |\C^-_j |^3\leq 2D n^{7/4}.\ee If
$\sum_{j\geq 2} |\C^+_j |^2 +\sum_{|\C_j|>\sqrt{n}/6} |\C^-_j
|^2\leq \frac{c^2}{32} n^{5/4}$, by Markov's inequality, we have
\be\label{markov11}\P\Big(\Big|\sum_{j\geq 2}\eps_j |\C^+_j |
+\sum_{|\C_j|>\sqrt{n}/6} \eps'_j |\C^-_j |\Big|\leq \frac{c}{4}
n^{5/8}\Big)\geq 1/2.\ee Otherwise
$$\frac{|\C_j|}{\big(\sum_{j\geq 2} |\C^+_j |^2 +\sum_{|\C_j|>\sqrt{n}/6}
|\C^-_j |^2\big)^{1/2}}>\eps$$ implies $$|\C_j|\geq {\eps c^2 \over 32} n^{5/8} \, ,$$
and since $\{|\C_j|\}$ also satisfy (\ref{eqn_thirdmo}), we learn that the Lindeberg
condition is satisfied. By Lindeberg-Feller theorem (see \cite{D},
(4.5)), we have
\begin{equation}\label{eqn_bigsmall}
\P\Big(\Big|\sum_{j\geq 2}\eps_j |\C^+_j | +\sum_{|\C_j|>\sqrt{n}/6}
\eps'_j |\C^-_j |\Big|\leq \frac{c}{4} n^{5/8}\Big)\geq \alpha>0.
\end{equation} Combining this and (\ref{markov11}) yields (\ref{eqn_bothsmall}).

To estimate $M_2$ let
$$b=\sum_{|\C_j|\leq\sqrt{n}/6}|\C^-_j| \quad \and \quad a= n^{-9/8} b^{1/2}  \, .$$
By Lemma \ref{lem_localclt}, for every $x\in b +2\mathbb{Z}$, we
have
$$\P\big(M_2=x\big|\A\big)\geq \frac{\sqrt{2}}{\sqrt{\pi}an^{5/8}}\Big(e^{-\frac{x^2}{2a^2 n^{5/4}}}-1/2\Big).$$
For all $x$ such that $|x|\leq c n^{5/8}$, we have
$$\frac{\sqrt{2}}{\sqrt{\pi}an^{5/8}}\Big(e^{-\frac{x^2}{2a^2
n^{5/4}}}-1/2\Big) \geq \frac{\sqrt{2}}{\sqrt{\pi}D
n^{5/8}}\Big(e^{-1/2}-1/2\Big)\geq \delta n^{-5/8}$$ where $\delta$
is a constant. So for every $x\in b+2\mathbb{Z}$ and $|x|\leq c
n^{5/8}$, we have
\begin{equation}\label{eqn_restclt}
\P\big(M_2=x\big|\A\big)\geq \delta n^{-5/8}.
\end{equation}
By (\ref{eqn_bothsmall}) and (\ref{eqn_restclt}), for every $x\in
n+2\mathbb{Z}$ with $|x-x_0|\leq \frac{c_1}{2} n^{5/8}$, we have
\begin{eqnarray*}
&&\P\big(M_1+M_2=x\big|\A\big)\\&\geq& \P\Big(|M_1-x_0|\leq
\frac{c_1}{2}n^{5/8}\, ,\
M_2 =(x-x_0)-(M_1-x_0)\big|\A\Big)\\
&\geq& \alpha\delta n^{-5/8}.
\end{eqnarray*}
This proves (\ref{eqn_condclt}), which concludes the whole
proof.\qed\\

To prove Lemma \ref{lem_localclt} we need the following two small
assertions. The first is Exercise 3.2 of \cite{D}.
\begin{lemma}\label{lem_durretclt}
If $\P(X\in b+h\Z)=1$, where $b$ is a complex number and $h>0$ is a
real number. Then for any $x\in b+h\Z$, we have
$$\P(X=x)=\frac{h}{2\pi}\int_{-\pi/h}^{\pi/h} e^{-itx}\phi(t)\d t,$$
where $\phi(t)$ is the characteristic function of $X$.
\end{lemma}

\begin{lemma}\label{cor_fancycosx}
For any $x$ in $\mathbb{R}$, let $m(x)$ be the integer that is
closest to $x$ (if $x-\frac{1}{2}$ is an integer, then we put
$m(x)=x-\frac{1}{2}$). Then for any $x$
$$|\cos x|\leq \exp\Big(-\frac{(x-m(\frac{x}{\pi})\pi)^2}{2}\Big).$$
\end{lemma}
\noindent{\bf{Proof.}} Since $m(\frac{x}{\pi})\pi\in
\{k\pi\}_{k\in\Z}$, we have $|\cos x|=|\cos
(x-m(\frac{x}{\pi})\pi)|$. Also, we have $-\frac{\pi}{2}\leq
x-m(\frac{x}{\pi})\pi\leq\frac{\pi}{2}$. Since $\cos x\leq
e^{-\frac{x^2}{2}}$ for all $x\in[-\frac{\pi}{2},\frac{\pi}{2}]$ we
have that
$$|\cos x|=\cos
\Big(x-m(\frac{x}{\pi})\pi\Big)\leq
\exp\Big(-\frac{(x-m(\frac{x}{\pi})\pi)^2}{2}\Big)\, .$$ \qed

\noindent{\bf{Proof of Lemma \ref{lem_localclt}:}} For simplicity we
will abbreviate $c(n)$ by $c$. Let
$$d_j=\frac{a_j}{c n^{5/8}}.$$ Then we have $\sum_{j=1}^{K_n}
d_j^2 =1$ and $\frac{X_n}{c n^{5/8}}= \sum_{j=1}^{K_n}\eps_j d_j$.
Since $a_j=O(n^{1/2})$, we have that $d_j=O(n^{-1/8})$. Thus it
satisfies Lindeberg condition (see \cite{D}). Consequently, we have
that
$$\frac{X_n}{c n^{5/8}}\stackrel{d}{\rightarrow} N(0,1).$$ Denote the
characteristic function of $\frac{X_n}{c n^{5/8}}$ by $\phi_n(t)$. A
straightforward computation gives that
\begin{equation}\label{eqn_charact} \phi_n(t)=\Big(\cos\frac{t}{c n^{5/8}}\Big)^{qn}\prod_{j=qn+1}^{K_n}\cos
(t d_j) \, ,
\end{equation} and we have $\phi_n(t)\rightarrow e^{-\frac{t^2}{2}}$
for all fixed $t\in\mathbb{R}$. Taking $h=\frac{2}{cn^{5/8}}$ in
Lemma \ref{lem_durretclt}, for $x\in b(n)+2\Z$ we have
\begin{equation}\label{eqn_durretuse}
\P\Big(X_n=x\Big)=\frac{1}{\pi c n^{5/8}}\int^{\frac{\pi}{2}c
n^{5/8}}_{-\frac{\pi}{2}c n^{5/8}} e^{-itx}\phi_n(t)\d t.
\end{equation} Let $M$ be a large constant to be selected later. Note that
$\phi_n(t)$ is an even function so
\begin{eqnarray}
\int^{\frac{\pi}{2}c n^{5/8}}_{-\frac{\pi}{2}c n^{5/8}}
e^{-itx}\phi_n(t)\d t&=&\int^{M}_{-M} e^{-itx}\phi_n(t)\d t+
2\int^{\frac{\pi}{2}c n^{5/8}}_{M} e^{-itx}\phi_n(t)\d t\nonumber
\\
&\geq&\int^{M}_{-M} e^{-itx}\phi_n(t)\d t-2\int^{\frac{\pi}{2}c
n^{5/8}}_{M} |\phi_n(t)|\d t.\label{eqn_cltsplit}
\end{eqnarray} We will first bound from above the second term of (\ref{eqn_cltsplit}). Let $m_j(t)=m\Big(\frac{td_j}{\pi}\Big) \frac{\pi}{d_j}$,
i.e., $m_j(t)$ is the element in
$\Big\{k\frac{\pi}{d_j}\Big\}_{k\in\Z}$ that is closest to $t$. Note
that by Lemma \ref{cor_fancycosx}, we have
\begin{equation*}
\cos (t d_j)\leq \exp\Big\{-\Big[t d_j-m\big(\frac{t
d_j}{\pi}\big)\pi\Big]^2
\Big/2\Big\}=\exp\Big\{-d_j^2\frac{(t-m_j(t))^2}{2}\Big\}.
\end{equation*}
For large enough $n$, we have $\frac{1}{c^2 n^{5/4}}\geq
\frac{1}{2}\frac{1}{c^2 n^{5/4}-qn}$. Thus, we get
\begin{equation*}
\Big|\cos (t d_j)\Big|\leq \exp\Big\{-\frac{a_j^2}{c^2
n^{5/4}-qn}\cdot \frac{(t-m_j(t))^2}{4}\Big\}.
\end{equation*}
Since $\sum_{j=qn+1}^{K_n}\frac{(a_j)^2}{c^2 n^{5/4}-qn}=1$ and
$e^{-x}$ is a convex function, we have by Jensen's inequality that
\begin{eqnarray}
\prod_{j=qn+1}^{K_n}\Big|\cos (td_j)\Big|&\leq&
\exp\Big\{-\sum_{j=qn+1}^{K_n}\frac{a_j^2}{c^2
n^{5/4}-qn}\frac{(t-m_j(t))^2}{4}\Big\}\nonumber\\
&\leq& \sum_{j=qn+1}^{K_n}\frac{a_j^2}{c^2
n^{5/4}-qn}\exp\Big(-\frac{(t-m_j(t))^2}{4}\Big) .\label{eqn_bigcos}
\end{eqnarray}
Recall that $|t|\leq \frac{\pi}{2}c n^{5/8}$ and $|\cos(x)|\leq
e^{-\frac{x^2}{2}}$ for $x\in[-\frac{\pi}{2},\frac{\pi}{2}]$, whence
\begin{equation}\label{eqn_smallcos}
\Big|\cos\frac{t}{c n^{5/8}}\Big|^{qn}\leq\exp \Big(-\frac{q
t^2}{2c^2 n^{1/4}}\Big).
\end{equation}
Plugging (\ref{eqn_bigcos}) and (\ref{eqn_smallcos}) into
(\ref{eqn_charact}), we get
\begin{equation*}
|\phi_n(t)|\leq \sum_{j=qn+1}^{K_n}\frac{a_j^2}{c^2
n^{5/4}-qn}\exp\Big(-\frac{(t-m_j(t))^2}{4}-\frac{q t^2}{2c^2
n^{1/4}}\Big).
\end{equation*}
Hence, we have
\begin{equation}\label{eqn_intsupM}
\int^{\frac{\pi}{2}c n^{5/8}}_{M} |\phi_n(t)|\d t\leq
\sum_{j=qn+1}^{K_n}\frac{a_j^2}{c^2 n^{5/4}-qn}\int^{\infty}_{M}
\exp\big(-\frac{(t-m_j(t))^2}{4}-\frac{q t^2}{2c^2 n^{1/4}}\Big)\d
t.
\end{equation}

We will divide the integral into two parts such that the first part
converges to $0$ as $M$ goes to infinity and the second part is
bounded by a constant. Recall that $m_j(t)=0$ for
$t\in[-\frac{\pi}{2d_j},\frac{\pi}{2d_j}]$, so for any
$j\in[qn+1,K_n]$, we have
\begin{eqnarray} & &\int^{\infty}_{M}
\exp\Big(-\frac{(t-m_j(t))^2}{4}-\frac{q t^2}{2c^2 n^{1/4}}\Big)\d
t\nonumber\\ &=& \int^{\frac{\pi}{2d_j}}_{M}
\exp\Big(-\frac{t^2}{4}-\frac{q t^2}{2c^2 n^{1/4}}\Big)\d
t\nonumber\\&+&
\sum_{\ell=1}^{\infty}\int^{\frac{\pi}{2d_j}(2\ell+1)}_{\frac{\pi
}{2d_j}(2\ell-1)} \exp\Big(-\frac{(t-m_j(t))^2}{4}-\frac{q t^2}{2c^2
n^{1/4}}\Big)\d t\label{complex}\
\end{eqnarray} The first term of the right hand side of (\ref{complex}) is bounded
by $\int^{\infty}_{M} e^{-\frac{t^2}{4}}\d t$. For the second term,
note that for $t\geq \frac{\pi}{2d_j}(2\ell -1)$, we have
$$\exp\Big(-\frac{q
t^2}{2c^2 n^{1/4}}\Big)\leq \exp\Big(-\frac{q\pi^2 n(2\ell-1)^2}{8
a_j^2}\Big)$$ and
\begin{eqnarray}
\int^{y+\frac{\pi}{d_j}}_y e^{-\frac{1}{4}(t-m_j(t))^2}\d t &=&
\int^{\frac{\pi}{2d_j}}_{-\frac{\pi}{2d_j}}
e^{-\frac{1}{4}(t-m_j(t))^2}\d t\nonumber\\
&=& \int^{\frac{\pi}{2d_j}}_{-\frac{\pi}{2d_j}} e^{-\frac{t^2}{4}}\d
t \leq 2\sqrt{\pi}.\label{eqn_int2pi}
\end{eqnarray} for any $y$, since $\frac{1}{4}(t-m_j(t))^2$ is a periodic
function. Thus, we get
\begin{equation}\label{eqn_biggerM}
\int^{\infty}_{M} \exp\Big(-\frac{(t-m_j(t))^2}{4}-\frac{q t^2}{2c^2
n^{1/4}}\Big)\d t \leq \int^{\infty}_{M} e^{-\frac{t^2}{4}}\d t+
\sum_{\ell=1}^{\infty}2\sqrt{\pi}\exp\Big(-\frac{q\pi^2
n(2\ell-1)^2}{8 a_j^2}\Big).
\end{equation}
Recall that $a_j\leq \sqrt{qn/2}$, hence
\begin{eqnarray*}
\sum_{\ell=1}^{\infty}\exp\Big(-\frac{q\pi^2 n(2\ell-1)^2}{8
a_j^2}\Big)&\leq& \sum_{\ell=1}^{\infty}\exp\Big(-\frac{\pi^2
(2\ell-1)^2}{4}\Big)\\
&\leq& \frac{e^{-\pi^2 / 4}}{1-e^{-\pi^2/2}}\leq \frac{1}{8}.
\end{eqnarray*}
Plug into (\ref{eqn_biggerM}), we get
$$\int^{\infty}_{M} \exp\Big(-\frac{(t-m_j(t))^2}{4}-\frac{q t^2}{2c^2
n^{1/4}}\Big)\d t \leq
2\sqrt{\pi}\Big(1-\Phi\big(\frac{M}{\sqrt{2}}\big)\Big)+
\frac{\sqrt{\pi}}{4},$$ where $\Phi(\cdot)$ is the distribution
function of $N(0,1)$. Plugging back into (\ref{eqn_intsupM}), we get
\begin{equation}\label{eqn_bigMpart}
\int^{\frac{\pi}{2}c n^{5/8}}_{M} |\phi_n(t)|\d t\leq
2\sqrt{\pi}\Big(1-\Phi\big(\frac{M}{\sqrt{2}}\big)\Big)+
\frac{\sqrt{\pi}}{4}.
\end{equation}\\

Now we go back to the first term of the right hand side of
(\ref{eqn_cltsplit}). Recall that $\phi_n(t)$ converge to
$e^{-t^2/2}$ for all $t$. We have the following estimate:
\begin{eqnarray}
& &\int^{M}_{-M} e^{-itx}\phi_n(t)\d t\nonumber\\ &=&
\int^{\infty}_{-\infty} e^{-itx}e^{-t^2/2}\d t-
\Big(\int^{-M}_{-\infty}+\int^{\infty}_{M}\Big) e^{-itx}e^{-t^2/2}\d
t +\int^{M}_{-M}
e^{-itx}(\phi_n(t)-e^{-t^2/2})\d t\nonumber\\
&\geq& \sqrt{2\pi}e^{-x^2/2} -2\int^{\infty}_{M} e^{-t^2/2}\d t
-\int^{M}_{-M} \big|\phi_n(t)-e^{-t^2/2}\big|\d t.\label{eqn_center}
\end{eqnarray}
Note that the second term of the left most side of
(\ref{eqn_center}) converges to 0 as $M\rightarrow \infty$ and for
fixed $M$ we have $\int^{M}_{-M} |\phi_n(t)-e^{-\frac{t^2}{2}}|\d
t\rightarrow 0$ by the Dominated Convergence Theorem. Plugging these
and (\ref{eqn_bigMpart}) into (\ref{eqn_cltsplit}), we get
$$\liminf_{n\rightarrow\infty}\int^{\frac{\pi}{2}c
n^{5/8}}_{-\frac{\pi}{2}c n^{5/8}} e^{-itx}\phi_n(t)\d
t\geq\sqrt{2\pi}e^{-x^2/2}-\frac{\sqrt{\pi}}{2},$$ which concludes
the whole proof. \qed


\subsection{Starting at the $[0,n^{3/4}]$ regime: Proof of Theorem \ref{startzero}}
\begin{thm}\label{pushup2}
Let $I=[-A n^{2/3}, An^{2/3}]$ where $A$ is a fixed large
constant.Then there exist positive constants $K,a,q$ such that
\begin{equation}\label{eqn_willupqq}
\Pr\Big(\tau_a \leq Kn^{1/4}\, |\, X_0\in I\Big)\geq q
\end{equation} where $\tau_a=\inf \{t\geq 0: X_t\geq an^{3/4}\}$.
\end{thm}

\begin{thm}\label{pushdownuse}
For a constant $A$ put $I=[-An^{2/3},\, An^{2/3}]$ and $\tau =\inf
\{t\geq 0:X_t\in I\}$. Then there exist constant $c>0$ such that for
sufficiently large $A$, we have
$$\P(\tau
>t) \leq \frac{ 2|X_0|}{ct\sqrt{n}}.$$
\end{thm}

\noindent{\bf{Proof of Theorem \ref{startzero}:}} Let $A,c$ be
constants such that the assertion of Theorem \ref{pushdownuse}
holds and write $I=[-A n^{2/3}, An^{2/3}]$. Since $X_0\leq
n^{3/4}$, by Lemma \ref{pushdownuse} with $t=\frac{4n^{1/4}}{c}$, we
have
$$\Pr(\tau \leq\frac{4n^{1/4}}{c})
> \frac{1}{2}$$ where $\tau=\min\{t:X_t\in I\}$.  By Theorem
\ref{pushup2} and the strong Markov property, we get that $X_t$
exceeds $a n^{3/4}$ within $(K+\frac{4}{c})n^{1/4}$ steps with
probability at least $\frac{q}{2}$. By Theorem \ref{allplus}, we can
couple $X_t$ and with the stationary chain $Y_t$ within $O(n^{1/4})$
steps such that they meet each other with probability $\Omega(1)$. Applying Lemma \ref{thm_coupling} concludes the proof. \qed\\

We now proceed to the proof of Theorem \ref{pushup2}. We begin with
some lemmas.

\begin{lemma}\label{pushup1}
Let $A$ be a large constant and put $I=[-A n^{2/3}, An^{2/3}]$. For
any $q\in(0,1)$, there exist a state $Z=Z(q)\in I$ and constants
$a=a(q)>0$, $K=K(q)>0$ such that
\begin{equation}\label{eqn_willup}
\Pr\Big(\tau_a \leq Kn^{1/4}|X_0=Z\Big)\geq q
\end{equation}  where $\tau_a =\inf \{t\geq 0:X_t\geq an^{3/4}\}$.
\end{lemma}
\noindent{\bf{Proof.}} Let $Y_t$ be a SW chain with
$Y_0\stackrel{d}{=}\pi_n$. By Theorem \ref{thm_statregime}, there
exists a constant $B$ such that
\begin{equation}\label{eqn_station}
\pi_n([B^{-1}n^{3/4},\,Bn^{3/4}])\geq 1-\frac{1-q}{4}.
\end{equation}
We will prove the lemma for $a=B^{-1}$ and $K=\frac{6B}{c}$ where
$c$ is the constant in Theorem \ref{pushdownuse}. Write
$J=[B^{-1}n^{3/4},\,Bn^{3/4}]$, then by (\ref{eqn_station}) we have
\begin{equation}\label{eqn_bothInBig}
\Pr\Big(Y_0\in J \and Y_{Kn^{1/4}}\in J\Big)\geq 1-\frac{1-q}{2}.
\end{equation} Put $\tau =\inf \{t\geq 0:Y_t\in I\}$. We have
$$\Pr(\tau \leq Kn^{1/4}|Y_0\in J) \geq 1- \frac{2 Bn^{3/4}}{c\sqrt{n} K n^{1/4}}
=\frac{2}{3}$$ by Theorem \ref{pushdownuse}. Thus we have
$$\Pr\Big(Y_0\in J,\, \tau \leq Kn^{1/4}\Big)\geq \Big(1-\frac{1-q}{4}\Big)\frac{2}{3}\geq \frac{1}{2}.$$ Let $$\delta=\max_{W\in I}\P(\tau_a\leq Kn^{1/4}|X_0=W).$$ By the strong Markov property
\begin{eqnarray*}
&&\P\Big(Y_0\in J,\, \tau \leq Kn^{1/4},\, Y_{Kn^{1/4}} \leq
an^{3/4}\Big)\nonumber\\&=&\P\Big(Y_0\in J,\, \tau \leq
Kn^{1/4}\big)\P\big(Y_{Kn^{1/4}}\leq an^{3/4}|\tau,
Y_{\tau}\big)\geq\frac{1-\delta}{2} \, ,
\end{eqnarray*}
since $\tau_a > Kn^{1/4}$ implies that $Y_{Kn^{1/4}} \leq
an^{3/4}$. We deduce that \be\label{pushdown1}\P\Big(Y_0\in
J,\,Y_{Kn^{1/4}}\not\in J \Big)\geq\frac{1-\delta}{2}.\ee Combining
(\ref{eqn_bothInBig}) and (\ref{pushdown1}) we get that $\delta \geq q$, concluding our proof. \qed\\

\begin{lemma}\label{criticalgnp}
Consider the random graph $G(\frac{n+X_0}{2},\frac{2}{n})$ where $X_0\in[-An^{2/3},A
n^{2/3}]$ for some large constant $A$. Then the intersection of the following events occurs with probability at least
$\delta=\delta(A)>0$:
\begin{itemize}
\item $|\C_1|+|\C_2|\in[4An^{2/3},8An^{2/3}],$ and
$|\C_2|>\frac{4A}{3}n^{2/3}$,
\item $\C_1$ and $\C_2$ are trees, and
\item $\sum_{j \geq 3} |\C_j|^2 \leq n^{4/3}$.
\end{itemize}
\end{lemma}

\noindent{\bf{Proof.}} Let $\A$ be the event
\begin{itemize}
\item $|\C_1|+|\C_2|\in[4An^{2/3},8An^{2/3}],|\C_2|>\frac{4A}{3}n^{2/3},$
\item $\C_1$ and $\C_2$ are trees.
\end{itemize} By Theorem 5.20 of \cite{JLR}, we have $\P(\A)\geq \delta=\delta(A)>0$. Conditioned on $\A$ and on $\C_1$ and $\C_2$ the remaining graph, $\{\C_j\}_{j\geq3}$, is distributed
$G(\frac{n+X_0}{2}-|\C_1|-|\C_2|,\frac{2}{n})$ conditioned to the
event that it does not have components larger than $|\C_2|$. By
Theorem 7 of \cite{NP1} the complement of this event has probability
decaying exponentially in $A$. Let $\{\C_j'\}$ be the component size
in the unconditioned space
$G(\frac{n+X_0}{2}-|\C_1|-|\C_2|,\frac{2}{n})$. We have
$$\E\sum_{j \geq 1}
|\C_j'|^2=\big (\frac{n+X_0}{2}-|\C_1|-|\C_2|\big )\E|\C(v)| \, .$$
Since $X_0 \leq An^{2/3}$ and $|\C_1|+|\C_2|\geq 4An^{2/3}$, Theorem
7 of \cite{NP1} gives that
$$ \E |\C(v)| \leq O(e^{-cA}) n^{1/3} \, ,$$
and so
$$\E\sum_{j \geq 1} |\C_j'|^2\leq O(e^{-cA}) n^{4/3} \, .$$
The lemma now follows since in the conditioned space, the event we condition on has probability exponentially close to $1$. \qed \\

\begin{lemma}\label{uniformbd} Let $I=[-An^{2/3}, An^{2/3}]$ for some large $A$. There exist a constant $c=c(A)>0$ such
that \be\label{localcouple} \prob ( X_1 = x\,|\,X_0 \in I ) \geq c
n^{-2/3}\ee for any $x\in n+2\Z$ with $x \in I$ and $x>0$.
\end{lemma}

\noindent{\bf{Proof.}} Write $\A$ for the event
of the assertion of Lemma \ref{criticalgnp} in $\{\C_j^+\}$, so that
$\P(\A)\geq\delta(A)>0$. In $G(\frac{n-|X_0|}{2},\frac{2}{n})$ we
have by Theorem \ref{gnp2moments} that $$ \E \sum _{j \geq 1}
|\C_j^-|^2 \leq D n^{4/3}$$ where $D=D(A)$ is a constant. We have by
Markov's inequality that $$ \P\Big(\big|\sum _{j \geq 3} \epsilon_j
|\C_j^+|\big| \leq D n^{2/3}\big|\A\Big)\geq1-1/D^2 \, ,$$ and
$$\P \big(\big|\sum _{j \geq 1} \epsilon_j |\C_j^-|\big| \leq D
n^{2/3}\big)\geq1-1/D^2\, ,$$ and these two events are independent. Thus,
the following event which we denote by $\B$ happens with probability $\Omega(1)$.
\begin{itemize}
\item $|\C_1^+|+|\C_2^+| \in [4An^{2/3}, 8An^{2/3}],|\C_2^+|>\frac{4A}{3}n^{2/3}$,
\item $\C_1^+$ and $\C_2^+$ are trees,
\item $\big|\sum_{j \geq 3} \epsilon_j |\C_j^+|+\sum_{j \geq 1} \epsilon_j
|\C_j^-|\big|\leq 2Dn^{2/3}.$
\end{itemize} Note that if a negative spin is assigned to $\C_2^+$ then $$X_1=|\C_1^+|-|\C_2^+|+\sum_{j \geq 3} \epsilon_j |\C_j^+|+\sum_{j \geq 1} \epsilon_j|\C_j^-|.$$ Thus $$\P(X_1=x)\geq \frac{1}{2}\P\big (|\C_1^+|-|\C_2^+|+\sum_{j \geq 3} \epsilon_j |\C_j^+|+\sum_{j \geq 1} \epsilon_j |\C_j^-|=x \big ).$$ So we only need to show that for any
$x\in[-An^{2/3},An^{2/3}]$ we have
\be\label{criticalgnp1}\P(|\C_1^+|-|\C_2^+|=x\,|\,\B)\geq c
n^{-2/3},\ee for some constant $c=c(A)>0$. For any
$m\in[4An^{2/3},8An^{2/3}]$ let $l=\frac{m+x}{2}$. By Cayley's
formula we have that
\begin{eqnarray}&&\P(|\C_1^+|-|\C_2^+|=x\,\big|\,|\C_1^+|+|\C_2^+|=m, \C_1\cup\C_2,\B)\nonumber\\&=&\P(|\C_1^+|=l,|\C_2^+|=m-l\,\big|\,|\C_1^+|+|\C_2^+|=m, \C_1\cup\C_2,\B)\nonumber\\&=& { {m \choose l} l^{l-2}
(m-l)^{(m-l)-2} \over \sum_{k=\frac{4A}{3}n^{2/3}}^{m/2} {m \choose
k} k^{k-2} (m-k)^{(m-k)-2}}.
\end{eqnarray} Let $$a(k)= {m
\choose k} k^{k-2} (m-k)^{(m-k)-2}.$$ By Stirling's formula, there
are two constants $c$ and $C$ such that for large enough $n$ and any
$k_1, k_2\in[\frac{4A}{3}n^{2/3},\frac{m}{2}]$, we have
$$c\leq\frac{a(k_1)}{a(k_2)}\leq C.$$ This implies $${ {m
\choose k} k^{k-2} (m-k)^{(m-k)-2} \over
\sum_{k=\frac{4A}{3}n^{2/3}}^{m/2} {m \choose k} k^{k-2}
(m-k)^{(m-k)-2}}\geq cn^{-2/3}$$ which proves (\ref{criticalgnp1}).\qed\\

\noindent{\bf{Proof of Theorem \ref{pushup2}:}} Let $A$ be large and
$q\in (0,1)$ will be chosen later very close to $1$. Let $Z$ be the
site and $K>0$ the number satisfying the assertion of Lemma
\ref{pushup1}. Let $\{\widetilde{X}_t\}$ be an independent SW chain
starting at $Z$ and $\widetilde{\tau}_a$ is as in Lemma
\ref{pushup1}. Then we have
\begin{equation}\label{eqn_TauUp}
\Pr\Big(\widetilde{\tau}_a \geq Kn^{1/4}\, \mid \, \widetilde{X}_0 =
Z \Big)\leq 1-q \, . \end{equation}
Let $c>0$ be the constant from Lemma \ref{uniformbd}. This lemma
implies that we can couple $X_t$ and $\widetilde{X}_t$ such that
$X_1 = \widetilde{X}_1$ with probability at least $c$. From that
point we can couple such that the two processes stay together with
probability $1$.
\begin{equation}\label{eqn_coupleafter}
\Pr(X_t=\widetilde{X}_t \hbox{ for } t\geq 1)\geq c \, .
\end{equation} Thus, we have
\begin{eqnarray*}
&&\P(\tau_a\leq Kn^{1/4})\geq\P(X_t=\widetilde{X}_t \hbox{ for }
t\geq 1 \and \widetilde{\tau}_a\leq Kn^{1/4}) \geq c - (1-q) \, ,
\end{eqnarray*} so we choose $q \geq 1-c/2$ and conclude the proof. \qed \\

To prove Theorem \ref{pushdownuse} we consider yet another
modification of the SW dynamics $\{X_t'\}$. For any $X_0'$, in the
supercritical random graph $G(\frac{n+|X_0|'}{2},\frac{2}{n})$, let
$\C_{\delta\epsilon n}$ be the component discovered by the
exploration process at time $\delta\epsilon n$ where
$\epsilon=\frac{X_0'}{n}$ and $\delta$ is a small constant (see
Lemma \ref{critgnp}). We assign positive spin to this component and
random spins to all other components in
$G(\frac{n+|X_0|'}{2},\frac{2}{n})$ and all components in
$G(\frac{n-|X_0|'}{2},\frac{2}{n})$. Let $X_1'$ be the sum of spins
after this assigning process.

The reason we require this change is that we were not able to obtain
the bounds of Theorem \ref{critgnp} for $\C_1$, but only for
$\C_{\delta \epsilon n}$ which is very likely to be $\C_1$. This
will become evident in the proof. We first state a key lemma and
then use it to prove Theorem \ref{pushdownuse}.

\begin{lemma}\label{downdriftlem}
For any constant $A$ put $I=[-An^{2/3},\, An^{2/3}]$. Then there
exists a constant $c>0$ such that for sufficiently large $A$ we have
\begin{equation}\label{eqn_driftone} \E\Big(|X'_{1}| \1_{\{X'_1
\not\in I\}} +X'_{1} \1_{\{X'_1 \in I\}} \, \mid \, |X'_0| > An^{2/3} \Big)\leq |X'_0| -c\sqrt{n} \, .
\end{equation}
\end{lemma}

\noindent{\bf{Proof of Theorem \ref{pushdownuse}:}} Notice that of
$|X_0|=|X_0'|$ then $|X_1|\stackrel{d}=|X_1'|$, and so
$|X_t|\stackrel{d}=|X_t'|$ for all $t\geq1$. Thus, we only need to
prove the assertion of the Theorem for $\{X_t'\}$. For simplicity of
notation we write $X_t$ for $X_t'$. Assume that $|X_0| > An^{2/3}$ otherwise the assertion is trivial. We begin by noticing that \be\label{eqn_case11}\E\Big(|X_{t+1}|
\1_{\{\tau
>t+1\}} +X_{\tau} \1_{\{\tau\leq t+1\}}\Big| \F_t\Big)\1_{\{\tau\leq t\}}=\E
(X_{\tau}\1_{\{\tau\leq t\}}|\F_t)=X_{\tau}\1_{\{\tau\leq t\}} \, .\ee
By Lemma \ref{downdriftlem} we have
\begin{eqnarray}\E\Big(|X_{t+1}| \1_{\{\tau
>t+1\}} +X_{\tau} \1_{\{\tau\leq t+1\}}\Big| \F_t\Big)\1_{\{\tau\geq t+1\}}&=&\E\Big(|X_{t+1}| \1_{\{\tau
>t+1\}} +X_{\tau} \1_{\{\tau=t+1\}}\Big| \F_t\Big)\nonumber\\&\leq&|X_t|\1_{\{\tau>t\}}-c\sqrt{n}\1_{\{\tau
>t\}}.\label{eqn_case22}\end{eqnarray} Thus, We
have that $\{X_t\}$ satisfies the following inequality:
\begin{eqnarray}
\E\Big(|X_{t+1}| \1_{\{\tau >t+1\}} +X_{\tau} \1_{\{\tau\leq
t+1\}}\Big| \F_t\Big)\leq \ |X_{t}| \1_{\{\tau >t\}} +X_{\tau}
\1_{\{\tau\leq t\}} -c\sqrt{n}\1_{\{\tau >t\}}.
\label{eqn_contract1m}
\end{eqnarray} Taking expectations of both sides of (\ref{eqn_contract1m}), we get
$$\E\Big(|X_{t+1}| \1_{\{\tau >t+1\}} +X_{\tau}\1_{\{\tau\leq
t+1\}}\Big) \leq \E\Big(|X_{t}| \1_{\{\tau >t\}} +X_{\tau}
\1_{\{\tau\leq t\}}\Big) -c\sqrt{n}\Pr (\tau >t)\, .$$ Summing over
$t$ from $0$ to $k-1$, we get
\begin{eqnarray}
\E\Big(|X_{k}| \1_{\{\tau >k\}} +X_{\tau} \1_{\{\tau\leq k\}}\Big)
&\leq& |X_0| -\sum_{t=0}^{k-1} c\sqrt{n}\Pr (\tau >t)\nonumber\\
&\leq& |X_0| -k c\sqrt{n}\Pr (\tau >k)\label{eqn_11}.
\end{eqnarray}
We also have $$\E\Big(|X_{k}| \1_{\{\tau >k\}} +X_{\tau}
\1_{\{\tau\leq k\}}\Big)\geq \E\Big(X_{\tau} \1_{\{\tau\leq
k\}}\Big)\geq -An^{2/3}\geq -|X_0|.$$ Combining this with
(\ref{eqn_11}), we have
$$k c\sqrt{n}\Pr (\tau >k)\leq 2 |X_0|$$ which implies the required result.\qed\\

\begin{lemma}\label{momentcontrol}
Let $X$ be a random variable. Then for any $b<0$ and positive
integer $k$, we have
$$\E\Big( |X| \textbf{1}_{(X\leq b)}\Big) \leq\frac{\E|X|^k}{|b|^{k-1}}.$$
\end{lemma}

\noindent{\bf{Proof of Lemma \ref{momentcontrol}:}} We have
\begin{eqnarray*}
\E|X|^k\geq\int_{-\infty}^b(-x)^kdF(x)\geq
|b|^{k-1}\int_{-\infty}^b(-x)dF(x)=|b|^{k-1}\E|X|\textbf{1}_{(X\leq
b)}.
\end{eqnarray*}\qed\\

\noindent{\bf{Proof of Lemma \ref{downdriftlem}:}} For simplicity write again $X_t$ for $X'_t$. Notice that
$$|X_{1}| \1_{\{X_1 \not\in I\}} +X_{1} \1_{\{X_1 \in
I\}}=X_{1} + 2|X_{1}|\1_{\{X_{1}< -A n^{2/3}\}} \, .$$
We first bound $\E X_1$ from above. Recall that in our modified chain we have
$$\E X_1 =\E |\C_{\delta\eps n}| \, .$$
By part (i) of Theorem
\ref{critgnp}, for sufficiently large $A$ and $|X_0| \geq An^{2/3}$, we
have
\begin{eqnarray*}
\E |\C_{\delta\eps n}| \leq 2\frac{X_0}{n}\frac{n+X_0}{2}
-c\Big(\frac{X_0}{n}\Big)^{-2} = X_0 +\frac{X_0^2}{n}
-\frac{cn^2}{X_0^2}.
\end{eqnarray*}
If $X_0\leq \sqrt[4]{\frac{c}{2}} n^{3/4}$, then we have
$\frac{X_0^2}{n}\leq \frac{cn^2}{2X_0^2}$. In this case we have
\begin{equation}\label{eqn_driftpart1}
\E X_1 =\E |\C_{\delta\eps n}|\leq X_0 - \frac{cn^2}{2X_0^2}.
\end{equation}
If $X_0> \sqrt[4]{\frac{c}{2}} n^{3/4}$, then by Lemma
\ref{driftestimate}, we have
\begin{equation}\label{eqn_driftpart2}
\E X_1\leq X_0 -\frac{X_0^2}{6n}.
\end{equation}

Next we bound $\E|X_{1}|\1_{\{X_1< -An^{2/3}\}}$ from above. Let
$$M=\sum _{\C^+_j \neq \C^+_{\delta \eps n}} \eps_j |\C^+_j| +
\sum_{j\geq 1} \eps'_j|\C^-_j|.$$ Then $$X_1 = | \C_{\delta \eps n}
| + M.$$ Since $| \C_{\delta \eps n} |>0$, if $X_1 < -An^{2/3}$, then $M < -An^{2/3}$ and $|X_1| \leq -M$. Thus,
\begin{eqnarray*}
\E\Big(|X_{1}|\1_{\{X_{1}< -A n^{2/3}\}}\Big) \leq \E
\Big((-M)\1_{\{X_{1}< -A n^{2/3}\}}\Big)\leq \E\Big((-M)\1_{\{M< -A
n^{2/3}\}}\Big) .
\end{eqnarray*}
Lemma \ref{momentcontrol} with $k=4$ gives
\be\label{eqn_Mbound}\E\Big((-M)\1_{\{M< -A n^{2/3}\}}\Big)\leq
\frac{\E M^4}{(A n^{2/3})^3}.\ee  We also have
\begin{eqnarray}\E M^4&\leq& \sum _{\C^+_j \neq \C^+_{\delta \eps n}} \E
|\C^+_j|^4 + \sum_{j\geq 1} \E |\C^-_j|^4 + 6 \Big [ \sum _{\C^+_j
\neq \C^+_{\delta \eps n}} \E|\C^+_j|^2 \Big ] \Big [\sum_{j\geq 1}
\E|\C^-_j|^2\Big ]\nonumber\\&+&6\Big [ \sum _{\C^+_j,\C^+_i \neq
\C^+_{\delta \eps n},i\neq j} \E |\C^+_i|^2|\C^+_j|^2+\sum_{i,j\geq
1,i\neq j} \E|\C^-_i|^2|\C^-_j|^2 \Big ]\label{eqn_M4}
.\end{eqnarray} By (ii) of Theorem \ref{critgnp} and Lemma
\ref{subgnpmoments}, we have
\begin{equation}\label{eqn_M41}
\E M^4=  O\Big(\frac{n^6}{X_0^5}\Big) +
O\Big(\frac{n^4}{X_0^2}\Big)=O\Big(\frac{n^4}{X_0^2}\Big)
\end{equation} since $|X_0|\geq An^{2/3}$. Plugging into (\ref{eqn_Mbound}), we have
\begin{eqnarray}\label{eqn_abo}
\E\Big(|X_{1}|\1_{\{X_{1}< -A n^{2/3}\}}\Big)=
O \big ( \frac{n^2}{A^3 X_0^2} \big )
\end{eqnarray}
If $X_0\leq \sqrt[4]{\frac{c}{2}} n^{3/4}$, then combining
(\ref{eqn_abo}) with (\ref{eqn_driftpart1}) for large enough $A$, we
get
$$\E X_{1} + 2\E|X_{1}|\1_{\{X_{1}< -A n^{2/3}\}}\leq X_0 -\frac{cn^2}{4X_0^2} \leq X_0 -c\sqrt{n}.$$
If $X_0> \sqrt[4]{\frac{c}{2}} n^{3/4}$, then combining
(\ref{eqn_abo}) with (\ref{eqn_driftpart2}) for large enough $A$, we
get $$\E X_{1} + 2\E|X_{1}|\1_{\{X_{1}< -A n^{2/3}\}}\leq X_0
-\frac{X_0^2}{6n}+O(A^{-3}\frac{n^2}{X_0^2})\leq X_0-c\sqrt{n}.$$
Combining these two cases finishes the proof. \qed

\subsection{The lower bound on the mixing time}\label{critlowerbound}
Recall that in this section $X_t$ is the original magnetization
chain we defined in (\ref{nextstep}).

\noindent{\bf{Proof of the lower bound of part (ii) of Theorem
\ref{mainthm}:}} Suppose $X_t'$ is a modified magnetization chain
and $\pi'$ is its stationary distribution. By Theorem
\ref{thm_statregime}, we can choose an interval $[a_1 n^{3/4},\,
a_2 n^{3/4}]$ with $0<a_1<a_2$ such that
\begin{equation}\label{eqn_pi34}
\pi'(a_1 n^{3/4}, a_2 n^{3/4})> \frac{3}{4}.
\end{equation} Suppose $X_0'=3a_2n^{3/4}$. By Theorem \ref{cor_1over2}, there exists a constant $k$ such that
$$\P(\tau >kn^{1/4} )\geq \frac{1}{2}$$
where $\tau$ is the first time that $X_t$ exit $[a_2 n^{3/4}, 4a_2
n^{3/4}]$. This implies $$\P\Big(X_{kn^{1/4}}'\geq
a_2 n^{3/4}\Big)\geq \frac{1}{2} \, .$$ By Theorem
\ref{thm_statregime}, we have $\pi'(-a_2 n^{3/4}, -a_1 n^{3/4})$
converges to $0$ as $n$ goes to infinity. Also, by
(\ref{nextstep1}), we have
$\P(X_{kn^{1/4}}'\in[-a_2n^{3/4},-a_1n^{3/4}])= O(n^{-1/12})$.
Combining these, we get that
$$\pi'[(a_1n^{3/4},a_2n^{3/4}),(-a_2n^{3/4},-a_1n^{3/4})]-\P(X_{kn^{1/4}}'\in[(a_1n^{3/4},a_2n^{3/4}),(-a_2n^{3/4},-a_1n^{3/4})])>\frac{1}{4}$$
for large enough $n$. Recall that $X_t\stackrel{d}=|X_t'|$, so this
is equivalent to
$$\pi(a_1n^{3/4},a_2n^{3/4})-\P(X_{kn^{1/4}}\in(a_1n^{3/4},a_2n^{3/4}))>
\frac{1}{4},$$ i.e., \be\label{totalva}\Big\|
X_{kn^{1/4}}-\pi\Big\|_{TV}> \frac{1}{4}.\ee \qed

\section{Fast mixing of the Swendsen-Wang process on trees}

In this section we provide an upper bound estimate of the mixing
time of the Swendsen-Wang process on any tree with $n$ vertices. We
will prove in a more general setting for The Swendsen-Wang process
for the $q$-state ferromagnetic Potts model. Recall that Ising model
is the case $q=2$.

For any given graph $G=(V,E)$, consider the set $S=\{0,1\}^{|E|}$ of
all edge configuration $\eta:E\rightarrow\{0,1\}$. We consider the
following Markov chain $\sigma_t$ on $S$. At each step, we first
color each component independently and uniformly from the $q$
colors. Then we add all edges that connect vertices with the same
color. Finally, delete each existing edge with probability $(1-p)$
to get a new state in $S$. It is easy to see that this process the
dual of the Swendsen-Wang process for the $q$-state ferromagnetic
Potts model on vertices configurations and the stationary
distribution of $\sigma_t$ is the random cluster model. For any two
Swendsen-Wang chains, if we can couple the corresponding edge models
so that they are the same(i.e., they have same clusters) at some
time, we therefore couple the original Swendsen-Wang process at the
same time. Consequently, any upper bound of the mixing time of this
edge model implies the same upper bound on ferromagnetic Potts
model.

There is an exploration process on trees to present $\sigma_t$.
Notice that on trees each edge with state $0$ connects two separate
components. For any given $\eta\in S $, we color each components
independently and uniformly from the $q$ colors, starting from the
root. We add edges connects vertices with the same color. Notice
that this procedure is equivalent to setting every edge originally
has configuration $0$ with configuration $1$ with probability
$\frac{1}{q}$ and maintain configuration $0$ otherwise. Thus, the
process $\sigma_t$ can be described as follows: First change each
edge of $0$ to $1$ with probability $\frac{1}{q}$ and stay $0$
otherwise, independently for each of them. Then, change each edge of
$1$ ,including those who have changed from $0$ to $1$ in the
previous step, to $0$ with probability $1-p$, and stay $1$
otherwise, independently for each of them. Each bit evolves
independently as a Markov chain on $\{0,1\}$, with transition matrix
\begin{equation}\label{eqn_matrix}
\p=\left(%
\begin{array}{cc}
  1-\frac{p}{Q} & \frac{p}{Q} \\
  1-p & p \\
\end{array}%
\right).
\end{equation}

\noindent{\bf{Proof of Theorem \ref{thm_tree}:}} The transition
matrix (\ref{eqn_matrix}) gives that we can couple every single edge
with probability at least $1-p+\frac{p}{q}\geq \frac{1}{q}$. Using
the path coupling method of Bubley and Dyer (see Theorem 14.6 and
Corollary 14.7 of \cite{LPW}), we have $$t_{mix}\leq \frac{\log
n+\log 4}{-\log p(1-{1\over q})}.$$\qed

\section*{Acknowledgements} We are very grateful to Jian Ding for carefully reading the paper and for numerous useful suggestions.

\end{document}

%% file: gnplemmas.tex
\section{Random graph estimates}\label{secgnplemmas}
In this section we prove some facts about random graphs which will be used in the proof. These lemmas might also be of sperate interests in random graph theory. Recall that $G(m,p)$ is obtained from the
complete graph on $m$ vertices by retaining independently each
edge with probability $p$ and deleting it with probability $1-p$.
We denote by $\C_j$ the $j$-th largest component of $G(m,p)$.

\subsection{The exploration process}\label{exploresec} We recall an exploration
process, due to Karp and Martin-L\"of (see \cite{K} and \cite{M}),
in which vertices will be either {\em active, explored} or {\em
neutral}. After the completion of step $t \in \{0,1,\ldots , m\}$
we will have precisely $t$ explored vertices and the number of the
active and neutral vertices is denoted by $A_t$ and $N_t$
respectively. Fix an ordering of the vertices $\{v_1, \ldots,
v_m\}$. In step $t=0$ of the process, we declare vertex $v_1$
active and all other vertices neutral. Thus $A_0=1$ and
$N_0 = m-1$. In step $t \in \{1,\ldots,m\}$, if $A_{t-1}>0$, then let
$w_t$ be the first active vertex; if $A_{t-1}=0$, let $w_t$ be the
first neutral vertex. Denote by $\eta_t$ the number of neutral
neighbors of $w_t$ in $G(m,p)$, and change the status of these
vertices to active. Then, set $w_t$ itself explored.

Denote by $\F _{t}$ the $\sigma$-algebra generated by $\{\eta_1,
\ldots , \eta_t\}$. Observe that given $\F_{t-1}$ the random
variable $\eta_t$ is distributed as Bin$(N_{t-1} - {\bf
1}_{\{A_{t-1}=0\}},p)$ and we have the recursions \be \label{nrec}
N_t = N_{t-1} - \eta_t - {\bf 1}_{\{A_{t-1}=0\}} \, , \qquad t
\leq m \, ,\ee and

\be \label{arec}
 A_t= \left \{
\begin{array}{ll}
A_{t-1} + \eta_t - 1, & A_{t-1} > 0 \\
\eta_t, & A_{t-1} = 0 \, , \qquad t
\leq m \, . \\
\end{array} \right .
\ee  As every vertex is either neutral, active or explored,

\be \label{nrec2} N_t = m - t - A_t \, , \qquad t \leq m \, . \ee

At each time $j\leq m$ in which $A_j = 0$, we have finished
exploring a connected component. Hence the random variable $Z_t$
defined by
$$ Z_t = \sum _{j=1}^{t-1} {\bf 1}_{\{A_j = 0\}} \, ,$$
counts the number of components completely explored by the process
before time $t$. Define the process $\{Y_t\}$ by $Y_0 = 1$ and
$$ Y_t = Y_{t-1} + \eta_t - 1 \, .$$ By (\ref{arec}) we have that
$Y_t = A_t - Z_t$, i.e. $Y_t$ counts the number of active vertices
at step $t$ minus the number of components completely explored
before step $t$. 
%
\begin{lemma}\label{ytbound}
For any $t$ we have \be\label{ytupper}Y_t \stackrel{d}{\leq}
\hbox{\rm Bin}(m-1,1-(1-p)^t)+1-t,\ee and \be\label{ytlower}Y_t \stackrel{d}{\geq}
\hbox{\rm Bin}(m-t-1,1-(1-p)^t)+1-t.\ee
\end{lemma}
\noindent{\textbf{Proof.}} For each vertex $v$ at each step of the process we examine precisely one of its edges emanating from it unless the vertex is active or explored at this step. Thus, all the vertices for which the process discovered an open edge emanating from them between time $1$ and $t$ are active, except for at least $t$ of them which are explored. The probability of a vertex having no open edges explored from it between time $1$ and $t$ is precisely $(1-p)^t$. This shows (\ref{ytupper}).

The reason this bound is not precise is that it is possible that a neutral vertex turns to be active because there were no more active vertices at this step of the exploration process. This, however, can only happen at most $t$ times between time $1$ and $t$ and this gives the lower bound (\ref{ytlower}).
%
%
%
%
\qed \\

At each step we marked as
explored precisely one vertex. Hence, the component of $v_1$ has
size $\min \{t \geq 1 : A_t = 0\}$. Moreover, let $t_1 < t_2 \ldots$
be the times at which $A_{t_j}=0$; then $(t_1, t_2 - t_1, t_3 - t_2,
\ldots)$ are the sizes of the components. Observe that $Z_{t} =
Z_{t_j}+1$ for all $t \in \{t_j+1, \ldots, t_{j+1}\}$. Thus
$Y_{t_{j+1}} = Y_{t_j}-1$ and if $t \in \{t_j +1, \ldots,
t_{j+1}-1\}$ then $A_t
> 0$, and thus $Y_{t_{j+1}}< Y_t$. By induction we conclude that
$A_t=0$ if and only if $Y_t < Y_s$ for all $s<t$. In other words $A_t=0$ if and only if $\{Y_t\}$ has hit a new record minimum at time $t$. By
induction we also observe that $Y_{t_j} = -(j-1)$ and that for $t
\in \{t_j +1, \ldots t_{j+1}\}$ we have $Z_t = j$. Also, by our
previous discussion for $t \in \{t_j +1, \ldots t_{j+1}\}$ we have
$\min _{s \leq t-1} Y_s = Y_{t_j} = -(j-1)$, hence by induction we
deduce that $Z_t = - \min_{s \leq t-1} Y_s + 1$. Consequently, \be
\label{formula} A_t = Y_t - \min _{s\leq t-1} Y_s + 1 \, .\ee

\begin{lemma} \label{lowerbound} For all $p\leq {2 \over m}$ there
exists a constant $c>0$ such that for any integer $t>0$,
$$ \prob \Big (N_{t} \leq m - 5t \Big) \leq e^{-ct} \, .$$
Where we recall $N_t$ is the number of neutral points in exploration process at time $t$.
\end{lemma}
\noindent The proof of Lemma \ref{lowerbound} can be found in Lemma 3 of
\cite{NP}. \\

\subsection{Random graph lemmas for non-critical cases}

Let $G(m,p)$ be the random graph where $p=\frac{\theta}{m}$.

\begin{lem}\label{lem_subdeviation}
Suppose $\theta<1$ is a constant. Then we have
\begin{equation}\label{eqn_subdeviation}
\E(\sum_{j\geq 1} |\C_j|^2)\leq\frac{\m}{1-\theta}.
\end{equation}
\end{lem}

\noindent{\bf{Proof.}} Observe that
$\sum_{j\geq 1} |\C_j|^2=\sum_v|\C(v)|$ since in the right hand side each component $\C(v)$ is counted precisely $|\C(v)|$ times. Hence
\be\label{eqn_useful}\E(\sum_{j\geq 1}
|\C_j|^2)=\E\sum_v|\C(v)|=mE|\C(v)|.\ee In the exploration process, we
can couple $Y_t$ with a process $W_t$ with i.i.d. increment of
$bin(n,\theta/m)-1$ and $W_0=1$ such that $W_t\geq Y_t$. Thus, the
hitting time of $0$ for $Y_t$ which equals to $|\C(v)|$ is bounded
from above by the hitting time of $0$ for $W_t$. For $W_t$, we have
$\E\tau=1/(1-\theta)$ by Wald's Lemma. This concludes the proof.\qed \\

For $\theta>1$ let $\beta=\beta(\theta)$ be the unique positive
solution of the equation
\begin{equation}\label{eqn_giantsize}
1-e^{-\theta x}=x.
\end{equation} In \cite{P} it was proved that
$\frac{|\C_1|-\beta m}{\sqrt{m}}$ converges in distribution to a
normal distribution. We were unable to deduce from that result
moderate deviation estimates, and we provide them in the following
lemma.

\begin{lem}\label{realsupergnp} There exists constants $c=c(\theta)>0$ and universal constant $C$
such that for any $A>0$ we have
\begin{equation}\label{realsupercriticalcon}
\P(\big||\C_1|-\beta m\big|\geq A\sqrt{m})\leq
Ce^{-cA^2}.
\end{equation}
\end{lem}

\noindent{\bf{Proof.}} Assume $A\leq\sqrt{m}$ otherwise this probability is $0$. Let $\xi=\xi(\theta)>0$ be a large constant that we will determine later. We
will show that for some $c>0$ \be\label{deviation1}\P(Y_{\beta
m+A\sqrt{m}}\geq -cA\sqrt{m})\leq e^{-cA^2},\ee and that
\be\label{deviation2}\P(\bigcup_{cA\sqrt{m}\leq t\leq \beta
m-\xi A\sqrt{m}}Y_t<0)\leq Ce^{-cA^2} \, . \ee If these two events do not occur, then there
exists a component of size in $[\beta m-(\xi+c)A\sqrt{m}, \beta
m+A\sqrt{m}]$. The remaining graph is a subcritical random graph
and it is a classical result that the probability that it contains a
component of size $\Theta(m)$ decays exponentially in $m$, and this
will conclude the proof.

The proof of (\ref{deviation1}) is based on the stochastic upper bound of
$Y_t$ in (\ref{ytupper}). Plugging in $t=\beta
m+A\sqrt{m}$ and using the fact that $1-x\geq e^{-x-x^2}$ for small enough $x$ we get
\begin{eqnarray*}
\P(Y_{\beta
m+A\sqrt{m}}\geq -cA\sqrt{m}) &\leq &\P\big({\rm{Bin}}(m, 1-(1-{\theta\over m})^{\beta m+A\sqrt{m}})\geq\beta m+(1-c)A\sqrt{m}\big)\\&\leq& \P\big({\rm{Bin}}(m, 1-e^{-\theta \beta - A \theta m^{-1/2} - A \theta^2 m^{-3/2}})\geq \beta m + (1-c) A\sqrt{m} \big ) \, .
\end{eqnarray*}
A quick calculation using the fact that $1-e^{-\theta \beta} = \beta$ gives that the expected value of this binomial random variable is at most
$$ \beta m  + A \theta e^{-\theta \beta} \sqrt{m} + O(1) \, .$$
Since $\theta e^{-\theta \beta} < 1$ it follows that we can choose $c$ so small so that this expectation is less than $\beta m +(1-2c)A \sqrt{m}$, and then Azuma-Hoeffding inequality (see for instance Theorem 7.2.1 of \cite{AS}) gives that
$$ \P(Y_{\beta m+cA\sqrt{m}}\geq -A\sqrt{m}) \leq e^{-cA^2} \, .$$

We now turn to prove (\ref{deviation2}). For this we will divide $[cA\sqrt{m},\beta
m-\xi A\sqrt{m}]$ into two subintervals $[cA\sqrt{m},\delta\beta m]$ and $[\delta\beta m, \beta m-\xi A\sqrt{m}]$ where $\delta>0$ is a small constant that will be chosen later. For convenience write $\alpha={t\over m}$. For any $t \in[cA\sqrt{m},\delta\beta m]$ we
have by (\ref{ytlower}) and the fact that $1-x\leq e^{-x}$ for all $x\geq 0$ that
\begin{eqnarray*}
\P(Y_t<0)&\leq&\P\Big({\rm{Bin}}\big((1-\alpha)m,1-(1-{\theta\over m})^{\alpha m}\big) \leq \alpha m\Big)\\&\leq&\P\Big({\rm{Bin}}\big((1-\alpha)m,1-e^{-\theta\alpha}\big)\big) \leq  \alpha m\Big) \, .
\end{eqnarray*}
Since $1-e^{-x} \geq x-x^2$ for all $x \geq 0$ we deduce that the expectation of the last binomial is at least
$$ (1-\alpha)(\theta \alpha - \theta^2 \alpha^2) > \alpha \, ,$$
since $\theta>1$ when $\alpha = t/m \leq \delta$ and $\delta=\delta(\theta)>0$ is chosen small enough. By a standard large deviation estimate (see for instance, Corollary A.1.14 of \cite{AS}), we have that $$\P(Y_t<0)\leq e^{-c \alpha m}\, ,$$ for some $c=c(\theta)$ and all $t\in[cA\sqrt{m},\delta\beta m]$. It follows from the union bound that
\be\label{deviation21}\P(\bigcup_{cA\sqrt{m}\leq t\leq \delta\beta
m}Y_t<0)=O(e^{-cA\sqrt{m}}).\ee

For the interval $[\delta\beta m,\beta m-\xi A\sqrt{m}]$ we will use the process $\widetilde{Y}_t$ which approximates $Y_t$ introduced by Bollobas and Riordan \cite{BR}. We write $$D_t=\E(\eta_t-1|F_{t-1}),$$ and define $$\Delta_t=\eta_t-1-D_t.$$ Let $y_t=m-t-m(1-p)^t$ and define the approximation process by \be \label{approxprocess} \widetilde{Y}_t=y_t+\sum_{i=1}^t(1-p)^{t-i}\Delta_i \, .\ee
In \cite{BR} Lemma 3 it is proved that for any $p >0$ and any $1\leq t\leq n$ we have
\begin{equation} \label{quoteBR}
|Y_t-\widetilde{Y}_t|\leq ptZ_t.
\end{equation}
Put $\tau=\min\{t\geq \delta\beta n, A_t=0\}.$ We have
\begin{eqnarray}\label{sblit}
\P(\tau<\beta m-\xi A\sqrt{m})&\leq&\P(|Y_{\tau}-\widetilde{Y}_{\tau}|\geq\theta A\sqrt{m})\nonumber\\&+&\P(|Y_{\tau}-\widetilde{Y}_{\tau}|<\theta A\sqrt{m},\tau<\beta m-\xi A\sqrt{m}).
\end{eqnarray}
By (\ref{quoteBR}) the first term has the upper bound
\begin{eqnarray*}
\P(|Y_{\tau}-\widetilde{Y}_{\tau}|\geq\theta A\sqrt{m})\leq\P(Z_{\tau}\geq A\sqrt{m})=O(e^{-cA\sqrt{m}}),
\end{eqnarray*} since $Z_{\tau}\geq A\sqrt{m}$ implies that there exists at least one time $t$ in $[A\sqrt{m},\delta\beta m]$ such that $Y_t<0$. The bound follows immediately from our estimate in (\ref{deviation21}).

To bound the second term of (\ref{sblit}) observe that on $[\delta\beta m,\beta m-\xi A\sqrt{m}]$, the minimum of $y_t$ is attained at the right end of the interval with value $(1-\theta e^{-\theta\beta})\xi A\sqrt{m}(1+o(1))$. Thus if we choose $\xi=\xi(\theta)$ large enough such that $(1-\theta e^{-\theta\beta})\xi>\theta$ and write $c=(1-\theta e^{-\theta\beta})\xi-\theta$, we have
\begin{eqnarray}\label{largedeviation}
\P(|Y_{\tau}-\widetilde{Y}_{\tau}|<\theta A\sqrt{m},\tau<\beta m-\xi A\sqrt{m})\leq\P\big(\sum_{i=1}^{\tau}(1-p)^{\tau-i}\Delta_i<-c A\sqrt{m}\big),
\end{eqnarray}since $Y_{\tau}\leq 0$ by definition.
Notice that $\tau$ is at most $m$, thus it suffices to bound from above $$\P\Big(\max_{1\leq t\leq m}(-\sum_{i=1}^t(1-p)^{t-i}\Delta_i)>c A\sqrt{m}\Big).$$ Notice that $(1-p)^t=\Theta(1)$, hence it is equivalent to bound $$\P\Big(\max_{1\leq t\leq m}(-\sum_{i=1}^t(1-p)^{-i}\Delta_i)>c A\sqrt{m}\Big).$$
Let $a>0$ be a small number (eventually we will take $a=\Theta(m^{-1/2})$).  Direct computation and the fact that $1+x>e^{x-x^2}$ for negative $x$ when $|x|$ is small enough and some Taylor expansion yield \begin{eqnarray}\label{submartingale}\E(e^{-a(1-p)^{-i}\Delta_i}|F_{i-1})&=&(1+p(e^{-a(1-p)^{-i}}-1))^{N_{i-1}-{\bf1}_{\{A_{i-1}=0\}}} e^{a(1-p)^{-i}p(N_{i-1}-{\bf1}_{\{A_{i-1}=0\}})}\nonumber\\&\geq&e^{\frac{a^2p}{3}(N_{i-1}-{\bf1}_{\{A_{i-1}=0\}})(1-p)^{-2i}}\geq1,\end{eqnarray} when $a$ is small enough. Thus we conclude $e^{-a\sum_{i=1}^t(1-p)^i\Delta_i}$ is a submartingale. By Doob's maximal inequality (see \cite{D}) we have $$\E\big(\max_{1\leq t\leq m}e^{-a\sum_{i=1}^t(1-p)^i\Delta_i}\big)^2\leq 4\E e^{-2a\sum_{i=1}^m(1-p)^i\Delta_i}.$$ On the other hand, the fact that $1+x\leq e^x$ for all $x$ yields
\begin{eqnarray}\label{submartingale1}\E(e^{-a(1-p)^{-i}\Delta_i}|F_{i-1})\leq e^{(\frac{a^2}{2}+O(a^3))p(N_{i-1}-{\bf1}_{\{A_{i-1}=0\}})(1-p)^{-2i}}.\end{eqnarray} Since $N_t\leq m$ we get $$\E e^{-2a\sum_{i=1}^m(1-p)^i\Delta_i}\leq4 e^{C m a^2},$$ where $C=C(\theta)$. By Markov's inequality, we have $$\P\big(\max_{1\leq t\leq m}(-\sum_{i=1}^t(1-p)^{-i}\Delta_i)>c A\sqrt{m}\big)\leq e^{C ma^2-2ac A\sqrt{m}} \, .$$ Choosing $a=\frac{c A\sqrt{m}}{Cm}$ to minimize the right hand side, we conclude $$\P\big(\max_{1\leq t\leq m}(-\sum_{i=1}^t(1-p)^{-i}\Delta_i)>\theta A\sqrt{m}\big)\leq4 e^{-cA^2},$$ for some $c=c(\theta)$ which is a continuous function of $\theta$, and this concludes the proof.\qed

\begin{corollary}\label{lem_giant}
Suppose $\theta\in[a,b]$ where $a>1$. Then, there exists a constant
$C=C(a,b)$ such that for $G(m, {\theta\over m})$ we have
\begin{equation}\label{eqn_firstmoment}
\Big|\mathbb{E} |\mathcal{C}_1 |-\beta (\theta) m\Big| \leq
C \sqrt{m},
\end{equation}
\end{corollary}

\noindent{\bf{Proof.}} It follows immediately by integrating Lemma \ref{realsupergnp}. \qed

\begin{corollary}\label{lem_superdeviation}
Suppose $\theta\in[a,b]$ where $a>1$. There exists a constant
$C=C(a,b)$ such that for $G(m,{\theta\over m})$ we have
\begin{equation}\label{eqn_superdeviation} \E(\sum_{j\geq 1}
|\C_j|^2)\leq (\E|\C_1|)^2+Cm.
\end{equation}
\end{corollary}

\noindent{\bf{Proof.}} Notice
that \be\label{difference11}\E\sum_{j\geq 1}
|\C_j|^2-(\E|\C_1|)^2=\Big(\E|\C_1|^2-(\E|\C_1|)^2\Big)+\E\sum_{j\geq
2} |\C_j|^2.\ee We have that
$$\E|\C_1|^2-(\E|\C_1|)^2=\E(|\C_1|-\E|\C_1|)^2\leq \E(|\C_1|-\beta
m)^2.$$ By integrating Lemma
\ref{realsupergnp} we get
$$\E|\C_1|^2-(\E|\C_1|)^2\leq Cm.$$ For
supercritical random graph $G(m,{\theta\over m})$, it is a classical
result that $|\C_1|\in((\beta-\eps)m,(\beta+\eps)m)$ with
probability at least $1-e^{-c_{\eps}m}$ for fixed $\eps$.
Conditioned on this and the vertex set of $\C_1$, the other
components are distributed as $G(m-|\C_1|, {\theta\over m})$ (which
is subcritical) restricted to the event that it does not contain any
component larger than $|\C_1|$. This event happens with
probability at most $e^{-cm}$. Thus we obtain $$\E(\sum_{j\geq 2}
|\C_j|^2)\leq (1+o(1))\frac{m}{1-(1-\beta)\theta},$$ \qed

\begin{lem}\label{lem_isolated}
Let $M=\sum_{v\in V} \1_{\{v\ is\ isolated\}}$ be the number of
isolated vertices in $G(m,\theta/m)$ where $\theta>0$ is a constant.
There exists a constant $C>0$ such that $$\P(M\geq
Cm)=1-O(\frac{1}{m}).$$
\end{lem}
\noindent{\bf{Proof.}} We have
$$\E M=\sum_{v\in V}\P(v\,\,{\rm{ is\,\, isolated.}})=\m (1-\frac{\theta}{\m})^{\m-1},$$ and
\begin{eqnarray*}
\E M^2 &=& \sum_{v\in V}\P(v\,\, \rm{is\,\, isolated.})+\sum_{v,w\in
V}\P(v,w\,\, \rm{are\,\, both\,\, isolated.})\\
&=& \m(1-\frac{\theta}{\m})^{\m-1}+\m(\m-1)
(1-\frac{\theta}{\m})^{2\m-3}.
\end{eqnarray*}
Thus, we obtain
$$\E (M-\E M)^2=O(m).$$
By Markov's inequality,
$$\P(M\leq \half\E M)\leq \P\Big((M-\E M)^2\geq \frac{1}{4}(\E M)^2\Big)\leq
\frac{\E (M-\E M)^2}{\frac{1}{4}(\E M)^2}=O(\frac{1}{\m}).$$ Since
$\E M=\Theta(m)$, we finished the proof.\qed \\


\subsection{Random graph lemmas for the near-critical case} \label{nearcritsec}

In \cite{PW}, Pittel and Wormald study the near-critical random graph $G(m,p)$ where $p={1 +\eps \over m}$ with $\eps = o(1)$ but $\eps^3 m \to \infty$. A direct corollary of Theorem 6 of \cite{PW} shows that in this regime ${|\C_1| - 2\eps m \over \sqrt{m/\eps}}$ converges in distribution to a normal random variable (see also \cite{BR} for a recent simple proof of this fact), and that a local central limit theorem holds. Unfortunately, once cannot deduce from that precise bounds on the average size of $|\C_1|$ and moderate deviations estimates on $|\C_1| - 2\eps m$. The following two theorems give these estimates.

\begin{theorem} \label{supergnp1} Consider $G(m,p)$ with $p={1 + \eps \over
m}$ where $\eps=o(1)$ and there exists a large constant $A>0$ such
that $\eps^3 m \geq A \log m$. Then we have that
$$ \E |\C_1| \leq 2\eps m - {8 \over 3}\eps^2 m + O(\eps^3 m),$$ and
there exists a constant $C>0$ such that
$$ \E |\C_1| \geq 2\eps m - C(\eps^{-2} + \eps^2 m) \, .$$
\end{theorem}


\begin{theorem} \label{gnpout2} Consider $G(m,p)$ with $p={1 + \eps \over
m}$ where $\eps^3 m \geq 1$. Then there exists some $c>0$ such that
$$ \prob \Big ( \Big | |\C_1| - 2\eps m \Big | > A \sqrt{ {m \over
\eps} } \Big ) = O\big (e^{-cA^2} \big )  ,$$
for any $A$ satisfying $2 \leq A \leq \sqrt{\eps^3 m}/10$. \\
\end{theorem}

\begin{corollary} \label{supergnpmoments} Consider $G(m,p)$ with $p={1 + \eps \over
m}$ where $\eps^3 m \geq 1$, then
$$\E  \Big | |\C_1| - 2\eps m \Big |^k  \leq C \Big ( {m \over
\eps } \Big ) ^{k/2} \, .$$
\end{corollary}

\begin{theorem}\label{giantclt} For any large constant $A$ and small $\delta>0$ there exists a constant $q_1(A,\delta)>0$ such that the following hold. Consider $G(m,p)$ with $p={1 + \eps \over
m}$ where $\eps\in [A^{-1} m^{-1/4}, A m^{-1/4}]$, then
$$\P\big (|\C_1|\in [2\eps m-\delta m^{5/8},\,2\eps m+\delta m^{5/8}] \big )\geq
q_1> 0 \, .$$
\end{theorem}
\noindent Theorem \ref{giantclt} is a direct corollary of Theorem 6 of \cite{PW} which provides a central limit theorem for the giant component. Next we provide some moment estimates of component sizes in the subcritical and supercritical regime.

\begin{theorem} \label{subgnpmoments} Consider $G(m,p)$ with $p={1 -\eps
\over m}$. Then we have
\begin{enumerate}
\item[(i)] $\E \sum _{j \geq 1} |\C_j|^k = O(m \eps^{-2k+3})$ for any fixed $k\geq 2$,
\item[(ii)] $\E\sum_{i,j}|\mathcal{C}_i|^2|\mathcal{C}_j|^2 = O(m^2 \eps^{-2})$,
\item[(iii)] If $\eps^3 m \geq 1$, then $\E \sum _{j \geq 1} |\C_j|^2 \geq c m \eps^{-1}$.
\end{enumerate}
\end{theorem}

\begin{theorem} \label{gnp2moments} Consider $G(m,p)$ with $p={1 + \eps \over
m}$ with $\eps >0$ and $\eps^3m\geq 1$ for large $m$. Then we have
\begin{enumerate}
\item[(i)] $\E |\C(v)|^k = O(\eps^{k+1} m^k)$, for any fixed $k \geq 2$.
\item[(ii)] $\E \sum _{j \geq 2} |\C_j|^k = O(m \eps^{-2k+3})$,
\item[(iii)] $\E\sum_{i,j\geq2}|\mathcal{C}_i|^2|\mathcal{C}_j|^2 = O(m^2 \eps^{-2})$.
\end{enumerate}
\end{theorem}
%

\begin{theorem} \label{gnpsmalls} Consider $G(m,p)$ with $p={1 - \eps \over
m}$ where $\eps \in [A^{-1} m^{-1/4},\,A m^{-1/4}]$. Then, for any
small positive constant $\delta$, there exist $K=K(A,\delta)$ and
$q_2=q_2(A,\delta)$ such that
$$\P\Big(\sum _{|\C_j |\leq\delta\sqrt{m}} |\C_j|^2 \geq K m^{5/4}\Big)\geq q_2>0.$$
\end{theorem}

In the following theorem we derive estimates on the expected cluster size valid as long as $\eps^3 m \geq 1$. We believe these estimates should hold for $\C_1$ but we were not able to prove that. The difficulty rises because when $\eps^3 m$ is large but does not grow at least logarithmically, it is hard to rule out the possibility that $\C_1$ is discovered after time $\delta \eps n$ for some fixed $\delta>0$. Luckily, for the main proof it suffices to have these estimate for $\C_{\delta \eps m}$, that component discovered at time $\delta \eps m$, rather than $\C_1$. This becomes evident in the proof of Theorem \ref{pushdownuse}.

\begin{theorem}\label{critgnp} Consider $G(m,p)$ with $p={1 + \eps \over m}$ and assume $\eps^3 m \geq 1$. For some fixed $\delta>0$ let $\C_{\delta \eps m}$ be the component which is discovered by the exploration process at time $\delta \eps
m$ (in other words, the length of the excursion of $Y_t$ containing
the time $\delta \eps m$). Then there is some small value of $\delta>0$ such that
\begin{itemize}
\item[(i)] $\E |\C_{\delta \eps m}| \leq 2\eps m - c\eps^{-2}\,$.
\item[(ii)] $ \E \sum _{\C_j \neq \C_{\delta \eps m}} |\C_j|^k
\leq C m\eps^{-2k+3}\, ,\ \hbox{ {\rm for} }\ k=2,4 \, ,$
\end{itemize}
where $C$ and $c$ are positive universal constants.
\end{theorem}

Before proceeding to the proofs of the theorems stated in this section, we first require some preparations about processeses with i.i.d. increments.

\subsubsection{Processes with i.i.d. increments} Fix some small $\eps>0$ and let $p={1+ \eps \over m}$ for some integer $m>1$. Let $\{\beta _j\}$ be a sequence of random
variables distributed as Bin$(m,p)$. Let $\{W_t\}_{t \geq 0}$ be a
process defined by
$$ W_0 = 1, \qquad W_t = W_{t-1} + \beta_t - 1 \, .$$
Let $\tau$ be the hitting time of $0$, i.e.
$$ \tau = \min _t \{ W_t =  0\} \, .$$

\begin{lemma} \label{tail} We have
\be \label{survprob} \prob ( \tau = \infty ) = 2 \eps - {8 \over 3}
\eps^2 + O(\eps^3) \, ,\ee and there exists constant $C, c>0$
such that for all $T\geq \eps^{-2}$ we have \be \label{tailbd} \prob
(T \leq \tau < \infty ) \leq C\Big ( \eps^{-2} T^{-3/2} e^{-
{(\eps^2 - c \eps^3) T \over 2}}\Big ) \, .\ee
\end{lemma}

We say that $t_0$ is a {\em record minimum} of $\{W_t\}$ if $W_t >
W_{t_0}$ for all $t < t_0$.

\begin{lemma} \label{numrecmin} Denote by $Z^w$ the number of
record minima of $W_t$. Then
$$ \E Z^w = {\eps^{-1} \over 2} + O(1) \, , \,\, \and \,\,\E (Z^w)^2 = O(\eps^{-2}) \, .$$
\end{lemma}

\begin{lemma} \label{recmin} Denote by $\gamma$ the random variable
$$ \gamma = \max \, \{ t \, : \, t \hbox{ {\rm is a record
minimum of} } W_t \} \, .$$ Then we have
$$ \E \gamma = O(\eps^{-2}) \, .$$
\end{lemma}

\noindent For the subcritical case we have the following.

\begin{lemma} \label{tail2} Assume $\eps<0$ in the previous setting. There exists constant $C_1, C_2, c_1, c_2> 0$ such that for all $T\geq \eps^{-2}$ we have
$$ \prob (\tau \geq T ) \leq C_1\Big ( \eps^{-2} T^{-3/2} e^{- {(\eps^2 - c_1 \eps^3) T \over
2}}\Big ) \, ,$$ and
$$ \prob (\tau \geq T ) \geq c_2\Big ( \eps^{-2} T^{-3/2} e^{- {(\eps^2 + C_2 \eps^3) T \over
2}}\Big ) \, .$$ Furthermore, for any fixed $k\geq 1$
$$ \E \tau^k = O ( \eps^{-2k+1} ) \, .$$
\end{lemma}

\noindent The proof of Lemma \ref{tail2} can be found in \cite{NP}
Lemma 4.

For the proof of Lemma \ref{tail} we will use the following
proposition due to Spitzer (see \cite{S1}).
\begin{prop} \label{spitzer}
Let $a_0, \ldots, a_{k-1} \in {\mathbb Z}$ satisfy $\sum
_{i=0}^{k-1} a_i = -1$. Then there is precisely one $j \in \{0,
\ldots, k-1\}$ such that for all $r \in \{0, \ldots, k-2\}$
$$ \sum _{i=0}^{r} a_{(j+i) {\rm \  mod \ } k} \geq 0 \, .$$
\end{prop}

\noindent {\bf Proof of Lemma \ref{tail}.} Let $\beta$ be a random
variable distributed as Bin$(m,p)$ and let $f(s) = \E s^{\beta}$.
It is a classical fact (see \cite{AN}) that $1- \prob ( \tau =
\infty )$ is the unique fixed point of $f(s)$ in $(0,1)$. For $s
\in (0,1)$ we have
$$ \E s^{\beta} = \Big [ 1 - p(1-s) \Big ]^m = 1 - (1+\eps)(1-s) +
{ (1+\eps)^2 (1-s)^2 \over 2 } -  {(1+\eps)^3(1-s)^3 \over 6} +
O\Big ( (1-s)^4 \Big ) \, ,$$ since $(1-x)^m = 1-mx+{m^2 x^2 \over
2}- {m^3 x^3\over 6} + O(m^4x^4)$. Write $q = 1-s$ and put $\E
s^{\beta} = s$. We get that
$$ 1-(1+\eps)q + {(1+2\eps)q^2 \over 2} - {q^3 \over 6} + O(q^4) +
O(\eps q^3) + O(\eps^2 q^2) = 1-q \, .$$ Solving this gives that
that $q = 2\eps - {8 \over 3}\eps^2 + O(\eps^3)$, as required.

We now turn to proving (\ref{tailbd}). By Proposition
\ref{spitzer}, $\prob (\tau = t) = {1 \over t} \prob (W_t = 0)$.
Since $\sum _{j=1}^t \beta_j$ is distributed as a Bin$(nt,p)$
random variable we have
$$ \prob (W_t = 0) = { mt \choose t-1}p^{t-1}(1-p)^{mt-(t-1)} \, .$$
Replacing $t-1$ with $t$ in the above formula only changes it by a
multiplicative constant which is always between ${1/2}$ and $2$. A
straightforward computation using Stirling's approximation gives

\be \label{stirling} \prob (W_t = 0) = \Theta \Big \{ t^{-1/2}
(1+\eps)^t \Big (1+ {1 \over m-1} \Big )^{t(m-1)}\Big (1 - {1+\eps
\over m} \Big)^{t(m-1)} \Big \} \, .\ee Denote $x = (1+\eps) \Big
( 1 + {1 \over m-1} \Big )^{m-1} \Big (1- {1+\eps \over m} \Big
)^{m-1} $, then
$$ \prob (\tau \geq T) = \sum _{t \geq T} \prob (\tau = t) = \sum _{t \geq
T}{1 \over t} \prob (W_t = 0) = \Theta \Big ( \sum _{t \geq T}
t^{-3/2} x^t \Big ) \, .$$ This sum can be bounded above by
$$ T^{-3/2} \sum _{t \geq T} x^t = T^{-3/2}{x^T \over 1-x} \, .$$
%
Observe that as $m \to \infty$ we have that $x$ tends to
$(1+\eps)e^{-\eps}$. By expanding $e^{-\eps}$ we find that
$$ x = (1+\eps)(1-\eps+{\eps^2 \over 2}) + \Theta(\eps^3) = 1 -
{\eps^2 \over 2} + \Theta(\eps^3) \, .$$ Using this and the
previous bounds on $\prob (\tau = \infty)$ we conclude the proof of
(\ref{tailbd}). \qed \\

\noindent {\bf Proof of Lemma \ref{numrecmin}.} This follows immediately since $Z^w$ is a geometric random variable with success probability
$p=P(\tau=\infty)=2\eps-\frac{8}{3}\eps^2+O(\eps^3)$ by (\ref{survprob}) of Lemma \ref{tail}. \qed \\

\noindent {\bf Proof of Lemma \ref{recmin}.} At each record minimum the process has probability $\Theta(\eps)$ of never going below its current location by (\ref{survprob}) of Lemma \ref{tail}. It is a classical fact that the expected size of each excursion between record minimum, on the event that it is finite, is $O(\eps^{-1})$. Thus, by Wald's Lemma
$$ \E (\gamma) \leq  C\eps^{-1} \E Z^w =O(\eps^{-2})
.$$ \qed

%
%

\subsubsection{Exploration process estimates} In this section we study the process $Y_t$ defined in Section \ref{exploresec} and provide some useful estimates.

\begin{lemma} \label{ztbound} For $p={1 +\eps \over m}$ we have
$$ \prob \Big ( Y_t \geq -45 \eps^2 m \,\, \hbox{ {\rm for
all} } 1 \leq t \leq 3 \eps m \Big ) \geq 1-5 e^{-48\eps^3 m} \, .$$
\end{lemma}
\noindent {\bf Proof.} Denote by $\gamma$ the stopping time
$$ \gamma = \min \{ t \, : \, N_t \leq m-15\eps m \} \, ,$$
and consider the process $\{W_t\}$ which has i.i.d.~ increments
distributed as Bin$(m-15\eps m,p)-1$ and $W_0=1$. Then we can
couple the processes $\{Y_t\}$ and $\{W_t\}$ such that $Y_{t \wedge
\gamma} \geq W_{t \wedge \gamma}$ and hence on the event
$\gamma>3\eps m$ we have \be \label{boundmin} \min_{t \leq 3\eps m}
Y_t \geq \min _{t \leq 3\eps m} W_t \, .\ee

Note that the expectation of the increment of $W_t$ is
$-15\eps-15\eps^2$, thus for any positive $\alpha>0$ the process
$-\alpha W_t$ is a submartingale whence $\exp(-\alpha W_t)$ is a
submartingale as well. We put $\alpha = 8 \eps$ and apply Doob's
maximal $L^2$ inequality (see \cite{D}) yields that
$$ \E \Big [ \max_{t \leq 3\eps m} e^{-16 \eps W_t}\Big ] \leq  4 \E
\Big [ e^{-16\eps W_{3\eps m} } \Big ] \, .$$

Since $W_{3\eps m}$ is distributed as Bin$(3\eps(1-15\eps)m^2,p)-3\eps
m+1$ we obtain by direct computation that
$$ \E \Big [ \max_{t \leq 3\eps m} e^{-16 \eps W_t}\Big ] \leq
4 e^{672 \eps^3 m} \, .$$ Markov's inequality implies that
$$ \prob \Big ( \exists \, t \leq 3\eps m \with W_t \leq -45\eps^2 m \Big
) \leq \prob \Big ( \max_{t \leq 3\eps m} e^{-16 \eps W_t} \geq
e^{720\eps^3 m} \Big ) \leq 4 e^{-48\eps^3 m} \, .$$ Note that if
there exists $t \leq 3\eps m$ with $Y_t \leq -45\eps^2 m$ then by
(\ref{boundmin}) either $\gamma \leq 3\eps m$ or there exists $t\leq
3\eps m$ such that $W_t \leq -45\eps^2 m$. Lemma \ref{lowerbound}
shows that $\prob (\gamma \leq 3\eps m) \leq e^{-c\eps m} = o(
e^{-48\eps^3 m} )$ and this concludes the proof
of the lemma. \qed \\

\noindent We now use the estimates of the previous lemma to
amplify Lemma \ref{lowerbound}.

\begin{lemma} \label{sharplowerbound} For $p={1 +\eps \over m}$ there exists some fixed $c>0$ such that
have that
$$ \prob \Big ( \exists t \leq 3\eps m \with N_t \leq m - t - 50 \eps^2 m \Big ) \leq 9 e^{-c\eps^3 m}  \, .$$
\end{lemma}
\noindent{\bf Proof.} Let $\alpha_i$ be independent random variables
distributed as Bin$(m,p)$ and we couple such that $\eta_i \leq
\alpha_i$ for all $i$. By (\ref{nrec}) and the fact that $Z_t$ is
non-decreasing we have that for $t\leq 3\eps m$ \be \label{sharpstep1}
N_t \geq m - 1 - \sum_{i=1}^{t} \alpha_i - Z_{3\eps m} \, .\ee
Observe that if for some positive $k$ we have $Y_t \geq -k$ for all
$t\leq T$ then $Z_T \leq k$. Thus, Lemma \ref{ztbound} together with
the fact that $\{Z_t\}$ is increasing implies that
$$ \prob \Big ( Z_{3\eps m} \geq 45 \eps^2 m \Big ) \leq 5 e^{-48\eps^3 m} \,
.$$ We have that $\sum_{i=1}^{t} \alpha_i$ is distributed as
Bin$(mt, p)$ and has mean $t+\eps t$. The same argument using Doob's
maximal inequality, as in the proof of Lemma \ref{ztbound}, gives
that
$$ \prob \Big ( \exists t \leq 3\eps m \with \sum_{i=1}^{t} \alpha_i \geq t + 4\eps^2 m \Big)
\leq 4e^{-c \eps^3 m } \, ,$$ for some fixed $c>0$. The assertion of
the lemma follows by putting the last two inequalities into
(\ref{sharpstep1}). \qed \\

\begin{lemma} \label{ytlowerbound} Assume that $p={1 +\eps \over m}$ and that $\eps^3 m \geq 1$. Then there exist a constant $c>0$ such that
for any $a$ satisfying $1 \leq a \leq \sqrt{\eps^3 m}$ we have
$$ \prob \Big ( Y_t > 0 \,\, \hbox{ {\rm for
all} } a\sqrt{m/\eps} \leq t \leq 2\eps m - a \sqrt{m/\eps} \Big )
\geq 1 - 2 e^{-c a^2} \, .$$
\end{lemma}
\noindent {\bf Proof.} Denote by $\gamma$ the stopping time
$$ \gamma = \min \{ t \, : \, N_t < m-t-50\eps^2m  \} \, .$$ Lemma \ref{sharplowerbound} states that
$$\prob(\gamma \leq 3\eps m) \leq 9e^{-c\eps^3 m} \, ,$$
for some constant $c>0$. Let $\{W_t\}$ be a process with independent
increments distributed as Bin$(m-t-50\eps^2m,p)-1$ (note that the
increments are not identically distributed) and $W_0=1$. As usual we
can couple such that $Y_{t \wedge\gamma} \geq W_{t \wedge \gamma}$
for all $t$. Hence, if $\gamma \geq 2\eps m$ and there exists $t
\leq 2\eps m$ with $Y_t \leq 0$ then it must be that $W_t \leq 0$.
We conclude that it suffices to show the assertion of the
lemma to the process $\{W_t\}$ and this is our next goal.

For any $\alpha >0$ we have $$ \E \Big [ e^{-\alpha (W_t - W_{t-1})}
\mid W_{t-1} \Big ] = e^\alpha [ 1-p(1-e^{-\alpha})]^{m-t - 50
\eps^2 m} \, .$$ We use $1-x \leq e^{-x}$ with $x=p(1-e^{-\alpha})$
and $1-e^{-\alpha} \geq \alpha - \alpha^2$ for $\alpha$ small enough
(we will eventually take $\alpha=O(\eps)$) to get \be \label{wtstp0}
\E \Big [ e^{-\alpha (W_t - W_{t-1})} \mid W_{t-1} \Big ] \leq e^{
\alpha^2 ( 1+\eps) - \alpha (\eps - {t \over m}(1+\eps)
-50\eps^2(1+\eps)) } \, .\ee Thus, we learn that the process
$$ e^{-\alpha W_t} e^{ - (1+\eps) \alpha^2 t  - (1+\eps) \alpha {t^2 \over 2m} + \eps \alpha t  ( 1 - 50\eps(1+\eps))}  \, ,$$
is a supermartinagle. Write $$ f(t) = t \big [ -(1+\eps)\alpha^2 +
\eps \alpha(1-50\eps(1+\eps)) \big ] - t^2 {(1+\eps)\alpha \over 2m}
\, .$$ We apply the optional stopping theorem on the stopping time
$\tau = \min \{ t \geq \sqrt{m/\eps} : W_t = 0 \}$ and get that
$$ \E e^{f(\tau)}\leq 1 \, .$$
Direct calculation gives that when we put $\alpha = {1\over 3} \eps$
the function $f$ attains its minimum on the interval
$[a\sqrt{m/\eps}, \eps m]$ at $\tau = a\sqrt{m/\eps}$ for any
$a\in[1,\sqrt{\eps^3m}/3]$. Hence
$$ \prob \big ( a\sqrt{m/\eps} \leq \tau \leq \eps m \big ) \leq
\prob ( e^{f(\tau)} \geq e^{f(a\sqrt{m/\eps})} ) \, .$$ An immediate
calculation shows that $f(a\sqrt{m/\eps}) \geq ca\sqrt{m\eps^3}$ and
we learn by Markov's inequality that \be\label{tauest} \prob \big (
a\sqrt{m/\eps} \leq \tau \leq \eps m \big ) \leq e^{-ca\sqrt{m
\eps^3}} \leq e^{-ca^2} \, ,\ee
since $a \leq \sqrt{m \eps^3}$.

We are left to estimate $\prob (\eps m \leq \tau \leq 2\eps m -
a\sqrt{m/\eps})$. To that aim we define a new process $\{X_t\}_{t
\geq 0}$ by $X_t = W_{\eps m + t}$. By (\ref{wtstp0}), for positive
$\alpha$ we have that
\begin{eqnarray*} \E \Big [ e^{-\alpha (X_t - X_{t-1})} \, \mid \, X_{t-1} \Big ] \leq e^{
\alpha^2 ( 1+\eps) - \alpha (\eps - {t + \eps m \over m}(1+\eps)
-50\eps^2(1+\eps)) } \, .
\end{eqnarray*}
%
This together with a straight forward computation yields that the
process
$$ e^{-\alpha X_t} e^{ - \alpha^2 t(1+\eps) - \alpha \Big ( {(1+\eps)t^2 \over 2m}
+ 55 \eps^2 t \Big ) } \, ,$$ is a supermartingale.
Write $\tau$ for the stopping time
$$ \tau = \min \{ t\geq 0 \, :\, X_t =0 \}
\, .$$  Optional stopping yields that \begin{eqnarray}
\label{ytstp3} \E \Big [ e^{ - \alpha^2 \tau
(1+\eps) - \alpha \Big ( {\tau^2 \over 2m} +
55\eps^2 (\tau\wedge4\eps m) \Big ) } \Big ] \leq \E \Big [
e^{-\alpha X_0} \Big ] \leq e^{\alpha^2 \eps m (1+\eps) - \alpha
\Big ( {\eps^2 m \over 2} - 55\eps^3 m \Big ) } \, ,
\end{eqnarray}
where the last inequality is an immediate calculation with
(\ref{wtstp0}) and the fact that $X_0=W_{\eps m}$. Observe that the
exponent on the left hand side of the previous display is
$$ f(\tau) = -\alpha^2 \tau (1+\eps) - \alpha (\tau^2/2m + 55 \eps^2 \tau ) \,
,$$ which is a non-increasing function of $\tau$ on $[0,\infty)$.
Hence, for any $a\in [1,\sqrt{\eps^3m}]$ we get that
\be\label{ytstp4} \prob \big ( \tau \leq \eps m - a\sqrt{m/\eps}
\big ) \leq \prob \Big ( e^{ f(\tau) } \geq e^{f(\eps
m - a\sqrt{m/\eps})} \Big ) \, .\ee We have that
$$ f(\eps m - a\sqrt{m/\eps}) \geq  - 2\alpha^2 \eps m - \alpha \Big ( {\eps^2 m \over 2} - a \sqrt{\eps m} + {1 \over 2}a^2 \eps^{-1} +55 \eps^3 m \Big ) \, .$$
We use Markov inequality and (\ref{ytstp3}) to get
\begin{eqnarray*} \prob \Big ( e^{ f(\tau) } \geq e^{f(\eps m -
a\sqrt{m/\eps})} \Big ) &\leq& e^{4 \alpha^2 \eps m - \alpha \Big
( a \sqrt{\eps m} - {1 \over 2}a^2 \eps^{-1} - 110 \eps^3 m \Big ) } \\
&\leq& e^{4 \alpha^2 \eps m - c \alpha a \sqrt{\eps
m} } \, ,
\end{eqnarray*}
where in the last inequality we used our assumption on $a$ and
$\eps$. We choose $\alpha \approx a (\eps m)^{-1/2}$ that
minimizes the last expression. This yields
$$ \prob \Big ( e^{ f(\tau) } \geq e^{f(\eps m -
a\sqrt{m/\eps})} \Big ) \leq e^{-{c a^2}} \, .$$ We put this
into (\ref{ytstp4}), which together with (\ref{tauest}) yields the assertion of the lemma. \qed

\begin{lemma}\label{yepscond} Assume that $p={1 +\eps \over m}$. Write $\tau = \min \{t : Y_t = 0\}$, then for any small $\alpha>0$, we have
$$ \E \big [ e^{\alpha Y_{\eps^{-2}}} \mid \tau \geq \eps^{-2} \big ] \leq C e^{2\alpha\eps^{-1}+\alpha^2 \eps^{-2}} \, .$$
\end{lemma}
\noindent{\bf Proof.} We have that $\prob(\tau \geq \eps^{-2}) \geq c\eps$. To see this we perform the usual argument of bounding $Y_t$ below by a process of independent increments (until a stopping time, using Lemma \ref{lowerbound}) and using Lemma \ref{tail}. This has been done in this section several times so we omit the details. Thus, it suffices to bound from above $\E e^{\alpha Y_{\eps^{-2}}}{\bf 1}_{\{\tau \geq \eps^{-2}\}}$. Since we can bound $Y_t$ by a process $W_t$ which has i.i.d. Bin$(m,p)-1$ increments, it suffices to bound the same expectation for $W_t$. Write $\gamma = \min \{ t : W_t = 0 \hbox{ {\rm or} } W_t \geq \eps^{-1}\}$. We have
$$ \E e^{\alpha W_{\eps^{-2}}}{\bf 1}_{\{\tau \geq \eps^{-2}\}}\leq \E e^{\alpha W_{\eps^{-2}}}{\bf 1}_{\{ \tau \geq \eps^{-2}, \gamma \geq \eps^{-2} \} } + \E e^{\alpha W_{\eps^{-2}}}{\bf 1}_{\{ \tau \geq \eps^{-2}, \gamma < \eps^{-2} \} } \, .$$
For the first term on the right hand side we note that on $\gamma \geq \eps^{-2}$ we have that $W_{\eps^{-2}} \leq \eps^{-1}$, so
$$ \E e^{\alpha W_{\eps^{-2}}}{\bf 1}_{\{ \tau \geq \eps^{-2}, \gamma \geq \eps^{-2} \} } \leq C \eps e^{\alpha \eps^{-1}} \, .$$
For the second term we condition on $\{\tau\geq \eps^{-2}, \gamma <
\eps^{-2}\}$ (which implies $W_\gamma \geq \eps^{-1}$ and
$\gamma<\eps^2$) to get that
\be\label{yepscond.midstep} \E e^{\alpha W_{\eps^{-2}}}{\bf 1}_{\{ \tau \geq \eps^{-2}, \gamma < \eps^{-2} \} } \leq \prob(W_\gamma \geq \eps^{-1}) \E [ e^{\alpha W_\gamma}e^{\alpha (W_{\eps^{-2}} - W_\gamma)}  \mid W_\gamma \geq \eps^{-1},\gamma<\eps^{-2} ] \, .\ee
We have that $\prob(W_\gamma \geq \eps^{-1}) = O(\eps)$ by Lemma 7
of \cite{NP2}. We condition in addition on $W_\gamma$ and $\gamma$ and pull out the $e^{\alpha W_\gamma}$ factor. By the strong Markov property we have that conditioned on all these, the random variable $W_{\eps^{-2}} - W_\gamma$ is distributed as the sum of $\eps^{-2}-\gamma$ i.i.d. copies of Bin$(m,p)-1$ random variables. Thus,
$$ \E [ e^{\alpha (W_{\eps^{-2}} - W_\gamma)} \mid  W_\gamma , \gamma<\eps^2 ] \leq e^{-\alpha \eps^{-2}} [1+p(e^{\alpha}-1)]^{m\eps^{-2}} \, .$$
Furthermore, Lemma 5 of \cite{NP1} states that conditioned on
$W_\gamma \geq \eps^{-1}$ and $\gamma<\eps^2$ the distribution of
$W_\gamma - \eps^{-1}$ is bounded above by Bin$(m,p)$, whence
$$ \E [ e^{\alpha W_\gamma} \mid W_\gamma \geq \eps^{-1} , \gamma < \eps^{-2}] \leq e^{\alpha \eps^{-1}} [1+p(e^{\alpha}-1)]^m \, .$$
Putting this back into (\ref{yepscond.midstep}) gives
$$ \E e^{\alpha W_{\eps^{-2}}}{\bf 1}_{\{ \tau \geq \eps^{-2}, \gamma < \eps^{-2} \} } \leq C\eps e^{-\alpha(\eps^{-2} - \eps^{-1})} [1+p(e^{\alpha}-1)]^{m(\eps^{-2} + 1)} \, .$$
Putting all these together we get
\begin{eqnarray*} \E \big [ e^{\alpha Y_{\eps^{-2}}} \mid \tau \geq \eps^{-2} \big ] &\leq& C e^{\alpha \eps^{-1}} + C e^{-\alpha(\eps^{-2} - \eps^{-1})} [1+p(e^{\alpha}-1)]^{m(\eps^{-2}+1)} \\ &\leq& C e^{\alpha \eps^{-1}} + C e^{-\alpha(\eps^{-2} - \eps^{-1})} e^{(1+\eps)(\alpha+\alpha^2)(\eps^{-2}+1)} \, ,
\end{eqnarray*} The lemma follows now by an immediate calculation. \qed \\

\begin{lemma} \label{cvbound} Let $p={1 +\eps \over m}$ and assume $\eps^3 m\geq 1$. Then for any $\ell>0$, we have
$$ \prob( |\C(v)| \geq 2\eps m + \ell) \leq  C\eps e^{ {-c\ell^2 (2\eps m + \ell) \over m^2}} \, .$$
\end{lemma}
\noindent{\bf Proof.}
We assume that $\ell \geq 2\sqrt{m/\eps}$ since otherwise the exponential is of constant order and the assertion of the lemma follows simply from Lemma \ref{tail}. Recall that $|\C(v)|$ is distributed as the first hitting time $\tau$ of $Y_t$ at $0$.
We put $T=2\eps m +\ell$ and condition on $Y_{\eps^{-2}}$ and on
$\tau \geq \eps^{-2}$. That is, \be\label{cvbound.firststep}
\prob(\tau \geq 2\eps m + \ell) = \prob(\tau \geq \eps^{-2}) \E \big
[\prob(\tau \geq T \mid Y_{\eps^{-2}}, \tau \geq \eps^{-2} ) \big ]
\, .\ee Since $Y_t$ is bounded above by a process with increments
distributed as Bin$(m,p)-1$, we learn by Lemma \ref{tail} that
$\prob(\tau \geq \eps^{-2}) = O(\eps)$. The second term will give us the exponential in the assertion of the Lemma simply because $Y_T$ has small probability of being positive at this time. Indeed, since the increments of $Y_t$ are
stochastically bounded above by Bin$(m-t,p)-1$ we have that for any
small $\alpha>0$ $$ \E \Big [ e^{\alpha (Y_t-Y_{t-1})} \mid Y_{t-1}
\Big ] \leq e^{-\alpha} [ 1+p(e^{\alpha}-1)]^{m-t} \leq  e^{-\alpha
+ (1+\eps)(\alpha + \alpha^2)(1-t/m)} \, ,$$ since $e^\alpha -1 \leq
\alpha + \alpha^2$ for small enough $\alpha$. Summing this over $t$ ranging from $\eps^{-2}$ to $T$ gives
\begin{eqnarray*} \E \big [ e^{\alpha Y_T} \mid Y_{\eps^{-2}}, \tau \geq \eps^{-2} \big ] &\leq&  e^{-\alpha (T-\eps^{-2}) +
(1+\eps)(\alpha + \alpha^2)\big (T- \eps^{-2} - {T^2 -\eps^{-4} \over 2m} \big )} e^{\alpha Y_{\eps^{-2}}}
\\  &\leq& e^{\alpha^2 T(1+\eps) - \alpha \Big [ {T^2 - \eps^{-4} \over 2m }
- \eps T \Big ]} e^{\alpha Y_{\eps^{-2}}} \, .\end{eqnarray*}
Hence,
\begin{eqnarray*} \E \big [ e^{\alpha Y_T} \mid  \tau \geq \eps^{-2} \big ] &\leq& \E [ e^{\alpha Y_{\eps^{-2}}} \mid \tau \geq \eps^{-2} ] e^{\alpha^2 T(1+\eps) - \alpha \Big [ {T^2 - \eps^{-4} \over 2m }
- \eps T \Big ]}  \\ &\leq& C e^{\alpha^2 (T+\eps^{-2})(1+\eps) - \alpha \Big [ {T^2 - \eps^{-4} \over 2m }
- \eps T - 2\eps^{-1} \Big ]} \, ,\end{eqnarray*}
where the last inequality is due to Lemma \ref{yepscond}. Hence, by
Markov's inequality this is also an upper bound on $\prob(Y_T \geq 0 \mid  \tau \geq \eps^{-2})$ which is what we aim to estimate. We now choose $\alpha$
$$ \alpha = {{T^2 -\eps^{-4} \over 2m } - \eps T-2\eps^{-1}  \over 2(T+\eps^{-2})} \,
,$$ which is positive and of order $\ell /m$ since $\ell \geq 2\sqrt{m/\eps}$ and minimizes the above expectation. We get that
$$ \prob(Y_T \geq 0 \mid  \tau \geq \eps^{-2}) \leq C e^{-{c T(T-2\eps m)^2 \over m^2}} \, ,$$
for some $c>0$ by a straightforward calculation, concluding our proof. \qed \\

\subsubsection{Proof of near-critical random graph theorems.}
We are now ready to prove the Theorems stated in Section \ref{nearcritsec}. \\

\noindent{\bf Proof of Theorem \ref{supergnp1}.} We begin by proving the upper bound on $\E|\C_1|$. For any positive integer $\ell$ define by $X_\ell$ the random
variable
$$ X_\ell = \Big | \Big \{ v \, : \, |\C(v)| \geq \ell \Big \} \Big | \, .$$
Observe that if $|\C_1| \geq \ell$, then we must have that $|X_\ell|
\geq |\C_1|$. Thus for any positive integer $\ell$ we have \be
\label{upbound} \E |\C_1| \leq \ell\,  \prob (|\C_1| < \ell) + \E X_\ell \, .\ee

We take $\ell = {1 \over 20} \eps m$ and since Lemma
\ref{ytlowerbound} implies that $\prob( |\C_1| \leq \ell ) \leq
Ce^{-c\eps^3 m}$ and $\eps^3 m \geq A\log m$ we have that the first term on the right hand side
of (\ref{upbound}) is $o(1)$.
We now turn to bound the second term on the right hand side of
(\ref{upbound}). Since $\E X_\ell = m \prob ( |\C(v)| \geq \ell )$
it suffices to bound from above $\prob ( |\C(v)| \geq \ell )$.
Recall that $|\C(v)|$ is the hitting time of the process $\{Y_t\}$
at $0$. Let $\{W _t \}$ be a process with independent increments
distributed as Bin$(m,p)-1$ and $W_0=1$, as in Lemma \ref{tail}. Let
$\tau = \min_t \{W_t =0 \}$ be the hitting time of $W_t$ at $0$,
then it is clear that we can couple $W_t$ and $Y_t$ such that
$|\C(v)| \leq \tau$. Thus

$$\prob ( |\C(v)| \geq \ell ) \leq \prob ( \tau \geq \ell ) = \prob ( \tau = \infty ) + \prob ( \ell \leq \tau < \infty ) \, .$$
We now apply Lemma \ref{tail} with $T=\ell={1 \over 20} \eps m$ and
get by the previous display that \begin{eqnarray*} \prob (
|\C(v)| \geq \eps m /20 ) &\leq& 2\eps - {8
 \over 3} \eps^2 + O(\eps^3) + C_1 \eps^{-7/2}m^{-3/2} e^{-\eps^3 m /4}  \\&=& 2\eps - {8  \over 3} \eps^2 + O(\eps^3) \, ,\end{eqnarray*}
as long as $\eps^3 m \geq A \log m$ for large enough $A$. We
conclude that
$$ \E X_\ell \leq 2\eps m - {8 \over 3} \eps^2 m +
O(\eps^3 m) \, ,$$ which together with (\ref{upbound}) concludes the
proof of the upper bound on $\E |\C_1|$.

We turn to the proof of the lower bound on $\E|\C_1|$.  Recall that at
each record minimum of the process $\{Y_t\}$ we are starting the
exploration of a new component. Write
$$ \gamma = \max \Big \{ t \leq \eps m \, : \, Y_t \hbox{ {\rm is
at a record minimum}} \Big \} \, ,$$ and
$$ \tau = \min \{ t \geq 0 : Y_{\eps m+ t} < 0 \}  \, .$$ Then we have
that \be \label{thelowerbound} |\C_1| \geq \eps m - \gamma + \tau \,
.\ee Thus, in order to complete the proof we will provide an upper
bound on $\E \gamma$ and a lower bound on $\E \tau$. Let $\{W_t\}$
be a process defined as in Lemma \ref{recmin} with i.i.d. increments
distributed as Bin$(m(1-\eps/2),p)-1$. Define the stopping time
$\beta$ by
$$ \beta = \min \{ t \, : \, N_t \leq m(1-\eps/2) \} \, ,$$
then it is clear we can couple $\{Y_{t \wedge \beta}\}$ with
$\{W_{t \wedge \beta}\}$ such that the increments of the first are
larger than of the latter process. This guarantees that every
record minimum of the first process is also a record minimum of
the second, and thus if we put
$$ \gamma^w = \max \Big \{ t \, : \, W_t \hbox{ {\rm is
at a record minimum}} \Big \} \, ,$$ then we get that we can
couple such that
$$  \gamma {\bf 1}_{\{\hbox{{\rm no record minima at times} } [\eps m/10,\eps m] \}} \leq \gamma^w + \eps m {\bf 1}_{\{ \beta \leq \eps m/10 \}} \, .$$
Lemma \ref{ytlowerbound} shows that the probability that there is a
record minimum at some time between $\eps m/10$ and $\eps m$ decays
faster than $m^{-2}$ provided that $\eps^3 m \geq A\log m$ for $A$ large enough. Hence, taking expectations on both sides
and using Lemma \ref{recmin} and Lemma \ref{lowerbound} gives that $\E \gamma
= O(\eps^{-2})$.

We now turn to give a lower bound on $\E \tau$. We begin by
estimating $\E \tau^2$. As before, define the process $\{X_t\}_{t
\geq 0}$ by $X_t = Y_{\eps m + t}$ and note that $X_t - X_{t-1} =
\eta_{\eps m + t}$. For any $t$ such that $N_{t+\eps m} \geq
m-(t+\eps m)-50\eps^2 m$ we have
$$ \E [ X_{t+1} - X_{t} \, \mid \, \F_t ] \geq -{t(1+\eps) \over m} -55 \eps^2  \, .$$
Thus the process $\{ X_{t\wedge T} + {(t\wedge T)^2(1+\eps) \over
2m} + 55\eps^2 (t \wedge T)\}$ is a submartingale, where $T$ is
defined as
$$ T = \min\{t-\eps m:t\geq \eps m, N_t \leq m-t-50\eps^2 m\} \, .$$
Optional stopping yields that
\be\label{supergnp1.midstep} \E (\tau\wedge T)^2 \geq {2m\over 1+\eps} (\E X_0 - \E X_{\tau \wedge
T}) - {110\eps^2 m\over1+\eps} \E [\tau \wedge T] \, .\ee By Lemma
\ref{ytlowerbound}, we have $$\P(\tau<\eps m-a\sqrt{m/\eps})\leq
e^{-ca^2}.$$ Also by lemma \ref{ytlowerbound}, one can deduce
\begin{eqnarray*}
\P(\tau>\eps m+a\sqrt{m/\eps})&\leq&
\P(\tau>\eps m+a\sqrt{m/\eps}, Y_t>0
\hbox{{\rm \,\,for\,\, }}t\in[{a\over2}\sqrt{m/\eps},\eps
m]) + e^{-ca^2}\\&\leq& \P(|\C_1|>2\eps
m+{a\over2}\sqrt{m/\eps}) + e^{-ca^2} \\ &\leq& \P(X_{2\eps m a\sqrt{m/\eps}/2} > 2\eps
m) + e^{-ca^2},
\end{eqnarray*} where $X_{2\eps m a\sqrt{m/\eps}/2}$ is the number of vertices $v$ such that $|\C_v|\geq 2\eps m+{a\over2}\sqrt{m/\eps}$ as defined in the beginning of the proof. By Lemma \ref{cvbound}, we have $$\E X_{2\eps m a\sqrt{m/\eps}/2} \leq C m\eps e^{-ca^2} \, .$$ Plugging this into the previous inequality and using
Markov's inequality shows that $\E [\tau \wedge T] = O(\eps m)$
and \be\label{tauld}\P(|\tau-\eps m|>a\sqrt{m/\eps})\leq
C e^{-ca^2}.\ee Lemma \ref{sharplowerbound} shows that $\prob(T \leq 2\eps m) \leq m^{-2}$, and so
$\E X_{\tau \wedge T} = o(1)$ and $\E (\tau \wedge T)^2 = \E \tau^2 + o(1)$. We get that
$$ \E \tau^2 \geq 2m \E Y_{\eps m} - o(1) \, .$$
We bound from below $\E Y_{\eps m}$ using the approximating process $\widetilde{Y}_t$ defined in (\ref{approxprocess}). We have that $\E \widetilde{Y}_{\eps m} = \eps^2 m^2 + O(\eps^3 m)$ and using (\ref{quoteBR}) and Lemma \ref{numberofrecmins} we deduce the same estimate for $\E Y_{\eps m}$. This yields that
$$ \E \tau^2 \geq \eps^2 m^2 - C\eps^3 m^2 \, ,$$ for some $C>0$. Inequality (\ref{tauld}) gives that for some
$C>0$ we have
$$ {\rm Var}(\tau) \leq \E \Big [ (\tau - \eps m)^2 \Big ] \leq  {Cm \over
\eps} \, .$$ We conclude
$$ \E \tau = \sqrt{ \E \tau^2 - {\rm Var}(\tau)} \geq \eps m
\sqrt{ 1- C\eps - {C \over \eps^3 m} } \geq \eps m - C\eps^2 m -
C\eps^{-2} \, ,$$ since $\sqrt{1-x} \geq 1-x$ for $x \in (0,1)$.
Using this and our estimate on $\E \gamma$ in
(\ref{thelowerbound}) finishes the proof. \qed \\

%

%
%


\noindent {\bf Proof of Theorem \ref{gnpout2}.} Since component
sizes are excursions' length above past minima and $Y_0=1$, Lemma
\ref{ytlowerbound} immediately yields the bound \be\label{out2first}
\prob \Big ( |\C_1| \leq 2\eps m - A\sqrt{m\over \eps} \Big ) \leq
e^{-c A^2} \, ,\ee valid for any $A$ satisfying $1 \leq A \leq
\sqrt{\eps^3 m}$. For the upper bound we use Lemma \ref{cvbound} stating that
$$ \prob(|\C(v)| \geq 2\eps m + A \sqrt{m/\eps}) = O(\eps e^{-cA^2}) \, .$$
Write $X=|\{v : |\C(v)|\geq  2\eps m + A \sqrt{m/\eps}\}|$ so that $\E X = O( \eps m e^{-cA^2})$. As usual we have
$$ \prob(|\C_1| \geq 2\eps m + A \sqrt{m/\eps}) \leq \prob(X \geq 2\eps m) = O(e^{-cA^2}) \, ,$$
by Markov's inequality, concluding the proof. \qed \\

\noindent {\bf Proof of Corollary \ref{supergnpmoments}.}
Part (i) of the corollary follows immediately from Theorem \ref{gnpout2} by integration. Indeed,
\begin{eqnarray*}
&&\E \Big [ \Big | |\C_1| - 2\eps m \Big |^k \Big ]=\sum_\ell \ell^{k-1}\P(|\C_1-2\eps m|>\ell)\\
&\leq& \sum_{\ell=1}^{\sqrt{\frac{m}{\eps}}} \ell^{k-1}+C \sum_{\ell=\sqrt{\frac{m}{\eps}}}^{\eps m} \ell^{k-1} e^{{-c\ell^2 \eps \over m}} +\sum_{\ell \geq \eps m} \ell^{k-1} e^{{-c\ell^3 \over m^2}} \, ,
\end{eqnarray*} where we bounded the second sum on the right hand side using Theorem \ref{gnpout2} and the last sum using Lemma \ref{cvbound} (which is valid for all $\ell > 0$ and not limited by $\ell \leq \sqrt{\eps^3 m}$) and the usual Markov inequality on the variable $X=|\{v : |\C(v)|\geq  2\eps m + \ell\}|$. A quick calculation now shows each term is of order at most $(m/\eps)^{k/2}$, concluding our proof.
\qed \\

\noindent{\bf Proof of Theorem \ref{subgnpmoments}.} We begin by proving (i). As before, $|\C(v)|$ is stochastically dominated by the random variable $\tau$ defined in Lemma \ref{tail2}. This Lemma gives that for any fixed $k \geq 1$
$$ \E |\C(v)|^k = O ( \eps^{-2k+1} ) \, .$$
Number the vertices of $G(n,p)$ arbitrarily $v_1, \ldots, v_m$ and
observe that
$$ \sum _{j \geq 1} |\C_j|^k = \sum _{i=1}^m |\C(v_i)|^{k-1} \,
,$$ because each component $\C_j$ is counted in the sum in the
right hand size precisely $|\C_j|$ times. By symmetry we learn
that
$$ \E \sum _{j \geq 1} |\C_j|^k = m \E |\C(v)|^{k-1} = O( m \eps^{-2k+3}
) \, ,$$ finishing the first assertion of the theorem.

We proceed to prove (ii). Recall that $\sum_j|\mathcal{C}_j|^2=\sum_v|\C(v)|$. Thus,
\begin{eqnarray*}
\E\big(\sum_j|\mathcal{C}_j|^2\big)^2= \E\sum_{v,w}|\C(v)||\C(w)|1_{\{\C(v) = \C(w)\}}+\E\sum_{v,w}|\C(v)||\C(w)|1_{\{\C(v) \neq \C(w) \}} \, .
\end{eqnarray*}
The first term on the right hand side is simply $\sum_v \E |\C(v)|^2$ which equals $\E \sum_j |\C_j|^3$ and is upper bounded by $O(m \eps^{-3})$ by part (i) of the theorem. For the second term we note that we can write $\E\sum_{v,w}|\C(v)||\C(w)|1_{\{\C(v) \neq \C(w)\}}$ as
$$\E\sum_w|\C(w)|\sum_v|\C(v)|1_{\{v \not \in \C(w)\}}=\E\sum_w|\C(w)|\sum_{\{v \not \in \C(w) \}}|\C(v)| \, .$$
Conditioned on $\C(w)$ the distribution of the rest of the graph is also subcritical random graph with $\eps'$ bigger than $\eps$. Thus the estimate of part (i) of the theorem (together with the fact that $\sum_v |C(v)| = \sum_j |\C_j|^2$) can be applied and we may bound
$$ \E\sum_{v,w}|\C(v)||\C(w)|1_{\{\C(v) \neq \C(w)\}} \leq C m\eps^{-1} \E \sum_w |\C(w)| = O(m^2 \eps^{-2}) \, ,$$
which finishes the proof of (ii).

To prove part (iii) of the theorem, let $W_t$ be a process with
i.i.d. increment distributed as
Bin$(m-5\eps^{-2},\frac{1-\eps}{m})-1$ and $W_0=1$. Let
$$\tau=\min\{t:N_t<m-5\eps^{-2}\}.$$ As usual we can couple such that
$Y_{t\wedge\tau}\geq W_{t\wedge\tau}$. Let $\gamma=\min\{t:
W_t\leq0\}.$ For any $T$ We have
\begin{eqnarray*}
\P(\gamma\geq T)&=&\P(\gamma\geq T, \tau\leq T)+\P(\gamma\geq T,
\tau> T)\\&\leq&\P(\tau\leq T)+\P(|\C(v)|\geq T),
\end{eqnarray*} which implies \be\label{subgnpmoments.midstep} \P(|\C(v)|\geq T)\geq\P(\gamma\geq
T)-\P(\tau\leq T) \, .\ee Put $T=\eps^{-2}$ we have by Lemma
\ref{lowerbound} that $$\P(\tau\leq T)\leq e^{-c\eps^{-2}} \, .$$
Furthermore, Lemma \ref{tail2} shows that $$\P(\gamma\geq \eps^{-2})\geq
c\eps\, ,$$ for some constant $c>0$. Thus, by (\ref{subgnpmoments.midstep}) we get that
$$\P(|\C(v)|\geq \eps^{-2})\geq c\eps-e^{-c\eps^{-2}} \, ,$$ which
implies $\E|\C(v)|\geq c\eps^{-1}$ and concludes the proof. \qed \\

\noindent {\bf Proof of Theorem \ref{gnp2moments}.}
%
%
The proof of (i) is a calculation using Lemma \ref{cvbound}. We have
$$ \E |\C(v)|^k = \sum_{\ell=1}^{\eps^{-2}} \ell^{k-1} \prob (|\C(v)|\geq \ell) + \sum_{\ell=\eps^{-2}}^{10\eps m} \ell^{k-1} \prob (|\C(v)|\geq \ell) + \sum_{\ell=10\eps m}^{m} \ell^{k-1} \prob(|\C(v)|\geq k) \, .$$
For the first sum we use the estimate $\prob(|\C(v)|\geq l) \leq
O(\eps + \ell^{-1/2})$ appearing in the proof of Proposition 1 of
\cite{NP2}. We get
$$ \sum_{\ell=1}^{\eps^{-2}} \ell^{k-1} \prob (|\C(v)|\geq \ell) \leq C \sum_{\ell=1}^{\eps^{-2}} \ell^{k-1} (\eps+\ell^{-1/2}) = O(\eps^{-2k+1}) \, .$$
For the second sum, since $Y_t$ is bounded above by a process with
i.i.d. increments Bin$(m,p)$-1, each term is of order $\eps$ by
Lemma \ref{tail}. This gives the main contribution of
$O(\eps^{k+1} m^k)$.  Lastly, the third sum we bound using Lemma
\ref{cvbound} to get
$$ \sum_{\ell=10\eps m}^{m} \ell^{k-1} \prob(|\C(v)|\geq k) \leq C\eps \sum_{\ell=10\eps m}^{m} \ell^{k-1} e^{-c m^{-2} \ell^3} \, .$$
Since $\eps^3 m \geq 1$ we may bound the sum above by summing from $m^{2/3}$ to $m$. A straightforward calculation then gives that
$$ \sum_{\ell=10\eps m}^{m} \ell^{k-1} \prob(|\C(v)|\geq k) \leq C \eps (\eps m)^{k-1} m^{2/3} = O(\eps^{k+1} m^k)\, , $$
which finishes the proof of (i). We proceed to prove (ii). We have that \be\label{splittwo}\E \sum _{j \geq 2} |\C_j|^k=\E
\sum _{j \geq 2} |\C_j|^k\1_{\{|\C_1|< 1.5\eps m\}}+ \E \sum _{j
\geq 2} |\C_j|^k\1_{\{|\C_1|\geq 1.5\eps m\}}\, .\ee For the first
term of (\ref{splittwo}) we apply FKG inequality to get $$\E \sum
_{j \geq 2} |\C_j|^k\1_{\{|\C_1|< 1.5\eps m\}}\leq \E \sum _{j \geq
1} |\C_j|^k\1_{\{|\C_1|< 1.5\eps m\}}\leq \P(|\C_1|< 1.5\eps
m)\E\sum _{j \geq 1} |\C_j|^k.$$ By Theorem \ref{gnpout2}, we have
$$\P(|\C_1|< 1.5\eps m)\leq C e^{-c\eps^3 m}\, ,$$
and so
$$ \E \sum _{j \geq 2} |\C_j|^k\1_{\{|\C_1|< 1.5\eps m\}} \leq Ce^{-c\eps^3 m} m \E|\C(v)|^{k-1} \, .$$
By part (i) of the theorem this is at most $C\eps^k m^k e^{-c\eps^3 m}$ which is $O(m
\eps^{-2k+3})$ since $\eps \geq m^{-1/3}$. This shows the required
bound for the first term of (\ref{splittwo}).

To take care of the second term of (\ref{splittwo}) we condition on
$\C_1$ and note that the graph remaining is distributed as
$G(m-|\C_1|,p)$ conditioned on the event of not having a component
larger than $|\C_1|$. But since $|\C_1| \geq 1.5\eps m$ this random
graph is in the subcritical regime, and the probability of having
such a component is smaller than $1/2$ (in fact, it is exponentially
small). The required estimate follows by part (i) of Theorem
\ref{subgnpmoments}. This finishes the proof of (ii).

The proof of (iii) goes in similar lines of (ii). We have
$$ \E (\sum_{j\geq 2}|\C_j|^2)^2 = \E (\sum_{j\geq 2}|\C_j|^2)^2 {\bf 1}_{\{|\C_1|< 1.5\eps m \}} + \E (\sum_{j\geq 2}|\C_j|^2)^2 {\bf 1}_{\{|\C_1|\geq 1.5\eps m \}} \, .$$
As in the proof of (ii), to control the first term we use FKG inequality, extract $\prob(|\C_1|< 1.5\eps m)$ and bound the rest by $\E (\sum_{j \geq 1} |\C_j|^2)^2$ (instead of $j \geq 2$). The analysis performed in the proof of part (ii) of Theorem \ref{subgnpmoments} shows that $\E (\sum_{j \geq 1} |\C_j|^2)^2$ is controlled by $(\E \sum_{j \geq 1} |\C_j|^2)^2$. We get that
$$ \E (\sum_{j\geq 2}|\C_j|^2)^2 {\bf 1}_{\{|\C_1|< 1.5\eps m \}} \leq Ce^{-c\eps^3 m} m^4 \eps^4 = O(m^2 \eps^{-2}) \, .$$
To control the second term, as in the proof of (ii), we condition on $\C_1$ and use part (ii) of Theorem \ref{subgnpmoments} to estimate the remaining subcritical graph. This is done identically to part (ii) and we omit the details. \qed \\


\noindent {\bf Proof of Theorem \ref{gnpsmalls}.} We will use a second moment method. First we show that
$$ \E\sum _{|\C_j |\leq\delta\sqrt{m}} |\C_j|^2 \geq c m^{5/4} \, ,$$
for some $c=c(\delta)>0$. Indeed, we have \begin{equation}\label{niub}\E |\C(v)| {\bf
1}_{\{|\C(v)|\leq\delta\sqrt{m}\}}\geq \E |\C(v)| {\bf
1}_{\{{\delta\over2}\sqrt{m}\leq|\C(v)|\leq\delta\sqrt{m}\}}\geq{\delta\over2}\sqrt{m}\P({\delta\over2}\sqrt{m}\leq|\C(v)|\leq\delta\sqrt{m}).
\end{equation}
We proceed further by restricting to the case that $\C_v$ is tree.
Indeed, we have
\begin{eqnarray*}\P({\delta\over2}\sqrt{m}\leq|\C(v)|\leq\delta\sqrt{m})&\geq&\sum_{k=\delta/2
\sqrt{m}}^{\delta \sqrt{m}} \P(|\C(v)|=k, \C(v) \hbox{ {\rm is a tree}})\\&=&\sum_{k=\delta/2 \sqrt{m}}^{\delta \sqrt{m}} {m-1 \choose
k-1} k^{k-2} p^{k-1} (1-p)^{k(m-k) + {k \choose 2} - (k-1)}.
\end{eqnarray*} A quick calculation using Stirling's formula gives that for all such
$k$, each summand is of order $\Theta(m^{-4/3})$ and so the probability is of order at least $m^{-1/4}$ and the expectation in (\ref{niub}) is of order at least $m^{1/4}$. This gives the first moment estimate since
\begin{eqnarray*}
\E\sum _{|\C_j |\leq\delta\sqrt{m}} |\C_j|^2=\E\sum_{v:|\C(v) |\leq\delta\sqrt{m}}|C(v)| \, .
\end{eqnarray*}
We continue with the second moment estimate. By Theorem \ref{subgnpmoments} the second moment satisfies
$$ \E \big [ \sum_{j} |\C_j|^2 ]^2 = O(m^{5/2}) \, ,$$
and so the assertion of the Theorem follows by the inequality (see \cite{D})
$$ \prob(V > a) \geq {(\E V - a)^2 \over \E V^2} \, ,$$
valid for any non-negative random variable $V$ and $a < \E V$.  \qed \\

Now we turn to the proof of Theorem \ref{critgnp}. Recall that $Z_t$ counts the number of record minima of $\{Y_s\}$ before time $t$.

\begin{lemma} \label{numberofrecmins}  For any fixed $\delta \in (0,1/10)$, there exists
an universal constant $C>0$ such that as long as $\eps^3 m \geq 1$
we have
$$ \E Z_{\delta \eps m} \leq  {1 \over 2[(1-5\delta-5\delta \eps)\eps]} + O(1) \, ,$$
and
$$ \E Z^2_{\delta \eps m} = O(\eps^{-2}) \, .$$
\end{lemma}

\noindent {\bf Proof.} Define the stopping time $\tau$ by
$$ \tau = \min \{ t : N_t \leq m(1- 5 \delta \eps) \} \, ,$$
and $\{W_t\}$ to be the process with increments distributed as
Bin$(m(1- 5 \delta \eps), p)$ and $W_0=1$. As usual we can couple
such that $Y_{t \wedge \tau} \geq W_{t \wedge \tau}$ and that the
increments of the first process are always larger than of the
second. This guarantees that the number of record minimum of $Y_{t
\wedge \tau}$ is bounded from above by the record minimum of $W_{t
\wedge \tau}$. Denote by $Z^w$ the number of record minima of the
process $\{W_t\}$, then by the above discussion we have
$$ \E Z_{\delta \eps m} \leq \delta\eps m \E {\bf 1}_{\{ \tau < \delta \eps m \}} + \E Z^w \, .$$
The order of the first term can be arbitrarily small since
$\P(\tau<\delta\eps m)$ is exponentially small in $\eps m$ by Lemma
\ref{lowerbound}. Lemma \ref{numrecmin} bound the
second term by the required amount. This concludes the bound on $\E Z_{\delta \eps m}$. For
the second moment estimate, note that by the same argument, we have
$$ \E Z^2_{\delta \eps m} \leq \delta^2\eps^2 m^2 \E {\bf 1}_{\{ \tau <
\delta \eps m \}} + \E (Z^w)^2 \, ,$$
and the exponential decay of $\P(\tau<\delta\eps m)$ and Lemma \ref{numrecmin} concludes the proof. \qed \\

\begin{lemma} \label{hitrecmin} For any fixed $\delta \in
(0,1/10)$ denote by $\tau_\delta$ the stopping time
$$ \tau_\delta = \min _{t \geq \delta \eps m} \Big \{ t \hbox{
{\rm is a record minimum of} } Y_t \Big \} -\delta \eps m \, .$$
Then
$$ \E \tau_\delta \leq (2-\delta)\eps m - {1 \over 4\eps^{2}} \, .$$
\end{lemma}
\noindent {\bf Proof.} Define the process $\{X_t\}$ by $X_t =
Y_{\delta \eps m + t}$ so that
$$ \tau_{\delta} = \min \{ t \geq 0 \, : \, X_t = -Z_{\delta \eps m} \}
\, .$$ Let $\{W_t\}$ be a process defined by $W_0 = X_0$ and with
independent increments distributed as Bin$(m-t-\delta \eps m,
p)-1$ and let $\tau$ denote the stopping time $\displaystyle
\min_t~\{W_t~=~-Z_{\delta \eps m} \}$.  As usual, $X_t$ can be
stochastically bounded above by $W_t$ and hence $\E \tau_\delta
\leq \E \tau$ and we are left to estimate $\E \tau$. We have \be
\label{wdrift} \E [ W_t - W_{t-1} \, \mid \, \F_{t-1} ] =
(1-\delta)\eps - {t(1+\eps) \over m} - \delta \eps^2 \, .\ee Put
\begin{eqnarray*} f(t) &=& {t ^2 \over 2m} - (1-\delta)\eps t - (\delta -
\delta^2/2)\eps^2 m - \delta \eps^2 t +{t(1+\eps) + \eps t^2 \over 2m} \\
&=& {[t-(2-\delta)\eps m]^2 \over 2m} + \eps [t - (2-\delta)\eps m]
- \delta \eps^2 t +{t(1+\eps) + \eps t^2 \over 2m} \,
,\end{eqnarray*} then by (\ref{wdrift}) we deduce that $M_t = W_t
+ f(t)$ is a martingale. A direct calculation with (\ref{ytupper}) gives that
$$ \E W_0 = \E Y_{\delta \eps m} \leq -\delta \eps m + \delta \eps^2 m - \delta^2 \eps^2 m /2 + O(\eps^3 m) \, ,$$
and so we deduce that $\E M_0 \leq  C\eps^3 m$. Furthermore, we have that $\E \tau = O(\eps m)$ since after time $2 \eps m$ the process becomes subcritical with drift $-\eps$. Put $\bar{\tau} = \tau - (2-\delta)\eps m$, then by the above and optional stopping if follows that
$$ {\E \bar{\tau}^2 \over 2m} + \eps \E \bar{\tau} - \E Z_{\delta \eps m}
\leq  C \eps^3 m \, .$$
This and Lemma \ref{numberofrecmins} gives that \be \label{firstbd} \E
\bar{\tau} \leq {1 \over 2[(1-5\delta-5\delta \eps)\eps]} - {\E \bar{\tau}^2
\over 2\eps m} + O(\eps^2 m) \, .\ee

Next, we wish to derive a lower bound on $\E \bar{\tau}^2$. Put $T
= \delta m$, then for $t \leq T$ we have that
$$ \E \Big [ (M_t - M_{t-1})^2 \Big ] \geq 1 - \delta \, ,$$
hence the process
$$ M^2_{t \wedge T} - (1 - \delta) (t\wedge T) \, ,$$
is a submartingale and optional
stopping gives \be\label{opt1} (1-\delta) \E [\tau \wedge T] \leq \E
M^2_{\tau \wedge T} \, .\ee We now bound $\E M^2_{\tau
\wedge T}$ from above. We have
$$ W_{\tau \wedge T} = -Z_{\delta \eps m} {\bf 1} _{\{ \tau
\leq T \}} + W_T {\bf 1} _{\{ \tau > T \}} \, .$$ Thus,
$$ \E W^2_{\tau \wedge T} \leq \E Z^2_{\delta \eps m} + O(m^2)
\prob ( \tau > T) \, .$$ Since after time $\delta m/2$ the process is subcritical with constant negative drift we have that $\prob(\tau > T)$ decays exponentially in $m$. Lemma \ref{numrecmin} now yields that $\E
W^2_{\tau \wedge T} = \E Z_{\delta \eps m}^2 + o(1) = O(\eps^{-2})$. Next we estimate $\E f^2(\tau
\wedge T)$. Write $\mu = (2-\delta)\eps m$ and simplify $f(t)$ to
get
\begin{eqnarray*} f(t) &=& {(t-\mu)^2 (1-\eps) \over 2m} +
(t-\mu)\Big [ \eps - \delta \eps ^2 + {1+\eps \over 2m} - {\mu \eps
\over m} \Big ] -\delta \eps^2 \mu + {1+\eps \over 2m}\mu - {\eps
\over 2m} \mu^2
\\ &=& {(t-\mu)^2 (1-\eps) \over 2m} + (t-\mu)\Big [ \eps +
O(\eps^2) \Big ] + O(\eps^3 m) \, .\end{eqnarray*} Hence
\begin{eqnarray*} f^2(t) = {(t-\mu)^4 (1-\eps)^2 \over 4m^2} +
{(t-\mu)^3 (\eps + O(\eps^2)) \over m} + (t-\mu)^2 \eps^2(1+O(\eps))
+ (t -\mu) O(\eps^4 m) .
\end{eqnarray*}
Lemmas \ref{ytlowerbound} and \ref{cvbound} imply that $\E \bar{\tau}^k$ is of order
$(m/\eps)^{k/2}$ and hence the third term on the right hand side is
dominant so, \be \label{est1} \E f^2(\tau \wedge T) = (1+o(1))\eps^2 \E
\Big [ (\tau \wedge T - \mu)^2 \Big ]  \, .\ee We also use
Cauchy-Schwartz to estimate
$$ \Big | \E W_{\tau \wedge T} f(\tau \wedge T) \Big | \leq \sqrt{
\E W^2_{\tau \wedge T}} \sqrt {\E f^2(\tau \wedge T)} =
O(\sqrt{m/\eps}) = o(\E f^2(\tau \wedge T)) \, ,$$ since $\E
W^2_{\tau \wedge T} = O(\eps^{-2})$ and $\sqrt{m/\eps} = o(\eps m)$. We put this and (\ref{est1})
into (\ref{opt1}) and get that
$$ (1+o(1))\eps^2 \E \Big [ (\tau \wedge T - \mu)^2 \Big ] \geq
(1-\delta)\mu - (1-\delta)\E[\tau\wedge T - \mu] =
(1+o(1))(1-\delta)\mu \, ,$$ since $\E \bar{\tau} = O(\eps^{-1} m)$
and $\prob(\tau> T)$ decays exponentially in $m$. We learn that
$$ \E \bar{\tau}^2 \geq (1-o(1))(1-\delta)(2-\delta){ m \over \eps} \, . $$
Putting this into (\ref{firstbd}) gives that if $\delta>0$ is chosen
small enough (but fixed) and $m$ is large enough
$$ \E \bar{\tau} \leq - {1 \over 4\eps^2} \, ,$$
concluding the proof of the lemma. \qed\\

\noindent {\bf Proof of Theorem \ref{critgnp}.} Part (i) follows immediately from Lemma \ref{hitrecmin} since $|\C_{\delta \eps m}| \leq \delta \eps m + \tau_\delta$. To prove (ii) we proceed as in the proof of Lemma \ref{supergnpmoments} and write
$$ \E \sum_{\C_j \neq \C_{\delta \eps m}} |\C_j|^k = \E \sum_{\C_j \neq \C_{\delta \eps m}} |\C_j|^k{\bf 1}_{\{\C_{\delta \eps m} \leq 1.5\eps m \}} + \E \sum_{j \geq 2} |\C_j|^k {\bf 1}_{\{\C_{\delta \eps m} \geq 1.5\eps m \}} \, .$$
Lemma \ref{ytlowerbound} shows that $\prob(\C_{\delta \eps m} \leq 1.5\eps m) \leq Ce^{-c\eps^3 m}$ and so FKG inequality gives
$$ \E \sum_{\C_j \neq \C_{\delta \eps m}} |\C_j|^k{\bf 1}_{\{\C_{\delta \eps m} \leq 1.5\eps m \}} \leq Ce^{-c\eps^3 m} \E \sum_j |\C_j|^k = O(m\eps^{-2k+3}) \, ,$$
by part (i) of Lemma \ref{gnp2moments}. The second term is handled as in the proof of Lemma \ref{supergnpmoments} by conditioning on $\C_{\delta \eps m}$ and using Lemma \ref{gnp2moments} for the remaining subcritical graph.  \qed \\
